\documentclass{uofsthesis-cs}
\usepackage{amssymb}
\usepackage{amsmath}
\usepackage{amsthm} 
\usepackage{hyperref}
\usepackage{xcolor}
\hypersetup{
    colorlinks,
    linkcolor={black!50!black},
    citecolor={blue!50!black},
    urlcolor={blue!80!black}
}
\usepackage[sort&compress,noabbrev, nameinlink]{cleveref}
\usepackage[sort, numbers]{natbib}
\usepackage{url}

\newtheorem{thm}{Theorem}[section]
\newtheorem{lem}[thm]{Lemma}
\newtheorem{cor}[thm]{Corollary}

\newtheorem{defn}{Definition}[section]
\newtheorem{eg}[defn]{Example}
\newtheorem{rmk}[defn]{Remark}
\newtheorem{nott}[defn]{Notation}

\usepackage[normalem]{ulem}

\title{Generalized Metrics}

\author{Samer Assaf}

\degree{\PhD}         

\convocationdate{June/2016}

\department{Mathematics and Statistics}

\ptuaddress{Head of the Department of Mathematics and Statistics\\
Room 142 McLean Hall\\
106 Wiggins Road\\
University of Saskatchewan\\
Saskatoon, Saskatchewan\\
Canada\\
S7N 5E6}

\abstract{A distance on a set is a comparative function. The smaller the distance between two elements of that set, the closer, or more similar, those elements are. Fr\'echet axiomatized the notion of distance into what is today known as a metric. In this thesis we study several generalizations of Fr\'echet's axioms. These include partial metric, strong  partial metric,  partial $n-\mathfrak{M}$etric and strong partial $n-\mathfrak{M}$etric. Those generalizations allow for negative distances, non-zero distances between a point and itself and even the comparison of  $n-$tuples. We then present the scoring of a DNA sequence, a comparative function that is not a metric but can be modeled as a strong partial metric.
\\\indent Using the generalized metrics mentioned above we create topological spaces and investigate convergence, limits and continuity in them. As an application, we discuss contractiveness in the language of our generalized metrics and present Banach-like fixed, common fixed and coincidence point theorems.}

\acknowledgements{
The Ph.D. program was by far  my most challenging endeavor. Luckily, I did not go through it alone. I would like to take a moment to thank the silent heroes who helped me make this possible:
\\
\\\indent \textbf{Professor emeritus Edward Tymchatyn:} Prof. Tymchatyn was my supervisor during my Ph.D. program. He had an open door policy, always made sure to make time for me and  broadened my horizons by exposing me to topological conferences in which I  presented my findings. No question was too small for him to answer and, with an extensive knowledge in topology, he was an inexhaustible source of knowledge to be envied for. His attention to minute details, mathematical rigor and perfectionist approach provided lessons I hope to keep with me on my journey.
\\
\\\indent \textbf{Doctor Koushik Pal}: Dr. Pal is my co-author for the first two papers submitted. As an undergraduate student, I was used to being presented with a question and finding a way to prove it. At the graduate level however, things are not as straightforward. It is with him that I learnt that as mathematical investigators, we may actually change the question. This is a simple lesson but a very important one. It is with Dr. Pal that my journey as a researcher started.
\\
\\\indent \textbf{Associate Professor Christopher Dutchyn:} Prof. Dutchyn introduced me to the scoring of DNA strands mentioned in \Cref{C2b}. His valuable insight on DNA scoring scheme guided me as I investigated the subject. He also took time from his busy schedule to read and discuss my first results in partial metrics. Additionally,  he was my go-to reference in latex coding. The organized bibliography and referencing system are all thanks to him.
\\
\\\indent \textbf{Assistant Professor Katarzyna Kuhlmann:} Prof. Kuhlmann was the one who introduced me to the notion of generalized metrics. One day, while she was walking down the hall, I saw a partial metric paper in her hand and asked her about it. She took the time to introduce it, explain it and point me in the right direction for further research. She also took time to discuss dislocated metrics and present some examples. It is thanks to her positive attitude and willingness to help that ignited my interest in the topic which I ended up researching. 
\\
\\\indent \textbf{Professor Franz-Viktor Kuhlmann:} Prof. Kuhlmann organized the fixed point theorem seminars in our department. During those seminars, I was exposed to various techniques which sparked the main questions for my first paper. Additionally, he took the time to read the first paper and suggested valuable amendments.
\\
\\ \indent \textbf{Assistant Professor Walid Abou Salem:} There is more to the graduate program than the mathematics part and Prof. Abou Salem was a refreshing reference. A true teacher at heart, he understands the difficulties we experience as students and always provided valuable advice. He was a tremendous support on this journey.
\\
\\ \indent \textbf{Mathematics Department Office Administrators:} Simply stated, they are amazing. With an uplifting smile and the knowledge to navigate policies in the most efficient way, they made every problem, no matter how big, seem easily manageable. A detailed explanation of procedures, email reminders and ingenious solutions  were always happily provided in such a clear way that I had no problem staying on top of what was required of me.
\\
\\ \indent \textbf{Mathematics Department Systems Manager Richard Kondra:} From downloading a program I needed, to updating Bakoma packages, to adjusting screen problems, to setting lecture conferences and the list goes on and on. Mr. Kondra was our IT savior and thankfully, a very patient one.
\\
\\ \indent \textbf{Last and most certainly not least, my parents:} To my mother who spent countless hours nurturing the young mathematician  in me and my father who became an expatriate to support us, a deep and heartfelt thanks. Both are hugely responsible for the man I am today. I hope they know that I would not even dream of being here without their encouragement, patience and support. 
}

\dedication{This thesis is  dedicated to every graduate student who feels overwhelmed. It is normal to doubt your capabilities when you are dealing with professors who spent their entire career researching that field. Be patient: You do have what it takes, you just need experience. Persevere and enjoy the journey.
\vspace{10ex}
\\\indent \textit{``We fail to realize that mastery is not about perfection. It''s a process, a journey. The master stays on the path day after day, year after year. The master is willing to try and fail and try again, for as long as he/she lives.''}---George Leonard}

\begin{document}

\maketitle

\frontmatter

\chapter{Introduction}
\label{C1}
\setcounter{thm}{0}
\setcounter{defn}{0}
\section{Why do we generalize a metric?}
\label{C1a}
\indent \indent Given a set $X$, a metric is a function $d: X\times X \to
\mathbb{R}$ as defined in \Cref{C2a}. The point of a metric is for it to be a comparative function. For any two given elements $x$ and $y$ of  $X$, the value (also called the distance) $d(x,y)$ represents how close $x$ and $y$ are to each other. Here, we use ``close'' in a loose sense as its meaning may vary depending on what information  we need when comparing the elements of $X$. For example, when we say the distance between Paris and Saskatoon is 6944.32 km, the information we are looking for is the amount of space that needs  to be traveled to go from one to the other in a direct straight flight. In Computer Science, when we say the error between the actual value of $\pi$ and $3.14$ is of the order of $10^{-2}$, the information we are conveying is how good of an estimation is $3.14$ to the real value of $\pi$?
\\\indent Unfortunately some comparative functions, although very meaningful, do not adhere to the strict axioms of a metric presented in \Cref{C2a}. One example is the scoring presented in \Cref{C2c}. A scoring is a function used to compare the similarity between two DNA strands. We generalize the known metric in two steps depicted in the figure below:\begin{center}
Metric.
\\$\downdownarrows$
\\Partial metric and Strong partial metric.
\\ (\Cref{C2b} and \Cref{C2c})
\\$\downdownarrows$
\\Partial $n-\mathfrak{M}$etric and Strong Partial $n-\mathfrak{M}$etric.
\\(\Cref{C2e} and \Cref{C2f})
\\
\end{center}
First we generalize the metric into a partial metric by allowing the self distance (the distance between a point and itself) to take real values that may be different than zero. Then we generalize a partial metric that acts on pairs into a` partial $n-\mathfrak{M}$etric that acts on $n-$tuples. The Partial $n-\mathfrak{M}$etric is the most general case presented. As an application, we use these generalized metrics to form fixed point theorems. 
\section{Why do we study fixed point theory?}
\label{C1b}
\indent \indent Fixed point theory is a branch of topology that studies the conditions under which a function from a set to itself  has a fixed point.
\\\indent Advancements in fixed point theory enrich many scientific fields such as biology, chemistry, computer science, economics and game theory. Below we give one example from each field which illustrates the use  of  fixed point theory in that science.
\\\indent In biology, fixed point theory helps with studying  how cancer cells replicate. Statistical data is collected,  modeled, and  fixed point theorems are used to form an educated guess of how those cells will progress in the future.
\\\indent In chemistry, the branch that uses computer programs and simulations is called computational chemistry. This branch of chemistry investigates the charge density and electronic charge per volume of atoms, ions and molecules. The bigger the molecule studied, the more complex the problem becomes. Some of the problems are too complex to be solved analytically, hence, fixed point theory is used to develop an iterative process that will converge to a solution.
\\\indent In computer science, fixed point theory is used to check whether a program will stabilize or run indefinitely.
\\\indent In economics and game theory, fixed point theory is used  to prove that a certain simulation has at least one equilibrium point.
\section{How are fixed point theorems classified?}
\label{C1c}
\indent \indent There are two approaches to the study of fixed point theory. The first is an existential approach where, under some criteria, we ensure that a function has a fixed point. Unfortunately, the theory itself does not provide a way of constructing the fixed point. One example of such an approach is the Brouwer's Fixed Point Theorem which ensures the existence of a fixed point for a continuous function over a compact convex set in a Euclidian space. 
\\\indent The second approach is an iterative one where, under some criteria, we are able to construct a sequence whose terms converge to a limit which is a fixed point of the function in question. The main example is the Banach Fixed Point Theorem which asserts that  continuous functions on a complete metric space have a fixed point should those functions be contractive.
\section{What are we presenting in this thesis?}
\label{C1d}
\indent \indent In this thesis, we restrict ourselves to iterative approaches to fixed point theorems. Our major concern is that this method relies heavily on a metric to allow us to see if successive iterations are getting ``closer'' fast enough or not. In some cases, such as in comparing DNA strands, the use of a metric is too computationally taxing. One would like to use a comparison score which may not satisfy all the axioms of a metric. Generalized metrics are obtained by amending metric axioms to allow meaningful comparison scores to be used.\\
\\\indent We present four improvements to the metrics and fixed point theorems at hand:
\\
\\1) A metric assigns a number to each pair of points; the smaller the number, the closer the points are together. G\"ahler \cite{Gah1963,Gah1966}  presented the following question: Can we do a similar thing to assess a triplet? Will our new ``metric'' be stable enough to describe a $T_o$ topology? Many mathematicians have attempted to define such a ``metric''; some with more success than others. In this thesis we define and generalize n-metrics by assigning a number to each n-tuple and generate a $T_o$ topology  from that generalized metric.
\\
\\2)  Edelstein \cite{Ede1961,Ede1962,Ede1964} generalized the Banach Fixed Point Theorem to allow the continuous function of a metric space to itself  to be contractive only on the orbit of a certain point $x_o \in X$. We adapt Edelstein's idea and use it with  the new generalized metrics discussed above.
Our functions  may not even be continuous but are contractive on an orbit.\\
\\3) Matthews \cite{Mat1994,Mat2009} extended Banach's Fixed Point Theorem to "metrics" where self distances are not necessarily zero. Those ``metrics" are also known as partial metrics. His theorems, however, only applied where the limit has a self distance ( in this paper called a central distance) of zero. O'Neill \cite{One1996} extended the definitions of partial metrics to allow negative values. We extend Matthews' theorems to partial metrics in the sense of O'Neill. We also extend those theorems to allow the fixed point to have a self distance ( the distance between that fixed point and itself) to be any real number rather than be restricted to the value zero.  \\
\\4) Markov \cite{Mar1936} began investigating common fixed point theorems for commuting maps. His Theorem was later generalized by  Kakutani \cite{Kak1938} to give us the Markov-–Kakutani fixed-point theorem. Eilenberg and  Montgomery \cite{Eil1946} were among the first to investigate coincidence points for non-self maps, where the domain set is ordered.  We generalize the techniques used in fixed point theorems to common and coincidence point theorems. Our approach allows us to present common and coincidence point theorems where neither  commuting functions nor ordered spaces are required.
\section{How do we organize our work?}
\label{C1e}
\indent \indent We organize our work into six chapters:
\\\indent In  \Cref{C2}, we concentrate on generalizing the original notion of metric and partial metric. We define strong partial metrics, partial $n-\mathfrak{M}$etrics, and strong partial $n-\mathfrak{M}$etrics. We also present  other generalized metrics already found in the literature. We give examples of each and derive some  inequalities needed for the chapters to come.
\\\indent In  \Cref{C3}, we define open balls for each of the generalized metrics defined in \Cref{C2}. We show that these balls form a basis for a topology that is at least   $T_o$ and, in some generalized metric cases, even $T_1$. The $T_0$ property is important because we want our topologies to be able to distinguish points.
\\\indent In \Cref{C4}, we define  Cauchy sequences for each of our generalized metrics. We  investigate the properties of  the limits of  Cauchy sequences in each case.
We then move on to define a stronger version of the topological limit called the special limit, study its properties and define a complete space. We also introduce the notion of  Cauchy pairs. Those are pairs of sequences whose terms eventually get arbitrarily ``close" to one another.
\\\indent In \Cref{C5}, we present a variety of  contractive conditions  on   functions, or pairs of functions, sufficient for the  use of  iterative methods to obtain  Cauchy sequences.
\\\indent In \Cref{C6}, we investigate some criteria of  functions, such as a relaxed version of continuity,   and study their effect on the limits of the sequences found in  \Cref{C5}. Those properties will   make sure the iterations lead to a fixed point, common fixed point, or  a coincidence point as shown in  \Cref{C7}.
\\\indent \Cref{C7} is dedicated to  our main theorems. The proofs will be short and direct due to the extensive scaffolding created in the previous chapters. Our main theorems consist of  fixed point, common fixed point and coincidence point theorems in each of the generalized metric spaces. The fixed point theorems use a technique similar to Edelstein's by taking a function that is contracting on an orbit and provide additional constraints needed for the function to have a fixed point. The common fixed point theorems have a similar contractive approach.  Fisher \cite{Fis1979},  Yeh \cite{Yeh1979} and many others developed common fixed point theorems that rely on the two functions commuting ($f(g(x))=g(f(x))$), or a weaker form. We adapt our contractive approach  to generate common fixed point theorems that do not require the two functions to commute. Finally, we formulate a coincidence point theorem where the only requirement on the domain is for it to be a complete generalized metric space.
\section{Notations.}
\label{C1f}
Let $x_1,x_2,...,x_n$ and $a$ be elements of a set $X$.
\begin{center}
\begin{tabular}{c c}
$\mathbb{N}$:  & Set of natural numbers. \\
$\mathbb{R}$: & Set of real numbers. \\
$\mathbb{R}^{\ge 0}$: & Set of non-negative real numbers.\\
$\mathcal{P}(X)$: & Power set (set of all subsets)  of $X$. \\
$\mathbb{R}^{>0}$: & Set of positive real numbers. \\
$\langle a \rangle^n$: & The $n-$tuple $(a,a,a,...,a)$.\\
$ \langle x_i \rangle_{i=1}^n$: & The $n-$tuple $(x_1,x_2,x_3,...,x_n)$. \\
$(\langle x_i \rangle_{i=1}^k,\langle y_j \rangle_{j=1}^{n-k})$: & The $n-$tuple $(x_1,x_2,....,x_k,y_1,y_2,...,y_{n-k}).$ \\
\end{tabular}
\end{center}
It is very important not to confuse  
$$ \langle x_i \rangle_{i=1}^n =(x_1,x_2,x_3,...,x_n)$$
 with
$$\langle x_i \rangle^n=(x_i,x_i,...,x_i,x_i).$$
\chapter{Metrics and their generalizations}\label{C2}
\setcounter{thm}{0}
\setcounter{defn}{0}
Distances have been used in the world since times immemorial. Fr\'echet \cite{Fre1906} was the first to axiomatize a distance, he called it \'ecart. Hausdorff \cite{Hau1914} coined the term metric space and used metrics to define a topology.
M. Deza and E. Deza in their book \emph{ Encyclopedia of Distances} \cite{Enc2009} provided us with an unparalleled reference for a variety  of metrics that appear in the literature.
\\\indent In this chapter we present the various generalized metrics used in this thesis. Some of these generalizations are already found in the literature. We will state the motivation for each  and provide some examples. 
 In order not to make this chapter too cumbersome, we give proofs  only in non-trivial cases.
\section{Metric}
\label{C2a}
\begin{defn}\label{D2a1}  A \uline{\textbf{metric}} $d$
 on  a set $X$ is a function $d: X\times X \to\mathbb{R}$ satisfying the following axioms:
\\For all $x,y,z\in X,$
\\(m-lbnd):  $0\le d(x,y)$.\hspace{30ex} (non-negative axiom or lower boundary axiom)
\\(m-sym):  $d(x, y) = d(y, x)$.\hspace{24ex} (symmetry axiom)
\\(m-sep):  $d(x,y)=0 \iff x=y$.\hspace{18.5ex}(separation axiom)
\\(m-inq): $d(x, y)\le d(x, z) + d(z, y)$.\hspace{15ex} (triangular inequality axiom)
\end{defn}
\hspace{2ex} We remark that  in his paper Fr\'echet \cite{Fre1906}  assumed but did not actually state  (m-sym). The reader who would like more information on metrics and metric topologies may consult \cite{Mun2000}. 
\section{Partial Metric}
\label{C2b}
\indent\indent In \textbf{1992},  Matthews \cite{Mat1994} considered  finite sequences as  partially computed versions of  infinite sequences. He noticed that depending on the computational approach, from the computer's  perspective, two  computed sequences of the same original infinite sequence need not be  identical. Therefore, a distance between a sequence and itself need not be zero. Motivated by this observation, Matthews  generalized metrics into what he called  partial metrics where a distance between a point and itself need no longer be zero. This change led him to modify the initial metric axioms (see \Cref{D2a1}). His partial metrics however, allowed only non-negative values.
\\\indent In \textbf{1996},  O'Neill \cite{One1996} generalized Matthews' partial metric to allow negative values. When we mention a partial metric from here on, we will be referring to O'Neill's partial metric.
\begin{defn}\label{D2b1} A \textit{\uline{\textbf{partial metric}}} $p$ on  a set $X$ is a function  $p: X\times X \to\mathbb{R}$ satisfying the following axioms:
\\For all $ x, y, z\in X,$
\\(p-lbnd):  $p(x,x)\le p(x,y)$.
\\(p-sym):  $p(x, y) = p(y, x)$.
\\(p-sep):  $p(x,x)=p(x,y)=p(y,y) \iff x=y$.\\(p-inq): $p(x, y)\le p(x, z) + p(z, y)-p(z,z)$.
\end{defn}
It is easy to see that every metric is a partial metric with  self distance ( i.e. distance of a point to itself) equal to  zero. Every partial metric defines a metric as follows:
\begin{lem}\label{L2b1} (O'Neill \cite{One1996}: Lemma 2.7 ) Let $p$ be a partial metric defined on a set $X$. Define $d:X \times X \to \mathbb{R}$ as
follows: 
\\For all $x,y \in X$, $$d(x,y)=2p(x,y)-p(x,x)-p(y,y).$$ Then, $d$ is a metric on $X$.
\end{lem}
\textbf{Proof:} This Lemma is a special case of \Cref{T2e2}.
\\\indent We will show in \Cref{C3} that the   topological structure on a set induced  by a partial metric is coarser than that induced by the metric described in \Cref{L2b1}.
\\\indent We give some  examples of simple partial metrics.
\begin{eg}{\textbf{(Basic Partial Metric):}}
\label{E2b1}
\\Consider the set $X=\{x,y\}$. Let $p:X \times X\to \mathbb{R}$ be defined by:
$$p(x,x)=0,p(x,y)=1\text{ and }p(y,y)=1.$$
Then $p$ is a partial metric on $X$.
\end{eg}
\begin{eg}{\textbf{(Maximum Partial Metric \cite{Mat1994}):}}
\label{E2b2}
\\Consider $X \subseteq \mathbb{R}$. Let $p:X \times X\to \mathbb{R}$ be defined by setting for all $x,y\in X$,
$$p(x,y)=\max\{x,y\}.$$
Then $p$ is a partial metric on $X$.
\end{eg}
\textbf{Proof:}Let $ x,y,z \in X$.
\\\indent  Without loss of generality, we may assume that $x \le y$. Thus,  $$(\star):  p(x,x)=x \text{,  } p(y,y)=y \text{, and } p(x,y) =y .$$
\indent Proof of (p-lbnd): From $(\star)$ we deduce that $$p(x,x)=x \le y= p(x,y) \text{ and } p(y,y)=y =p(x,y) . $$
\indent Proof of (p-sym): From $(\star)$ we have $p(x,y)=y$ and $ p(y,x) = y$ , so, $p(x,y)= p(y,x) .$
\\
\\\indent Proof of (p-sep):
\\\indent $(\Leftarrow )$ Trivial.
\\\indent $(\Rightarrow)$: Suppose $$p(x,x)=p(x,y)=p(y,y).$$ Then from $(\star)$ it is trivial to see that $x=y$.
\\
\\\indent Proof of (p-inq): Knowing that $x\le y$, to prove (p-inq) we need to consider three cases corresponding to the three possible orderings of $\{x,y,z\}$:
\\\uline{Case 1}: Suppose  that $z\le x\le y$. Then $ p(x,z)+p(y,z)-p(z,z)=x+y-z =y+(x-z)\ge y=p(x,y)$.
\\
\\\uline{Case 2}: Suppose that  $x \le z \le y$. Then $ p(x,z)+p(y,z)-p(z,z)=z+y-z = y=p(x,y)$.
\\
\\\uline{Case 3}: Suppose that $x \le y \le z$. Then $ p(x,z)+p(y,z)-p(z,z)=z+z-z=z \ge y=p(x,y)$.
\\
\\Therefore, in all cases, $$p(x,y) \le p(x,z) + p(z,y)-p(z,z).\hspace{6ex}\square$$
\begin{eg}{\textbf{(Augmented Real Line):}}
\label{E2b3}
\\Consider the set $X=\mathbb{R}\cup\{a\}$ where $ a \notin \mathbb{R}$. Let $p:X \times X\to \mathbb{R}$ be defined by:
\\For all $x,y\in \mathbb{R}$,
$$p(a,a)=0,p(a,x)=\vert x \vert \text{ and }p(x,y)=\vert x-y \vert -1.$$
Then $p$ is a partial metric on $X$.
\end{eg}
\textbf{Proof:} Let $ x,y,z \in \mathbb{R}$.
\\\indent Proof of (p-lbnd): Three cases arise.\\\uline{Case 1:} $p(a,a)=0$ and $p(a,x)=\vert x \vert$, hence, $p(a,a)\le p(a,x)$.
\\\uline{Case 2:} $p(x,x)=-1$ and $p(a,x)=\vert x \vert$, hence, $p(x,x) \le p(a,x)$.
\\\uline{Case 3:} $p(x,x)=-1$ and $p(x,y)=\vert x-y\vert -1$, hence, $p(x,x) \le p(x,y)$.
\\
\\\indent Proof of (p-sym): Trivial.
\\\indent Proof of (p-sep):
\\\indent $( \Leftarrow )$ Trivial.
\\\indent $(\Rightarrow)$: Let $x,y \in \mathbb{R}$.
\\ Suppose that  $$p(x,x)=p(x,y)=p(y,y). $$ Then, $-1=\vert x-y\vert-1$ giving us that $x=y$.
\\Now suppose that $$p(a,a)=p(a,x)=p(x,x).$$ Then, $p(a,a)=p(x,x)$ which is a contradiction since $0\ne -1 .$
\\
\\\indent Proof of (p-inq): The proof of (p-inq) is given by considering four possible claims. Let $x,y,z \in X$. 
\\\uline{Claim 1}: $ p(a,a)\le p(a,x)+p(a,x)-p(x,x)$.
\\Proof of claim 1: $$p(a,x)+p(a,x)-p(x,x)=\vert x \vert +\vert x \vert -(-1)=2\vert x \vert +1\ge0=p(a,a) . $$
\uline{Claim 2}: $p(x,x)\le p(a,x)+p(a,x)-p(a,a)$.
\\Proof of claim 2: $$p(a,x)+p(a,x)-p(a,a)=2\vert x \vert \ge -1=p(x,x).$$
\uline{Claim 3}: $p(x,y)\le p(a,x)+p(a,y)-p(a,a)$.
\\Proof of claim 3: $$ p(a,x)+p(a,y)-p(a,a)=\vert x \vert +\vert y \vert \ge\vert x-y \vert \ge \vert x-y \vert -1=p(x,y).$$
\uline{Claim 4}: $p(x,y)\le p(x,z)+p(z,y)-p(z,z)$.
\\Proof of claim 4: We have
$$p(x,z)+p(z,y)-p(z,z)=\vert x-z \vert -1 +\vert y-z \vert -1 -(-1)$$
$$=\vert x-z\vert +\vert y-z\vert -1\ge \vert (x-z)-(y-z)\vert -1 \text{ \hspace{4ex}(by the triangular inequality for the absolute value in } \mathbb{R} \text{)} $$ $ = \vert x-y \vert -1=p(x,y).\hspace{6ex}\square$
\section{Strong Partial Metric}
\label{C2c}
\indent\indent DNA strands, proteins and words are all examples of finite sequences generated from a finite alphabet $\mathcal{A\mathfrak{}}$. A generic question is: Given two finite sequences $x=\langle x'_i\rangle_{i=1}^s$ and $y=\langle y'_j\rangle_{j=1}^t$, how similar are these two sequences?
\\\indent In the case of DNA for example, the alphabet $\mathcal{A}=\{C, G, A, T\}$. While studying mutation from a sequence $x=\langle x'_i\rangle_{i=1}^s$ to a sequence  $y=\langle y'_j\rangle_{j=1}^t$, it becomes important to come up with a measure that can effectively compare partial DNA strands. One such  measure is the following commonly used {\em scoring scheme} \cite{Bio2008} :
\\\indent We first align two given words by inserting gaps so that their lengths match. Formally we adjoin the symbol $-$, called a gap to form a new alphabet $\mathcal{A^\star}=\mathcal{A}\cup\{-\}$ where $- \notin\mathcal{A}$. Obviously more than one alignment of words is possible: We consider one such alignment $\mathcal{L}$ where we represent  $x$ in the alignment by $\langle x_i\rangle_{i=1}^n$ and $y$ by $\langle y_i\rangle_{i=1}^n$.  Fixing $\alpha ',\beta ', \gamma '\in\mathbb{R}$  we then do a letter-by-letter  comparison while assigning a score to each of  four distinct possibilities.
 Namely, a score of zero is assigned if both letters are $-$. The score $\gamma '$ is assigned if only one of the letters is in $\mathcal{A}$. The score  $\alpha '$ is assigned if the two letters match (are the same) and are in $\mathcal{A}$. The score  $\beta '$ is assigned  if the two letters mismatch (are distinct) and are both in $\mathcal{A}$. Then these scores are summed up to assign a total score  $h'_\mathcal{L}(x,y)$. Finally 
$$s'(x, y)=\max\{h'_\mathcal{L}(x,y)\vert \mathcal{L} \text{ is an alignment}\}$$
denotes the highest possible score of all alignments of $x$ and $y$. Then $s'(x,y)$ is used as a measure of similarity or dissimilarity of the two words. As an example, assume 
$$\alpha '= +1, \beta ' = -1 \text{  and } \gamma '= -2.$$
Therefore, the total score of the pair (CGATC, CAGA) for the particular alignment
\begin{center}\begin{tabular}{|c|c|c|c|c|c|c|}\hline
$-$ & $C$ & $G$ & $A$ & $-$ & $T$ & $C$ \\\hline
$-$ & $C$ & $-$ & $A$ & $G$ & $A$ & $-$ \\\hline
$0$ & $+1$ & $-2$ & $+1$ & $-2$ & $-1$ & $-2$ \\\hline
\end{tabular}
\end{center}
is $$0+1 - 2 + 1 -2 -1 -2 = -5.$$ It is not hard to show that the best possible score for this pair of words is $-2$ arising from the alignment (CGATC and C---AGA).
\\\indent In our conventional view of metrics, we need a scoring function $s$ in which closeness is indicated by smaller numbers in the usual ordering of $\mathbb{R}$. For that reason in \cite{Ass20151}, we investigate the function $s(x,y)=-s'(x,y)$ and from it derive the axioms for a strong partial metric.
 And hence in the example above, the score given by the best alignment becomes $+2$: 
 \begin{center}\begin{tabular}{|c|c|c|c|c|}\hline
$C$ & $G$ & $A$ & $T$ & $C$  \\\hline
$C$ & $-$ & $A$ & $G$ & $A$ \\\hline
$-1$ & $+2$ & $-1$ & $+1$ & $+1$ \\\hline
\end{tabular}
\end{center}
\begin{rmk}
\label{R2c1} As mentioned above, the difference between the scoring we are investigating and the one presented in \cite{Bio2008} is that we are taking the negative of their score i.e.
$$s(x,y)=-s'(x,y).$$
\end{rmk}
\newpage
\begin{defn}\label{D2c1} A \textbf{\uline{strong partial metric}} $s$ on a set $X$ is a function  $s: X\times X \to\mathbb{R}$  satisfying the following axioms:
\\For all $ x,y,z\in X,$
\\(s-lbnd):  $s(x,x)<s(x,y)$, for $ x\ne y$.
\\(s-sym):  $s(x, y) = s(y, x)$.
\\(s-inq): $s(x, y)\le s(x, z) + s(z, y)-s(z,z)$.
\end{defn}
\hspace{3ex}Notice that a separation axiom (s-sep)  is hidden in  (s-lbnd)  as:
$$s(x,y)\le s(x,x)\iff x=y.$$
Clearly a strong partial metric is a partial metric. Therefore, a  strong partial metric $s$ defines a metric $d$ as defined by \Cref{L2b1}.
\begin{eg}{(\textbf{Shifted Metric):}}
\label{E2c1}
\\Let $d:X \times X \to \mathbb{R}$ be a metric defined on a set $X$. For a real number  $r$, let $s_r:X \times X \to \mathbb{R}$ be defined as
$$s_r(x,y)=d(x,y)+r.$$
Then $s_r$ is a strong partial metric on $X$.
\end{eg}
\textbf{Proof:} Let $ x,y,z\in X$.
\\Proof of (s-lbnd): From (m-lbnd) and (m-sep) we get if $x \ne y$ then $d(x,x)=0<d(x,y)$. Hence, 
$$d(x,x)+r<d(x,y)+r$$
 giving us that 
$$s_r(x,x)<s_r(x,y).$$
\indent Proof of (s-sym): From (m-sym) we get
$$s_r(x,y)=d(x,y)+r=d(y,x)+r=s_r(y,x).$$
\indent Proof of (s-inq): From (m-inq) we know that
$$d(x,y)\le d(x,z)+d(z,y).$$ Adding $r$ to both sides we get that
$$d(x,y)+r \le d(x,z)+d(z,y)+r.$$ From (m-sep), we know that $d(x,x) = 0$. With the above we deduce that
\\ $$d(x,y) + r \le d(x,z)+r + d(y,z) +r+ d(x,x)-r.$$
\\Finally we use the definition of $s_r$ to obtain
$$s_r(x,y) \le s_r(x,z)+s_r(z,y)-s_r(z,z).\hspace{6ex}\square$$
\begin{eg}{\textbf{(Positive Real Line):}}
\label{E2c2}
\\Consider $X$ to be the set of all positive real numbers. Let $s:X \times X\to \mathbb{R}$ be the function defined by setting :
$$s(x,x)=x \text{ and } s(x,y)=x+y \text{ for }x\ne y.$$
Then $s$ is a strong partial metric on $X$.
\end{eg}
\textbf{Proof:} Let $x,y,z \in X$.
\\\indent The proof of (s-lbnd) and (s-sym) are quite straight forward.
\\\indent Proof of (s-inq):  There are four cases.
\\
\\\uline{Case} 1: If $x=y=z$ then $$s(x,y)=x=x+x-x=s(x,z)+s(y,z)-s(z,z).$$
\uline{Case 2}: If $x=y$ and $y \ne z$ then $$s(x,y)=x \le x + x + z= x+z + y + z -z=s(x,z) + s(y,z)-s(z,z).$$
\uline{Case 3}: If $x \ne y$ and $x=z$ then $$s(x,y)=x+y =x + y+z -z= s(x,z)+s(y,z)-s(z,z).$$
\uline{Case 4}: If $x \ne y, y \ne z,$ and $ z \ne x$ then $$s(x,y)=x+y\le x+z+y=x+z +y+z-z \le s(x,z)+s(y,z)-s(z,z).\hspace{6ex}\square$$
\indent As in Bio-informatics \cite{Bio2008}, we try to find out how similar (similarity as measured by the partial metric $s$) are the two finite sequences $x=\langle x'_i\rangle _{i=1}^s$ and  $y=\langle y'_j\rangle_{j=1}^t$  .
For the convenience of the reader, we restate the definition of our scoring function $s$.
\begin{defn}\label{D2c2}
Consider $X$ to be the set of finite sequences generated by a finite alphabet $\mathcal{A}$.
First we augment the alphabet $\mathcal{A}$ by adding  a gap element $-$ i.e. $\mathcal{A^\star}=\mathcal{A}\cup\{-\}$ where $-\notin\mathcal{A}$.
 Then for  $\alpha, \beta, \gamma \in \mathbb{R}$, we define a \textbf{\uline{scoring function}} $s: X \times X \to \mathbb{R}$ by following the blueprint presented in the introduction of \Cref{C2c}. For each $x,y \in X$, we use ``deletions" and ``insertions" (called \textit{InDels} or gaps) to align the two sequences so that their lengths match i.e. we represent $x$ by $\langle x_i\rangle _{i=1}^n$ and $y$ by $\langle y_i\rangle _{i=1}^n$.
We denote this alignment $\mathcal{L}$. We then compare the  $i^{th}$ terms and assign a score $h_{(\mathcal{L},i)}(x,y)$ in the manner below:
\\
\\
\begin{center}
\begin{tabular}{|c|c|c|}\hline
Letter-by-letter comparison& Terminology &  $h_{(\mathcal{L},i)}(x,y)$ \\\hline
$x_i \in \mathcal{A}$ and $y_i=-$ & Deletion & $\gamma$ \\\hline
$x_i = -$ and $y_i \in \mathcal{A}$  & Insert & $\gamma$  \\\hline
$x_i=y_i \in \mathcal{A}$ & Match &  $\alpha$ \\\hline
$x_i, y_i \in \mathcal{A}$ but $x_i \neq y_i$  & Mismatch & $\beta$ \\\hline
$x_i=y_i=-$  & Relay & $0$ \\\hline
\end{tabular}
\end{center}
$$h_\mathcal{L}(x,y)=s(x,y).$$
\indent A Relay will have no effect in comparing these two sequences. That is why it is given a score of zero. The best score $s(x,y)$ is attained when we get the alignment with the smallest possible score i.e. 
$$s(x,y) =\min \{ h_\mathcal{L}(x,y)\vert\mathcal{L} \text{ is an alignment }\}.$$
We call $s(x,y)$  the \textbf{\uline{score}} of the pair $(x,y)$.
\end{defn}
\indent \indent We should note that more than one possible alignment may give us the best score, but the score itself is unique.
\begin{rmk}
\label{R2c2} In \cite{Bio2008}, the ``best score" $s'(x,y)$ is given by taking the score of the alignment $\mathcal{L}$ that gives us the biggest value of $h_\mathcal{L}'(x,y)$. Thus in \cite{Bio2008}, $s'(x,y)=-s(x,y)$ and the relative letter-by-letter scores are taken as $\alpha'=-\alpha , \beta'=-\beta$ and $\gamma'=-\gamma$.
\end{rmk}
\begin{nott}\label{N2c1} From this point forward, we will denote $h_{(\mathcal{L},i)}(x,y)$ and $h_\mathcal{L}(x,y)$ by $h_i(x,y)$ and $h(x,y)$ respectively when it is clear to which alignment $\mathcal{L}$ we are referring.
\end{nott}
\indent \indent To ensure $s$ is a strong partial metric we need to require that a Match is strictly our most  favorable occurrence. A Mismatch is at least as good as  two InDels and a Relay is strictly better than an  InDel.
\begin{lem}\label{L2c1} Consider $X$ the set of finite sequences generated from the finite alphabet $\mathcal{A}$. Let $s:X\times X \to \mathbb{R}$ be the scoring function  in \Cref{D2c2}. If $\alpha < \min\{ \beta, \gamma ,0\}$, $\beta \le 2\gamma$ and  $\gamma$ is strictly positive then $s$ is a strong partial metric.
\end{lem}
\textbf{Proof:} Let $x,y,z \in X$. 
\\\indent Proof of (s-lbnd): We compare the sequences after optimally aligning them. I.e. $s(x,y) =h(x,y) $ where we represent $x$ and $y$ by  $\langle x_i\rangle_{i=1}^n$ and  $\langle y_i\rangle_{i=1}^n$ respectively. For all $i\in \{1,...,n\}$
\\
\begin{center}
\begin{tabular}{|c|c|c|c|}\hline
$(x_i,y_i)$ & $h_i(x,x)$ &$s_i(x,y)$ & Comparison\\\hline
$x_i=y_i \in \mathcal{A}$ & $\alpha$ & $\alpha$ & $\alpha \le \alpha$ \\\hline
$x_i\in \mathcal{A}$ and $y_i = -$ & $\alpha$ & $\gamma$ & $\alpha < \gamma$ \\\hline
$x_i=-$ and $y_i \in \mathcal{A}$ & $0$ & $\gamma$ & $0<\gamma$ \\\hline
$x_i, y_i \in \mathcal{A}$ but $x_i \neq y_i$ & $\alpha$ & $\beta$ & $\alpha < \beta$ \\\hline
\end{tabular}
\end{center}
\vspace{2ex}
Thus, for all $i \in \{1,.....,n\}$ $h_i(x,x) \le s_i(x,y).$ Now since $x \ne y$ then there exists an $i_o$ such that $s_i(x,y) \ne \alpha.$
\\Since a relay has a score of $0$ then for these particular representations of $x$ and $y$ we have $$s(x,x) = h(x,x)=\sum\limits_{i=1}^n h_i(x,x) < \sum\limits_{i=1}^n s_i(x,y)=s(x,y).$$ 
Hence, $s(x,x)<s(x,y)$ for $x \ne y.$
\\\indent Proof of (s-inq): To compare three sequences $x,y$ and $z$, we are going to optimally align $x$ with $z$ and $y$ with $z$ by adding the necessary InDels and Relays so that they all have the same length.
\\\indent Before giving the formal proof of (s-inq) we give a simple example illustrating this step. Consider $x=CGT$, $y=AGAGT$ and $z=CAGC$. Now for some scoring scheme, assume the optimal alignment of $x$ to $z$ and $y$ to $z$ is given below
\begin{center}\begin{tabular}{|c|c|c|c|c|c|}\hline
$x:$ & $C$ & $-$ & $G$ & $-$ & $T$ \\\hline
$z:$ & $C$ & $A$ & $G$ & $C$ & $-$ \\\hline
\end{tabular}
\end{center}
and
\begin{center}
\begin{tabular}{|c|c|c|c|c|c|c|}\hline
$y:$ & $-$ & $A$ & $G$ & $A$ & $G$ & $T$ \\\hline
$z:$ & $C$ & $A$ & $G$ & $-$ & $-$ & $C$ \\\hline
\end{tabular}
\end{center}
\vspace{2ex}
\indent \indent By using Relays we can amend the above alignments such that the representations of $x, y$ and $z$ have the same length. Also note that a relay has a score of zero, hence, $s(x,z)$ and $s(y,z)$ may be computed using the representations of $x$, $y$ and $z$ in the box below. Here a gap in either representation of $z$ above appears as a gap in the representation of $z$ below.
\\
\begin{center}
\begin{tabular}{|c|c|c|c|c|c|c|c|}\hline
$x:$ & $C$ & $-$ & $G$ & $-$ & $-$ & $-$ & $T$ \\\hline
$y:$ & $-$ & $A$ & $G$ & $A$ & $G$ & $T$ & $-$ \\\hline
$z:$ & $C$ & $A$ & $G$ & $-$ & $-$ & $C$ & $-$ \\\hline
\end{tabular}
\end{center}
\vspace{2ex}
\indent\indent As the example above demonstrates, we may consider $x, y$ and $z$ having the same representations length  $\langle x_i\rangle_{i=1}^n, \langle y_i\rangle_{i=1}^n  \text{and } \langle z_i\rangle_{i=1}^n,$ respectively and optimally aligned to $z$. Then $$s(z,z)=h(z,z)$$
$$s(x,z)=h(x,z)=\sum\limits_{i=1}^n h_i(x,z)$$ and $$s(y,z)=h(y,z)=\sum\limits_{i=1}^n h_i(x,y).$$
\\ This leaves us with $x$ and $y$ not necessarily optimally aligned but $s(x,y) \le h(x,y)$.
\\\indent We now move to prove that $$h(x,y)\le h(x,z)+h(y,z)-h(z,z)=s(x,z)+s(y,z)-s(z,z).$$
\newpage
We do this term by term, thus ten cases arise:
\vspace{2ex}
\begin{center}
\begin{tabular}{|c|c|}

\hline
Case 1 & $x_i=y_i ,z_i\in \mathcal{A}$ with $z_i$ distinct from $x_i$ \\\hline
Case 2 & $x_i=y_i\in \mathcal{A}$, and $z_i=-$ \\\hline
Case 3 & $x_i, y_i, z_i \in \mathcal{A}$ and all three are distinct \\\hline
Case 4 & $x_i=z_i,y_i \in \mathcal{A}$ with $y_i$ distinct from $x_i$ \\\hline
Case 5 & $x_i$ and $y_i$ are distinct elements of $\mathcal{A}$, and $z_i=-$ \\\hline
Case 6 & $x_i=-$, and $y_i = z_i \in \mathcal{A}$ \\\hline
Case 7 & $y_i$ and $z_i$ are distinct elements of $\mathcal{A}$, and $x_i=-$\\\hline
Case 8 & $x_i=z_i=-$, and $y_i \in \mathcal{A}$ \\\hline
Case 9 & $x_i=y_i=-$, and $z_i \in \mathcal{A}$ \\\hline
Case 10 & $x_i=y_i=z_i=-$ \\\hline
\end{tabular}
\end{center}
\vspace{4ex}
\indent \indent We remind the reader that $\alpha < \min\{ \beta, \gamma ,0\}$, $\beta \le 2\gamma$ and  $\gamma>0$
\vspace{4ex}
\\\uline{Case 1}: $h_i(x,z)+h_i(y,z)-h_i(z,z)=\beta+\beta-\alpha$ and $h_i(x,y)=\alpha$.
\\ Hence, $ h_i(x,y) \le h_i(x,z)+h_i(y,z)-h_i(z,z)$ since $\alpha < \beta$.
\\
\\\uline{Case 2}: $h_i(x,z)+h_i(y,z)-h_i(z,z)=\gamma+\gamma-0$ and $h_i(x,y)=\alpha$.
\\ Hence, $ h_i(x,y) \le h_i(x,z)+h_i(y,z)-h_i(z,z)$ since $\gamma > 0$ and, therefore, $\alpha <0 < 2\gamma$.
\\
\\\uline{Case 3}: $h_i(x,z)+h_i(y,z)-h_i(z,z)=\beta+\beta-\alpha$ and $h_i(x,y)=\beta$.
\\ Hence, $ h_i(x,y) \le h_i(x,z)+h_i(y,z)-h_i(z,z)$ since $\alpha < \beta$.
\\
\\\uline{Case 4}: $h_i(x,z)+h_i(y,z)-h_i(z,z)=\alpha+\beta-\alpha$ and $h_i(x,y)=\beta$.
\\ Hence, $ h_i(x,y) \le h_i(x,z)+h_i(y,z)-h_i(z,z)$. 
\\
\\\uline{Case 5}: $h_i(x,z)+h_i(y,z)-h_i(z,z)=\gamma+\gamma-0$ and $h_i(x,y)=\beta$.
\\ Hence, $ h_i(x,y) \le h_i(x,z)+h_i(y,z)-h_i(z,z)$ since $ \beta \le 2\gamma$.
\\
\\\uline{Case 6}: $h_i(x,z)+h_i(y,z)-h_i(z,z)=\gamma+\alpha-\alpha$ and $h_i(x,y)=\gamma$.
\\Hence, $ h_i(x,y) \le h_i(x,z)+h_i(y,z)-h_i(z,z)$. 
\\
\\\uline{Case 7}: $h_i(x,z)+h_i(y,z)-h_i(z,z)=\gamma+\beta-\alpha$ and $h_i(x,y)=\gamma$.
\\ Hence, $ h_i(x,y) \le h_i(x,z)+h_i(y,z)-h_i(z,z)$ since $\alpha < \beta$.
\\
\\\uline{Case 8}: $h_i(x,z)+h_i(y,z)-h_i(z,z)=0+\gamma-0$ and $h_i(x,y)=\gamma$.
\\Hence, $ h_i(x,y) \le h_i(x,z)+h_i(y,z)-h_i(z,z)$.
\\
\\\uline{Case 9}: $h_i(x,z)+h_i(y,z)-h_i(z,z)=\gamma+\gamma-\alpha$ and $h_i(x,y)=0$.
\\ Hence, $ h_i(x,y) \le h_i(x,z)+h_i(y,z)-h_i(z,z)$ since $\gamma > 0$ and, therefore, $\alpha <0 < 2\gamma$.
\\
\\\uline{Case 10}: $h_i(x,z)+h_i(y,z)-h_i(z,z)=0+0-0$ and $h_i(x,y)=0$.
\\ Hence, $ h_i(x,y) \le h_i(x,z)+h_i(y,z)-h_i(z,z)$ since $\alpha < \beta$.
\\
\\In fact, Case $10$ need not arise when comparing three sequences. It is a useful tool however when using multiple alignment schemes discussed in  \Cref{C2f}. $\hspace{6ex}\square$
\begin{rmk}\label{R2c3}\Cref{L2c1} remains valid even if $\alpha$ and $\beta$ are functions with $\max\{\alpha\} <\min\{\beta,\gamma,0\}$, $\max\{\beta\} \le 2\gamma$, $\gamma >0$, and \uline{Case 3 } of the triangular inequality holds for distinct $x_i, y_i,z_i\in \mathcal{A}$. 
\end{rmk}
\begin{eg}{\textbf{(BLOSUM$62$):}}
\label{E2c3}
\\One example of the scoring scheme in \Cref{L2c1}is BLOSUM$62$ \cite{Bio2008}.
\end{eg}
\indent \indent We also note that in the alignment schemes, available in the bio-informatics literature, the letter-by-letter scores are restricted to:
$$\alpha < 0, \beta \ge 0 \text{ and }\gamma \ge 0.$$
Hence, $\beta \le 2\gamma$ and \uline{Case 3} are the remaining  requirements to check for $s$ to be a strong partial metric in those schemes.
\section{$G-$metric and $n-\mathfrak{M}$etric}
\label{C2d}
\indent \indent Many mathematicians attempted to generalize  a distance, a functions that assigns values to pairs, to a  function that assigns values to triplets or even to $n-$tuples. In \textbf{1963}, G\"ahler \cite{Gah1963,Gah1966}    attempted a generalization to triplets by modeling his axioms to mimic the area of the triangle whose vertices are three given points. While he called the function a $2-$metric, we will refer to it as G\"ahler's $2-$metric to avoid any confusion that may arise from later definitions.\begin{defn}\label{D2d1} A G\"ahler's  \textbf{\uline{$2-$metric}}  $\sigma$ on a set $X$ is  a function $\sigma: X \times X \times X \to \mathbb{R}$ satisfying the following axioms:
\\For all $ x,y,z,a\in X,$
\\(2-lbnd):  $0=\sigma(x,x,y)\le \sigma(x,y,z)$.
\\(2-sym):  $\sigma(x,y,z) = \sigma(\Pi\{x,y,z\})$, \hspace{3ex} where $\Pi$ denotes a permutation on $\{x,y,z\}$.
\\(2-sep):  If $x\ne y$ then there is at least one element $z\in X$ such that $\sigma(x,y,z)\ne 0$.
\\(2-inq): $\sigma(x,y,z)\le \sigma(x,y,a) + \sigma(x,a,z)+\sigma(a,y,z)$.
\end{defn}
\indent \indent In \textbf{1988}, Ha,  Cho and White \cite{Cho1988}  showed that G\"ahler's $2-$metric is not a generalization of a distance by giving an example  of a metric space which is not a G\"ahler's $2-$metric space and and an example of a G\"ahler's $2-$metric space which is not a metric space. To remedy this, in \textbf{1992}, Dhage \cite{Dha1992}  defined a $D-$metric by modeling his axioms on the perimeter of the triangle whose vertices are three given points.
\begin{defn}\label{D2d2} A \textbf{\uline{$D-$metric}}  $D$ on a set $X$ is  a function $D: X \times X \times X \to \mathbb{R}$ satisfying the following axioms:
\\For all $ x,y,z,a\in X,$
\\(D-lbnd): $0 \le D(x,y,z)$.
\\(D-sym):  $D(x,y,z) = D(\Pi\{x,y,z\})$, where $\Pi$ denotes a permutation on $\{x,y,z\}$.
\\(D-sep):  $D(x,y,z)=0 \iff x=y=z$
\\(D-inq): $D(x,y,z)\le D(x,y,a) + D(x,a,z)+D(a,y,z)$.
\end{defn}
\indent \indent In \textbf{2004}, Mustafa and Sims \cite{Mus2003}  showed that most of Dhage's claims about the structure of a $D-$metric space were incorrect. In \textbf{2006}, Mustapha and Sims \cite{Mus2006} modified Dhage's axioms by changing (D-inq) and introducing an additional boundary axiom. They called their function a  $G$-metric.
\begin{defn}\label{D2d3}  A \textbf{\uline{$G-$metric}}  $G$ on a set $X$ is  a function $G: X \times X \times X \to \mathbb{R}$ satisfying the following axioms:
\\For all $ x,y,z,a\in X,$
\\(G-lbnd):  $0 < G(x,x,y) \le G(x,y,z)$, \hspace{3ex} for $x \ne y$ and $y \ne z$ ,
\\(G-sym):  $G(x,y,z) = G(\Pi\{x,y,z\})$, \hspace{3ex} where $\Pi$ denotes a permutation on $\{x,y,z\}$.
\\(G-sep):  $G(x,y,z)=0 \iff x=y=z$,
\\(G-inq): $G(x,y,z)\le G(x,y,a) + G(a,a,z)$.
\end{defn}
\indent \indent In \textbf{2012},  Khan \cite{Kha2012} extended the Mustafa-Sims $G-$metric above into what he called a $K-$metric which is  a function $K: X^n \to \mathbb{R}^{\ge 0}$ for $n \ge 2$. The notation used is presented in \Cref{C1f}.
\begin{defn}\label{D2d4'} For $n\ge 2$, a \uline{\textbf{$K-$metric}} $K$ on a set $X$ is a function $K: X^n \to \mathbb{R}$ satisfying the following axioms:
\\For all  $(\langle x_i\rangle_{i=1}^{n},a)\in X^{n+1}$,
\\(K-lbnd):  $0 < K(\langle x_1\rangle^{n-1},x_2)\le K(\langle x_i \rangle_{i=1}^{n})\text{ \hspace{2ex}for all distinct elements } x_1,x_2,...,x_n $.
\\(K-sym):  $K(\langle x_i\rangle_{i=1}^{n}) = K(\langle x_{\pi(i)}\rangle_{i=1}^{n})$, \hspace{2ex} where $\pi$ is a permutation on $\{1,....,n\}$.
\\(K-sep):  $K(\langle x_i\rangle_{i=1}^n)=0 \iff x_1=x_2=x_3=...=x_n$,
\\(K-inq): $K(\langle x_i\rangle_{i=1}^{n})\le K(\langle x_i\rangle_{i=1}^{n-1},a) + K(\langle a \rangle^{n-1},x_n)$.
\end{defn}
The $K-$metric is modeled on the perimeter of an $n-$simplex (i.e. the sum of the length of the sides of the $n$-simplex).   For $n=2$, the $K-$metric axioms are just the usual metric axioms. For $n=3$, Khan's definition coincides with Mustafa and Sim's definition of a $G-$metric. We found that the axioms Khan proposed were unnecessarily restrictive. That is why in \cite{Ass20152}, we proposed the $n-\mathfrak{M}$etric. The reader should note the use of the scripted $\mathfrak{M}$ to differentiate it from other generalizations found in the literature.
\begin{defn}\label{D2d4} For $n\ge 2$, an \uline{\textbf{$n-\mathfrak{M}$etric}} $M$ on a set $X$ is a function $M: X^n \to \mathbb{R}$ satisfying the following axioms:
\\For all  $(\langle x_i\rangle_{i=1}^{n},a)\in X^{n+1}$,
\\(n-lbnd):  $0 \le M(\langle x_1\rangle^{n-1},x_2)$.
\\(n-sym):  $M(\langle x_i\rangle_{i=1}^{n}) = M(\langle x_{\pi(i)}\rangle_{i=1}^{n})$, \hspace{2ex} where $\pi$ is a permutation on $\{1,....,n\}$.
\\(n-sep):  $M(\langle x_1\rangle^{n-1},x_2)=0 \iff x_1=x_2$,
\\(n-inq): $M(\langle x_i\rangle_{i=1}^{n})\le M(\langle x_i\rangle_{i=1}^{n-1},a) + M(\langle a \rangle^{n-1},x_n)$.
\end{defn}
\begin{rmk}\label{R2d1} The $2-\mathfrak{M}$etric, not to be confused with G\"ahler's $2-$metric, is simply a metric. As mentioned below a $K-$metric (see \Cref{D2d4'}) is a special case of an $n-\mathfrak{M}$etric. Hence, a  $G-$metric (see \Cref{D2d3}) is a special case of a $3-\mathfrak{M}$etric.  \end{rmk}
\indent \indent For $n\ge 3$, our $n-\mathfrak{M}$etric axioms relax Khan's $K-$metric axioms in three ways. First we allow negative values (see \Cref{E2d2} ).  The second difference lies in  (n-lbnd)  which is a major weakening of (K-lbnd).
Third, our (n-sep) is  weaker than (K-sep).
\\\indent Mustafa and Sims \cite{Mus2006} proved several properties of  $G-$metrics. Not all of those properties hold for our $3-\mathfrak{M}$etric due to our weakened (n-lbnd) condition mentioned above. One property, however, 
$$G(a,a,b) \le 2G(b,b,a)$$
 still holds for a $3-\mathfrak{M}$etric. We were able to generalize it into a tool for term replacement in the $n-\mathfrak{M}$etric case in \Cref{T2d1}.

\begin{thm}{\textbf{(Term Replacement):}}
\label{T2d1}
\\Let $M$ be an $n-\mathfrak{M}$etric on a set $X$. For all $ \langle x_i \rangle_{i=1}^n,\langle y_i \rangle_{i=1}^n\in X^{n}$ and for $t\in\{1,.......,n\}$, $$M(\langle x_i \rangle_{i=1}^{n}) \le M(\langle y_j \rangle_{j=1}^{t},\langle x_i\rangle_{i=t+1}^n)+\sum\limits_{j=1}^tM(\langle y_j \rangle^{n-1},x_j).$$
\end{thm}
\textbf{Proof:} Let $ \langle x_i \rangle_{i=1}^n,\langle y_i \rangle_{i=1}^n\in X^{n}$. For $t=1$, the result follows by (n-inq) and (n-sym). 
\\Let $t \in \{2,....,n-1\}$ and  assume that the inequality holds for $t-1$. Then $$M(\langle x_i \rangle_{i=1}^{n}) \le M(\langle y_j \rangle_{i=1}^{t-1},\langle x_i\rangle_{i=t}^n)+\sum\limits_{j=1}^{t-1}M(\langle y_j \rangle^{n-1},x_j)$$
by (n-sym) $$=M(\langle y_j \rangle_{i=1}^{t-1},\langle x_i\rangle_{i=t+1}^n,x_t)+\sum\limits_{j=1}^{t-1}M(\langle y_j \rangle^{n-1},x_j)$$
by (n-inq) $$\le M(\langle y_j \rangle_{i=1}^{t-1},\langle x_i\rangle_{i=t+1}^n,y_t)+M(\langle y_t\rangle^{n-1},x_t)+\sum\limits_{j=1}^{t-1}M(\langle y_j \rangle^{n-1},x_j)$$
by (n-sym) $$M(\langle y_j \rangle_{i=1}^{t},\langle x_i\rangle_{i=t+1}^n)+\sum\limits_{j=1}^tM(\langle y_j \rangle^{n-1},x_j). \hspace{6ex} \square$$
\begin{rmk}\label{R2d2} The theorem above gives rise to important tools used in \Cref{C3c}. We present them in the corollaries below.
\end{rmk} 
\begin{cor}\label{C2d1} Let $M$ be an $n-\mathfrak{M}$etric on a set $X$. For all $ \langle x_i \rangle_{i=1}^n,\langle y_i \rangle_{i=1}^n\in X^{n}$ $$M(\langle x_i \rangle_{i=1}^{n}) \le M(\langle y_i \rangle_{i=1}^{n})+\sum\limits_{j=1}^nM(\langle y_j \rangle^{n-1},x_j).$$
\end{cor}
\textbf{Proof:} This is the case of \Cref{T2d1} when $t=n. \hspace{6ex} \square$
\begin{cor}\label{C2d2} Let M be an $n-\mathfrak{M}$etric on a set $X$. Then for  $ a,b \in X$ for $t\in\{1,.......,n\}$, $$M(\langle a\rangle^{t},\langle b\rangle^{n-t}) \le tM(\langle b\rangle^{n-1},a).$$
\end{cor}
\textbf{Proof:} By \Cref{C2d1} where $\langle x_i\rangle_{i=1}^n=(\langle a\rangle_{1}^{t},\langle b\rangle_{t+1}^{n})$ and $\langle y_j\rangle_{j=1}^n=\langle b\rangle^n$ we get $$M(\langle a\rangle^{t},\langle b\rangle_{t+1}^{n})\le M(\langle b \rangle^n)+\sum\limits_{j=1}^nM(\langle y_j \rangle^{n-1},x_j)$$
by (n-sep)
$$=0+\sum\limits_{j=1}^tM(\langle y_j \rangle^{n-1},x_j)+\sum\limits_{j=t+1}^nM(\langle y_j \rangle^{n-1},x_j)$$
$$=\sum\limits_{j=1}^tM(\langle b \rangle^{n-1},a)+\sum\limits_{j=t+1}^nM(\langle b \rangle^{n})$$
by (n-sep) again$$=\sum\limits_{j=1}^tM(\langle b \rangle^{n-1},a)=tM(\langle b \rangle^{n-1},a). \hspace{6ex} \square$$
\begin{cor}\label{C2d3}Let M be an $n-\mathfrak{M}$etric on a set $X$. Then for  $ a,b \in X,$ $$M(\langle a\rangle^{n-1},b) \le (n-1)M(\langle b\rangle^{n-1},a).$$
\end{cor}
\textbf{Proof:} This is the case of \Cref{C2d2} when  $t=n-1.\hspace{6ex} \square$
\\
\\\indent We now show that each $n-\mathfrak{M}$etric on a set $X$ naturally induces a metric on $X$. 
\begin{thm}{\textbf{(Metric from an $n-\mathfrak{M}$etric):}}
\label{T2d2}
\\Let M be an $n-\mathfrak{M}$etric on a set $X$. For $x,y \in X$ let $$d(x,y)=M(y,\langle x \rangle^{n-1})+M(x,\langle y \rangle^{n-1}).$$ Then $d$ is a metric on the set $X$.
\end{thm}
\textbf{Proof:} Let $ x,y,z \in X$.
\\\indent Proof of (n-lbnd): From (n-lbnd), we know that $M(x,\langle y \rangle^{n-1})\ge 0$ and $M(y,\langle x \rangle^{n-1})\ge 0.$
\\Hence, $d(x,y)=M(x,\langle y\rangle^{n-1})+M(y,\langle x \rangle^{n-1}) \ge 0.$
\\\indent Proof of (n-sym): Symmetry of $d$ follows from the symmetry of addition of real numbers.
\\\indent Proof of (n-sep):
\\\indent ($\Rightarrow$) From the definition of $d$, if $d(x,y)=0$ then $$M(x,\langle y \rangle^{n-1})+M(y,\langle x \rangle^{n-1})=0.$$
By (n-lbnd) 
$$M(x,\langle y \rangle^{n-1})=0=M(y,\langle x \rangle^{n-1}).$$
By (K-sep), $x=y$.
\\
\\\indent ($\Leftarrow$) If $x=y$ then  $M(x,\langle y \rangle^{n-1})=M(y,\langle x \rangle^{n-1})=M(\langle x \rangle^{n})=0$ by (n-sep).
\\Hence, 
$$d(x,y)=M(y,\langle x \rangle_{1}^{n-1})+M(x,\langle y \rangle_{1}^{n-1})=0.$$
\indent Proof of (n-inq): From (n-sym) and (n-inq) we get
$$M(x,\langle y \rangle^{n-1})=M(\langle y \rangle^{n-1},x)$$
$$\le M(\langle y \rangle^{n-1},z)+M(\langle z \rangle^{n-1},x)$$
$$=M(x,\langle z \rangle^{n-1})+M(z,\langle y \rangle^{n-1}).$$
\\Similarly
$$M(y,\langle x \rangle^{n-1})\le M(y,\langle z \rangle^{n-1})+M(z,\langle x\rangle^{n-1}).$$ Hence,
$$ d(x,y)=M(x,\langle y \rangle^{n-1})+M(y,\langle x\rangle^{n-1})$$
$$\le M(x,\langle z \rangle^{n-1})+M(z,\langle y \rangle^{n-1})+M(z,\langle x \rangle^{n-1})+M(y,\langle z \rangle^{n-1})$$
$$=M(x,\langle z \rangle^{n-1})+M(z,\langle x \rangle^{n-1})+M(z,\langle y \rangle^{n-1})+M(y,\langle z \rangle^{n-1})$$
$$=d(x,z)+d(z,y). \hspace{6ex} \square$$
\indent We shall show in \Cref{C3c} that the topology on $X$ induced by the $n-\mathfrak{M}$etric coincides with the topology on $X$ induced by the metric presented in \Cref{T2d2}. 
\begin{eg}{\textbf{(Unit $n-\mathfrak{M}$etric):}}
\label{E2d1}
\\Consider $X=\mathbb{R}$. Let $M:X^n\to \mathbb{R}$ be defined as follows:\\For all $ \langle x_i\rangle _{i=1}^{n} \in \mathbb{R}^n,M(\langle x_i \rangle_{i=1}^{n})=
\begin{cases}0 &$ if $  x_1=x_2=...=x_n. \\
1 & $otherwise.$
\end{cases}$
\\Then $M$ is an $n-\mathfrak{M}$etric.
\end{eg}
\begin{eg}{\textbf{($5-\mathfrak{M}$etric with Negative values):}}
\label{E2d2}
\\Consider $X=\{a,b\}.$ Let $M:X^5 \to \mathbb{R}$ be defined as follows:
\\For all $ \langle x_i \rangle_{i=1}^5 \in X^5$, $M(\langle x_i \rangle_{i=1}^5) = M(\langle x_{\pi(i)} \rangle_{i=1}^5).\text{ (Where } \pi \text{ is a permutation on } \{1,2,3,4,5\} \text{)}$  
\\Furthermore,
\\$M(a,a,a,a,a)=0 ,M(b,b,b,b,b)=0,$
\\$M(a,a,a,a,b)=3$,  $M(a,b,b,b,b)=4,$
\\$M(a,a,a,b,b)=-1$, and $M(a,a,b,b,b)=2$.
\\Then M is a $5-\mathfrak{M}$etric on the set $X$.
\end{eg}
\textbf{Proof:} (n-lbnd), (n-sym), and (n-sep) are direct results from the definition of $M$.
\\\indent Proof of (n-inq): There are ten cases:
\\$M(a,a,a,a,a)\le M(a,a,a,a,b)+M(b,b,b,b,a)$, since $(0\le 3+4).$
\\$M(a,a,a,a,b)\le M(a,a,a,b,b)+M(b,b,b,b,a)$, since $(3\le -1+4).$
\\$M(a,a,a,a,b)\le M(a,a,a,a,a)+M(a,a,a,a,b)$, since $(3\le 0+3).$
\\$M(a,a,a,b,b)\le M(a,a,b,b,b)+M(b,b,b,b,a)$, since $(-1\le 2+4).$
\\$M(a,a,a,b,b)\le M(a,a,a,a,b)+M(a,a,a,a,b)$, since $(-1\le 3+3).$
\\$M(a,a,b,b,b)\le M(a,b,b,b,b)+M(b,b,b,b,a)$, since $(2\le 4+4).$
\\$M(a,a,b,b,b)\le M(a,a,a,b,b)+M(a,a,a,a,b)$, since $(2\le -1+3).$
\\$M(a,b,b,b,b)\le M(b,b,b,b,b)+M(b,b,b,b,a)$, since $(4\le 0+4).$
\\$M(a,b,b,b,b)\le M(a,a,b,b,b)+M(a,a,a,a,b)$, since $(4\le 2+3).$
\\$M(b,b,b,b,b)\le M(a,b,b,b,b)+M(a,a,a,a,b)$, since $(0\le 4+3).\hspace{6ex}\square$
\begin{thm}{(\textbf{$n-\mathfrak{M}$etric from a metric):}}
\label{E2d3}
\\Every metric $d$ on a set $X$ naturally defines an $n-$metric  $M$ on $X$ as follows:\\For all $\langle x_i \rangle_{i=1}^n\in X^n$, let $$M(\langle x_i \rangle_{i=1}^n)=\sum\limits_{t=2}^{n}\sum\limits_{i=1}^{t-1}d(x_i,x_t).$$
\end{thm}
\textbf{Proof:} Let $(\langle x_j \rangle_{j=1}^n,a)\in X^{n+1}.$
\\\indent Proof of (n-lbnd): Since $d$ is a metric on $X$,
 for $x_1 \ne x_2$$$M(\langle x_1 \rangle^{n-1},x_2)=\sum\limits_{t=2}^{n-1}\sum\limits_{i=1}^{t-1}d(x_1,x_1)+\sum\limits_{i=1}^{n-1}d(x_1,x_2)$$
$$=\sum\limits_{t=2}^{n-1}\sum\limits_{i=1}^{t-1}0+\sum\limits_{i=1}^{n-1}d(x_1,x_2)$$
$$=(n-1)d(x_1,x_2)> 0.$$
\\\indent Proof of (n-sym): Follows from (n-sym).
\\
\\\indent Proof of (n-sep):
\\\indent$(\Leftarrow)$ If $x_1=x_2$, then  $M(\langle x_1 \rangle^{n-1},x_2)=\sum\limits_{t=2}^{n}\sum\limits_{i=1}^{t-1}d(x_1,x_1)=0$.
\\\indent$(\Rightarrow)$ Now if $M(\langle x_1 \rangle^{n-1},x_2)=0$, then $$\sum\limits_{t=2}^{n-1}\sum\limits_{i=1}^{t-1}d(x_1,x_1)+\sum\limits_{i=1}^{n-1}d(x_1,x_2)=0,$$ and, hence, $$(n-1)d(x_1,x_2)=0$$ i.e.$$d(x_1,x_2)=0.$$ By (n-sep) we get that $x_1=x_2.$
\\
\\\indent Proof of (n-inq): Using the definition of $M(\langle x_i \rangle_{i=1}^n)$ and (n-inq) on $\sum\limits_{i=1}^{n-1}d(x_i,x_n)$ we get
$$M(\langle x_i \rangle_{i=1}^n)=\sum\limits_{t=2}^{n}\sum\limits_{i=1}^{t-1}d(x_i,x_t)$$
$$=\sum\limits_{t=2}^{n-1}\sum\limits_{i=1}^{t-1}d(x_{i},x_{t})+\sum\limits_{i=1}^{n-1}d(x_i,x_n)$$
and, hence, by (d-inq)$$\le \sum\limits_{t=2}^{n-1}\sum\limits_{i=1}^{t-1}d(x_{i},x_{t})+\sum\limits_{i=1}^{n-1}(d(x_i,a)+d(a,x_n))$$
$$= \sum\limits_{t=2}^{n-1}\sum\limits_{i=1}^{t-1}d(x_{i},x_{t})+\sum\limits_{i=1}^{n-1}d(x_i,a)+ \sum\limits_{t=2}^{n-1}\sum\limits_{i=1}^{t-1}d(a,a)+\sum\limits_{i=1}^{n-1}d(a,x_n)$$
$$=M(\langle x_i \rangle_{i=1}^{n-1},a)+M(\langle a \rangle^{n-1},x_{n}). \hspace{6ex} \square$$
\section{$G_p-$Metric and Partial $n-\mathfrak{M}$etric}
\label{C2e}
\indent \indent In\textbf{ 2011},  Zand and Nezhad \cite{Zan2011} defined a function which they called a $G_p-$metric. A $G_p-$ metric  acts on triplets and is in fact a combination of the idea of    partial metrics and $G-$metrics. A partial metric acts on pairs and  has the property where the distance of a point to itself need not be zero. The $G-$metric on the other hand, acts on triplets but the distance of a triplet made up of the same point remains zero ($G(x,x,x)=0$). As a generalization of both, the $G_p-$metric acts on triplets and the distance of a triplet made up of the same point need not be zero.
\begin{defn}\label{D2e1} A \uline{\textbf{$G_p-$metric}} $G_p$ on a set $X$ is a function  $G_p:X^3 \to \mathbb{R}$ satisfying the following axioms:
\\For all $ x,y,z,a \in X,$
\\($G_p$-lbnd): $0 \le G_p(x,x,x) \le G_p(x,x,y)\le G_p(x,y,z).$
\\($G_p$-sym): $G_p(x,y,z)=G_p(\Pi\{x,y,z\})$,\hspace{3ex} where $\Pi$ denotes a permutation on $\{x,y,z\}$.
\\($G_p$-sep): $G_p(x,y,z)=G_p(x,x,x)=G_p(y,y,y)=G_p(z,z,z) \iff x=y=z$.
\\($G_p$-inq): $G_p(x,y,z) \le G_p(x,y,a) + G_p(a,a,z)-G_p(a,a,a)$.
\end{defn}
($G_p$-lbnd)  restricts $G_p$ to have exclusively nonnegative values. ($G_p$-lbnd) combined with ($G_p$-sym) forces $G_p(x,x,y)\le G_p(y,y,x)$  by taking $z=y$ and, hence, $G_p(x,x,y)=G_p(y,y,x)$.  In the previous section, and while investigating the $G-$metric, we  found  (G-lbnd) to be too restrictive. We proposed the $n-\mathfrak{M}$etric such that (n-lbnd) of the $3-\mathfrak{M}$etric was a weakening of (G-lbnd). We now proceed to generalize the $n-\mathfrak{M}$etric into what we call a partial $n-\mathfrak{M}$etric. \begin{defn}\label{D2e2} A \textbf{\uline{partial $n-\mathfrak{M}$etric}} $P$ on a set $X$ is a function $P:X^{n}\to\mathbb{R}$ satisfying the following axioms:
\\For all $(\langle x_i \rangle_{i=1}^n,a)\in X^{n+1},$
\\($P_n$-lbnd):  $ P(\langle x_1 \rangle^n) \le P(\langle x_1 \rangle^{n-1},x_2)$.
\\($P_n$-sym):  $P(\langle x_i \rangle_{i=1}^n) = P(\langle x_{\pi(i)} \rangle_{i=1}^n)$, \hspace{3ex} where $\pi$ is a permutation on $\{1,....,n\}$.
\\($P_n$-sep):  $P(\langle x_1 \rangle^{n-1},x_2)=P(\langle x_1 \rangle^n)$ and $P(\langle x_2 \rangle^{n-1},x_1)=P(\langle x_2 \rangle^n) \iff x_1=x_2$.
\\($P_n$-inq): $P(\langle x_i \rangle_{i=1}^n)\le P(\langle x_i \rangle_{i=1}^{n-1},a) + P(\langle a \rangle^{n-1},x_n)-P(\langle a \rangle^{n})$.
\end{defn}
\begin{rmk}\label{R2e1} A $G_p-$metric is a special case of a partial $3-\mathfrak{M}$etric. A partial $2-\mathfrak{M}$etric is the same as a partial metric.
\end{rmk}
\indent\indent The  cornerstone of  the above generalization is to allow  $P(\langle x_1 \rangle^n)$ to have non-zero values. The axioms were amended to flow with this adjustment. Note that in our proposed partial $3-\mathfrak{M}$etric, our ($P_n$-lbnd) is weaker than the ($G_p$-lbnd), allowing $P$ to have negative values. This weakening also allows  $P(x,x,y)$ and $P(y,y,x)$ to be  related only by the triangular inequality. i.e. $P(x,x,y)$ and $P(y,y,x)$ need not be equal. \\\indent We start by generalizing \Cref{T2d1} to the partial $n-\mathfrak{M}$etric case. 
\begin{thm}{\textbf{(Term Replacement):}}
\label{T2e1}
\\Let $P$ be a partial $n-\mathfrak{M}$etric on a set $X$. For all $ \langle x_i \rangle_{i=1}^n,\langle y_i \rangle_{i=1}^n\in X^{n}$ and for $t\in\{1,.......,n\}$, 
\\$$P(\langle x_i \rangle_{i=1}^{n}) \le P(\langle y_j \rangle_{j=1}^{t},\langle x_i\rangle_{i=t+1}^n)+\sum\limits_{j=1}^t[P(\langle y_j \rangle^{n-1},x_j)-P(\langle y_j\rangle^n)].$$
\end{thm}
\textbf{Proof:} Let $ \langle x_i \rangle_{i=1}^n,\langle y_i \rangle_{i=1}^n\in X^{n}.$ For $t=1$, the result follows from ($P_n$-inq) and ($P_n$-sym).
\\\indent Let $t \in \{2,....,n-1\}$ and assume the inequality holds for $t-1$. Then $$P(\langle x_i \rangle_{i=1}^{n}) \le P(\langle y_j \rangle_{i=1}^{t-1},\langle x_i\rangle_{i=t}^n)+\sum\limits_{j=1}^{t-1}[P(\langle y_j \rangle^{n-1},x_j)-P(\langle y_j\rangle^n)]$$
by ($P_n$-sym) $$=P(\langle y_j \rangle_{i=1}^{t-1},\langle x_i\rangle_{i=t+1}^n,x_t)+\sum\limits_{j=1}^{t-1}[P(\langle y_j \rangle^{n-1},x_j)-P(\langle y_j\rangle^n)]$$
by ($P_n$-inq) $$\le P(\langle y_j \rangle_{i=1}^{t-1},\langle x_i\rangle_{i=t+1}^n,y_t)+P(\langle y_t\rangle^{n-1},x_t)-P(\langle y_t \rangle^n)+\sum\limits_{j=1}^{t-1}[P(\langle y_j \rangle^{n-1},x_j)-P(\langle y_j\rangle^n)]$$
by ($P_n$-sym) $$=P(\langle y_j \rangle_{i=1}^{t},\langle x_i\rangle_{i=t+1}^n)+\sum\limits_{j=1}^{t}[P(\langle y_j \rangle^{n-1},x_j)-P(\langle y_j\rangle^n)]. \hspace{6ex} \square$$
\begin{rmk}\label{R2e2}
As does \Cref{T2d1}, \Cref{T2e1} gives rise to useful replacement tools which are presented in the three corollaries below.
\end{rmk}
\begin{cor}\label{C2e1}Let $P$ be a partial $n-\mathfrak{M}$etric on a set $X$. For all $ \langle x_i \rangle_{i=1}^n,\langle y_i \rangle_{i=1}^n\in X^{n}$ 
$$P(\langle x_i \rangle_{i=1}^{n}) \le P(\langle y_i \rangle_{i=1}^{n})+\sum\limits_{j=1}^n[P(\langle y_j \rangle^{n-1},x_j)-P(\langle y_j \rangle^n)].$$
\end{cor}
\textbf{Proof:} This is the case of  \Cref{T2e1} when  $t=n. \hspace{6ex} \square$
\begin{cor}\label{C2e2}Let P be a partial $n-\mathfrak{M}$etric on a set $X$. Then for  $ a,b \in X$ and $t\in\{1,.......,n\}$, $$P(\langle a\rangle^{t},\langle b\rangle^{n-t}) \le tP(\langle b\rangle^{n-1},a)-(t-1)P(\langle b \rangle ^n ).$$
\end{cor}
\textbf{Proof:} The proof follows from \Cref{C2e1} where $\langle x_i\rangle_{i=1}^n=(\langle a\rangle^t,\langle b\rangle^{n-t})$ and $\langle y_j\rangle_{j=1}^n=\langle b\rangle^n$ we get $$P(\langle a\rangle^t,\langle b\rangle^{n-t})\le P(\langle b \rangle^n)+\sum\limits_{j=1}^n[P(\langle y_j \rangle^{n-1},x_j)-P(\langle y_j \rangle^n)]$$
$$= P(\langle b \rangle^n)+\sum\limits_{j=1}^t[P(\langle y_j \rangle^{n-1},x_j)-P(\langle y_j \rangle^n)]+\sum\limits_{j=t+1}^n[P(\langle y_j \rangle^{n-1},x_j)-P(\langle y_j \rangle^n)]$$
$$=P(\langle b \rangle^n)+\sum\limits_{j=1}^t[P(\langle b\rangle^{n-1},a)-P(\langle b\rangle^n)]+\sum\limits_{j=t+1}^n[P(\langle b\rangle^{n-1},b)-P(\langle b\rangle^n)]$$
$$=P(\langle b \rangle ^n)+t[P(\langle b\rangle^{n-1},a)-P(\langle b\rangle^n)]$$
$$=tP(\langle b\rangle^{n-1},a)-(t-1)P(\langle b \rangle ^n ). \hspace{6ex} \square$$
\\\indent The above corollary hints to as why we were able to choose such a simple ($P_n$-sep) axiom. A very useful corollary is when $t=n-1$.
\begin{cor}\label{C2e3} Let P be a partial $n-\mathfrak{M}$etric on a set $X$. Then for  $ a,b \in X,$ $$P(\langle a\rangle^{n-1},b) \le (n-1)P(\langle b\rangle^{n-1},a)-(n-2)P(\langle b \rangle^n).$$
\end{cor}
\textbf{Proof:} This is the case of \Cref{C2d2} when $t=n-1.\hspace{6ex} \square$
\\\indent Similarly to the case of an $n-\mathfrak{M}$etric, every partial $n-\mathfrak{M}$etric on a set $X$ induces a metric on $X$. Thus \Cref{T2d2} is a special case of \Cref{T2e2} where for every $x\in X$, $P(\langle x \rangle ^n)=0$.  
\begin{thm}{\textbf{(Metric from  Partial $n-\mathfrak{M}$etric):}}
\label{T2e2} 
\\Let P be a partial $n-\mathfrak{M}$etric on a set $X$. For $x,y \in X$ let 
$$d(x,y)=P(y,\langle x \rangle^{n-1})-P(\langle x \rangle^{n})+P(x,\langle y \rangle^{n-1})-P(\langle y \rangle^{n}).$$
Then $d$ is a metric on the set $X$.
\end{thm}
\textbf{Proof:} Let  $x,y,z \in X$.
\\
\\Proof of (m-bnd): From ($P_n$-lbnd) we have $$P(\langle x \rangle^{n})\le P(y,\langle x \rangle^{n-1}) \text{ and } P(\langle y \rangle^{n})\le P(x,\langle y \rangle^{n-1}).$$ Therefore, $$d(x,y)=P(y,\langle x \rangle^{n-1})-P(\langle x \rangle^{n})+P(x,\langle y \rangle^{n-1})-P(\langle y \rangle^{n}) \ge 0.$$
\\\indent Proof of (n-sym): Symmetry of $d$ follows from the symmetry of addition of real numbers.\\
\\\indent Proof of (n-sep):
\\\indent $(\Leftarrow)$ From the definition of $d$ we get
$$d(x,x)=P(\langle x \rangle^{n})-P(\langle x \rangle^{n})+P(\langle x \rangle^{n})-P(\langle x \rangle^{n})=0.$$
\\\indent $(\Rightarrow)$ If $d(x,y)=0$, then $$P(y,\langle x \rangle^{n-1})-P(\langle x \rangle_{1}^{n})+P(x,\langle y \rangle_{1}^{n-1})-P(\langle y \rangle_{1}^{n})=0.$$
From ($P_n$-lbnd) we know that $$P(y,\langle x \rangle^{n-1})-P(\langle x \rangle_{1}^{n})\ge 0 \text { and }P(x,\langle y \rangle^{n-1})-P(\langle y \rangle^{n})\ge 0.$$
Therefore, $$P(y,\langle x \rangle^{n-1})-P(<x>^{n})=0 \text{ and }P(x,<y>^{n-1})-P(<y>^{n})=0,$$
which means that $$ P(y,<x>^{n-1})=P(<x>^{n}) \text{ and }P(x,<y>^{n-1})=P(<y>^{n}).$$ From ($P_n$-sep) we deduce that $x=y.$
\\
\\\indent Proof of (n-inq): From ($P_n$-inq) and ($P_n$-sym) we get $$d(x,y)=P(y,\langle x \rangle^{n-1})-P(\langle x \rangle^{n})+P(x,\langle y \rangle^{n-1})-P(\langle y \rangle^{n})$$
$$=P(\langle x \rangle^{n-1},y)-P(\langle x \rangle^{n})+P(\langle y \rangle^{n-1},x)-P(\langle y \rangle^{n})$$
$$\le P(\langle x \rangle^{n-1},z)+P(\langle z \rangle^{n-1},y)-P(\langle z \rangle^n)-P(\langle x \rangle^{n})+P(\langle y \rangle^{n-1},z)+P(\langle z \rangle^{n-1},x)-P(\langle z \rangle^n)-P(\langle y \rangle^{n})$$
$$=P(x,\langle z \rangle^{n-1})-P(\langle z \rangle^n)+P(z,\langle x \rangle^{n-1})-P(\langle x \rangle^{n})+P(y,\langle z \rangle^{n-1})-P(\langle z \rangle^n)+P(z,\langle y \rangle^{n-1})-P(\langle y \rangle^{n})$$
$$=d(x,z)+d(z,y). \hspace{6ex} \square$$
\begin{eg}{\textbf{($\{-1,1\}-$Discrete  Partial $n-\mathfrak{M}$etric):}}
\label{E2e1}
\\Consider $X$ to be any arbitrary set. Let $P:X^n \to \mathbb{R}$ be defined by:
\\For all $\langle x_i \rangle_{i=1}^n \in X^n, P(\langle x_i \rangle_{i=1}^n)=$
$\begin{cases}-1 &  \text{if } x_1=x_2=...=x_n.
\\1 &$ otherwise.$ \\
\end{cases}$
\\
\\Then $P$ is a  partial $n-\mathfrak{M}$etric on the set $X$.
\end{eg}
\begin{eg}{\textbf{(Maximum Partial $n-\mathfrak{M}$etric):}}
\label{E2e2}
\\Consider the set $X$ to be a subset of $\mathbb{R}$. For  all $\langle x_i \rangle_{i=1}^n \in X$ let $$P(\langle x_i \rangle_{i=1}^n)=\max \{x_i\}_{i=1}^n.$$
Then $P$ is a  partial $n-\mathfrak{M}$etric on the set $X$.
\end{eg}
\textbf{Proof:} Let $(\langle x_i \rangle_{i=1}^n,a) \in X^{n+1}$.
\\
\\\indent Proof of ($P_n$-bnd): $a \le \max\{a,x_1\}$ and, hence, $P(\langle a \rangle^n)\le P(\langle a \rangle^{n-1},x_1).$
\\
\\\indent Proof of ($P_n$-sym): The maximum of a finite set does not change  under the permutation of the set.
\\
\\\indent Proof of ($P_n$-sep):
\\\indent $(\Leftarrow$)is  trivial.
\\
\\\indent $(\Rightarrow)$ If $$P(\langle x_1 \rangle^{n-1},x_2)=P(\langle x_1 \rangle^n)\text{ and }P(\langle x_2 \rangle^{n-1},x_1)=P(\langle x_2 \rangle^n),$$
then $$\max\{x_1,x_2\}=x_1 \text{ and }\max\{x_1,x_2\}=x_2.$$
Therefore, $x_1=x_2$.
\\
\\Proof of ($P_n$-inq): Without loss of generality and due to ($P_n$-sym), we may assume that $x_1 \le x_2 \le....\le x_n$.
\\Hence, $P(\langle x_i \rangle_{i=1}^n)=x_n$. Three cases arise:
\\\uline{Case 1}: Suppose $a \le x_{n-1}\le x_n.$ Then $$P(\langle x_i \rangle_{i=1}^{n-1},a)+P(\langle a \rangle^{n-1},x_n)-P(\langle a \rangle^n)=x_{n-1}+x_n-a \ge x_n=P(\langle x_i \rangle_{i=1}^n).$$
\uline{Case 2}: Suppose $x_{n-1}\le a \le x_n.$ Then $$P(\langle x_i \rangle_{i=1}^{n-1},a)+P(\langle a \rangle^{n-1}, x_n)-P(\langle a \rangle^n)=a+x_n-a = x_n=P(\langle x_i \rangle_{i=1}^n).$$
\uline{Case 3}: Suppose $x_{n-1}\le x_n\le a $. Then $$P(\langle x_i \rangle_{i=1}^{n-1},a)+P(x_n,\langle a \rangle^{n-1})-P(\langle a \rangle^n)=a+a-a = a\ge x_n=P(\langle x_i \rangle_{i=1}^n).$$
\indent Therefore, for all $ (\langle x_i \rangle_{i=1}^n,a) \in X^{n+1}$,
$$P(\langle x_i \rangle_{i=1}^n)\le P(\langle x_i \rangle_{i=1}^{n-1},a)+P(\langle a \rangle^{n-1},x_n)-P(\langle a \rangle^n). \hspace{6ex} \square$$
\begin{thm}{ \textbf{(Partial $n-\mathfrak{M}$etric from a partial metric):}}
\label{E2e3}
\\Every partial metric $p$ on a set $X$ naturally defines a partial $n-$metric  $P$ on $X$ as follows:
\\For\ all $\langle x_i \rangle_{i=1}^n\in X^n,$ $$P(\langle x_i \rangle_{i=1}^n)=\sum\limits_{t=2}^{n}\sum\limits_{i=1}^{t-1}p(x_i,x_t).$$
\end{thm}
\textbf{Proof:} Let $(\langle x_i \rangle_{i=1}^n,a)\in X^{n+1}$.
\\\indent Proof of ($P_n$-lbnd): From the definition of $P$ we get
$$P(\langle x_1 \rangle^n)=\sum\limits_{t=2}^{n}\sum\limits_{i=1}^{t-1}p(x_1,x_1)=\sum\limits_{t=2}^{n-1}\sum\limits_{i=1}^{t-1}p(x_1,x_1)+\sum\limits_{i=1}^{n-1}p(x_1,x_1).$$
Now using ($P_n$-lbnd) we get
$$P(\langle x_1 \rangle^n) \le \sum\limits_{t=2}^{n-1}\sum\limits_{i=1}^{t-1}p(x_1,x_1)+\sum\limits_{i=1}^{n-1}p(x_1,x_2)=P(\langle x_1 \rangle^{n-1},x_2).$$
\\\indent Proof of ($P_n$-sym): ($P_n$-sym) follows from the above and (p-sym)
\\
\\\indent Proof of (P-ineq): Using the definition of $P(\langle x_i \rangle_{i=1}^n)$ we get $$P(\langle x_i \rangle_{i=1}^n)=\sum\limits_{t=2}^{n}\sum\limits_{i=1}^{t-1}p(x_i,x_t)=\sum\limits_{t=2}^{n-1}\sum\limits_{i=1}^{t-1}p(x_{i},x_t)+\sum\limits_{i=1}^{n-1}p(x_i,x_n)$$
by (p-inq)
$$\le \sum\limits_{t=2}^{n-1}\sum\limits_{i=1}^{t-1}p(x_{i},x_{t})+\sum\limits_{i=1}^{n-1}(p(x_i,a)+p(a,x_n)-p(a,a))$$
$$= \sum\limits_{t=2}^{n-1}\sum\limits_{i=1}^{t-1}p(x_{i},x_{t})+\sum\limits_{i=1}^{n-1}p(x_i,a)+\uline{\uline{0}}+\sum\limits_{i=1}^{n-1}p(a,x_n)-\sum\limits_{i=1}^{n-1}p(a,a)$$
$$= \sum\limits_{t=2}^{n-1}\sum\limits_{i=1}^{t-1}p(x_{i},x_{t})+\sum\limits_{i=1}^{n-1}p(x_i,a)+\uline{\uline{\ \sum\limits_{t=2}^{n-1}\sum\limits_{i=1}^{t-1}p(a,a)- \sum\limits_{t=2}^{n-1}\sum\limits_{i=1}^{t-1}p(a,a)}}+\sum\limits_{i=1}^{n-1}p(a,x_n)-\sum\limits_{i=1}^{n-1}p(a,a)$$
$$= \sum\limits_{t=2}^{n-1}\sum\limits_{i=1}^{t-1}p(x_{i},x_{t})+\sum\limits_{i=1}^{n-1}p(x_i,a) +  \sum\limits_{t=2}^{n-1}\sum\limits_{i=1}^{t-1}p(a,a)+\sum\limits_{i=1}^{n-1}p(a,x_n)-\uline{\uline{\uline{ \sum\limits_{t=2}^{n-1}\sum\limits_{i=1}^{t-1}p(a,a)-\sum\limits_{i=1}^{n-1}p(a,a)}}}$$
$$= \sum\limits_{t=2}^{n-1}\sum\limits_{i=1}^{t-1}p(x_{i},x_{t})+\sum\limits_{i=1}^{n-1}p(x_i,a) +  \sum\limits_{t=2}^{n-1}\sum\limits_{i=1}^{t-1}p(a,a)+\sum\limits_{i=1}^{n-1}p(a,x_n)-\uline{\uline{\uline{ \sum\limits_{t=2}^{n}\sum\limits_{i=1}^{t-1}p(a,a)}}}$$
$$=P(\langle x_i \rangle^{n-1},a)+P(\langle a \rangle^{n-1},x_{n})-P(\langle a \rangle^n).\hspace{6ex} \square$$
\section{Strong  Partial $n-\mathfrak{M}$etric}
\label{C2f}
\indent \indent In \Cref{C2c} we presented  the strong partial metric, a generalized metric aimed at simulating  scoring schemes set up to compare two finite sequences.\\\indent In bio-informatics \cite{Bio2008}, we have  scoring schemes which allow us to align and compare multiple DNA strands at the same time. We call those types of schemes multiple sequence alignment schemes. That is why in \textbf{2015} \cite{Ass20152}, we introduced a stronger version of the  partial $n-\mathfrak{M}$etric by combining ($P_n$-lbnd) and ($P_n$-sep) into a stronger axiom. The result is a generalized metric capable of emulating  a multiple sequence alignment schemes for a set of  finite sequences generated from a finite alphabet. We called it a strong partial $n-\mathfrak{M}$etric.
\begin{defn}\label{D2f1} A \textbf{\uline{strong  partial $n-\mathfrak{M}$etric}} $S$ on a set $X$ is a function $S:X^{n}\to\mathbb{R}$ satisfying the following axioms:
\\For all $ (\langle x_i \rangle_{i=1}^n,a)\in X^{n+1},$
\\($S_n$-lbnd):  $ S(\langle x_1 \rangle^n) <S(\langle x_1 \rangle^{n-1},x_2).$
\\($S_n$-sym):  $S(\langle x_i \rangle_{i=1}^n) = S(\langle x_{\pi(i)} \rangle_{i=1}^n)$, where $\pi$ is a permutation on $\{1,....,n\}$.
\\($S_n$-inq): $S(\langle x_i \rangle_{i=1}^n)\le S(\langle x_i \rangle_{i=1}^{n-1},a) + S(\langle a \rangle^{n-1},x_n)-S(\langle a \rangle^{n})$.
\end{defn}
\begin{rmk}\label{R2f1} A strong partial $2-\mathfrak{M}$etric is a strong partial metric.
\end{rmk}
\indent \indent Notice that an ($S_n$-sep)  is hidden in  ($S_n$-lbnd)  as
$$S(\langle x_1 \rangle^{n-1},x_2)= S(\langle x_1\rangle^n)\iff x_1=x_2.$$
\indent Clearly, a strong  partial $n-\mathfrak{M}$etric $S$ on a set $X$ is a  partial $n-\mathfrak{M}$etric on $X$. Hence, as in \Cref{T2e2}, $S$ induces a metric $d$ on $X$ defined as follows:
\\For\ all $x,y \in X$ $$d(x,y)=S(y,\langle x \rangle^{n-1})-S(\langle x \rangle^{n})+S(x,\langle y \rangle^{n-1})-S(\langle y \rangle^{n}).$$
\begin{eg}{\textbf{(Shifted $n-\mathfrak{M}$etric):}}
\label{E2f1}
\\Let $M:X ^{n}\to \mathbb{R}$ be an $n-\mathfrak{M}$etric defined on a set $X$. For any $ r \in \mathbb{R}$, let $S_r:X^{n} \to \mathbb{R}$ be defined by: $$S_r(\langle x_i \rangle_{i=1}^n)=M(\langle x_i \rangle_{i=1}^n)+r.$$
Then $S_r$ is a strong partial $n-\mathfrak{M}$etric on the set $X$.
\end{eg}
\textbf{Proof:}
Let $ (\langle x_i \rangle_{i=1}^n,a)\in X^{n+1}$.
\\
\\\indent Proof of ($S_n$-lbnd): For $a \ne x_1$ ,$$S_{r}(\langle a \rangle^n)=M(\langle a \rangle^n)+r=r.$$ 
From (n-lbnd) and (n-sep) and since $a\ne x_1$ we get $$M(\langle a \rangle^{n-1},x_1)>0.$$
Hence, $$S_{r}(\langle a \rangle^n)=r<M(\langle a \rangle^{n-1},x_1)+r=S_{r}(\langle a \rangle^{n-1},x_1).$$
\\\indent Proof of ($S_n$-sym): Follows directly from (n-sym).
\\
\\\indent Proof of ($S_n$-inq): From (n-inq) we know that $$M(\langle x_i \rangle_{i=1}^{n}) \le M(\langle x_i \rangle_{i=1}^{n-1},a) + M(\langle a \rangle^{n-1},x_n)$$
$$=M(\langle x_i \rangle_{i=1}^{n-1},a) + M(\langle a \rangle^{n-1},x_n)-M(\langle a \rangle^{n}).$$
By (n-sep), $M(\langle a \rangle^{n})=0$.
Therefore,
$$S_r(\langle x_i \rangle_{i=1}^{n})= M(\langle x_i \rangle_{i=1}^{n})+r$$
$$=M(\langle x_i \rangle_{i=1}^{n})+2r-r-M(\langle a\rangle^{n})$$
$$\le M(\langle x_i \rangle_{i=1}^{n-1},a) + M(\langle a \rangle^{n-1},x_n)+2r-r-M(\langle a \rangle^{n})$$
$$=[M(\langle x_i \rangle_{i=1}^{n-1},a)+r ]+ [M(\langle a \rangle_{1}^{n-1},x_n)+r]-[M(\langle a \rangle_{i=1}^{n})+r]$$
$$=S_r(\langle x_i \rangle_{i=1}^{n-1},a) + S_r(\langle a \rangle^{n-1},x_n)-S_r(\langle a \rangle^{n}). \hspace{6ex} \square$$
\begin{thm}{\textbf{(Strong partial $n-\mathfrak{M}$etric from a strong partial metric):}}
\label{E2f2}
\\Every strong partial metric $s$ on a set $X$ naturally defines a strong partial $n-$metric  $S$ on $X$ as follows:
\\For all $\langle x_i \rangle_{i=1}^n\in X^n,$ $$S(\langle x_i \rangle_{i=1}^n)=\sum\limits_{t=2}^{n}\sum\limits_{i=1}^{t-1}s(x_i,x_t).$$
\end{thm}
\textbf{Proof:} For the proof of ($S_n$-sym) and ($S_n$-inq), please refer to \Cref{E2e1}.
\\
\\\indent Proof of ($S_n$-lbnd): Let $ a$ and $x_{1}$ be two distinct elements of $ X$. Using the definition of $S(\langle a \rangle^n)$  we get
$$S(\langle a \rangle^n)=\sum\limits_{t=2}^{n}\sum\limits_{i=1}^{t-1}s(a,a)=\sum\limits_{t=2}^{n-1}\sum\limits_{i=1}^{t-1}s(a,a)+\sum\limits_{i=1}^{n-1}s(a,a).$$
By (s-lbnd) $s(a,a)<s(a,x_1)$ and, hence,
$$ S(\langle a \rangle^n)< \sum\limits_{t=2}^{n-1}\sum\limits_{i=1}^{t-1}s(a,a)+\sum\limits_{i=1}^{n-1}s(a,x_1)=S(\langle a\rangle^{n-1},x_1). \hspace{6ex} \square$$
\indent The above example is a general form of a multiple sequence alignment scheme \cite{Bio2008}, since each pairwise alignment scheme is a strong partial metric. (see \Cref{E2c3}.)
\chapter{Topology}
\label{C3}
\setcounter{thm}{0}
\setcounter{defn}{0}
\indent In \textbf{1914}, Hausdorff \cite{Hau1914} used set theory to generalize a Euclidean space while retaining concepts such as continuity, convergence and connectedness. He used the metric defined by Fr\'echet \cite{Fre1906} (see \Cref{D2a1}) to generate a topological space called a metric space.
The definitions and lemmas in this subsection and in \Cref{C3a}, with the exception of  \Cref{D3a3}, can be found in \cite{Mun2000} and every book on undergraduate topology.
\begin{defn}\label{D3-1} A \textbf{\uline{topology }} $\mathcal{T}$ on a set $X$ is a subset of $\mathcal{P}(X)$ satisfying the following axioms:
\\(R1): $ \emptyset \in \mathcal{T}$ and $X \in \mathcal{T.}$
\\(R2): If $\{O_i\}_{i\in \mathcal{I}}\subseteq \mathcal{T}$ then $\bigcup\limits_{i \in \mathcal{I}} O_i \in \mathcal{T}$. I.e. an arbitrary union of elements of $\mathcal{T}$ is an element of $\mathcal{T}$.
\\(R3): If  $\{O_i\}_{i=1}^n \subseteq \mathcal{T}$ and $n$ is a positive integer, then  $\bigcap\limits_{i=1}^{n} O_i \in \mathcal{T}$. I.e. a finite intersection of elements of $\mathcal{T}$ is an element of $\mathcal{T}$.
\\
\\We call the pair $(X,\mathcal{T})$ a \textbf{\uline{topological space}}.
\end{defn}
\begin{defn}\label{D3-2} The \textbf{\uline{open sets}} of a topological space $(X,\mathcal{T})$ are the elements of $\mathcal{T}$. While the \textbf{\uline{closed sets}} are the complement of the  elements of $\mathcal{T. }$
I.e.
\begin{center}
$O \in X$ is called an \uline{\textbf{open set}} of $X$ $\iff O\in \mathcal{T}$.
\\$F \in X$ is called a \textbf{\uline{closed set}} of $X$ $\iff$ $X-F \in \mathcal{T}$.
\end{center}
\end{defn}
\begin{defn}\label{D3-3} Let $(X,\mathcal{T})$ be a topological space. For each $x$ in $X$, every open set containing $x$ is called an\ \textbf{\uline{open neighborhood}} of $x$.
\end{defn}
\begin{lem}\label{L3-1} Let $(X,\mathcal{T})$ be a topological space and the set $U\subseteq X$. $U$ is an open set if and only if for each $x\in U$ there exists an open neighborhood $V_x$ of $x$ such that $V_x\subseteq U$.
\end{lem}
\indent When considering a topological space one usually does not need to deal with all open sets. \Cref{L3-1} makes it sufficient  to deal with a sub-collection of open sets called a basis. A key feature of a basis of a topology is that it can be used to generate all open sets \cite{Sch1964}.
\newpage
\begin{defn}\label{D3-4} A \textbf{\uline{basis}} $\mathcal{B}$ on a set $X$ is a subset of $\mathcal{P}(X)$ satisfying the following axioms:
\\
\\(B1): For all $ x \in X$, there exists $B \in \mathcal{B}$ such that $x \in B $. 
\\(B2): For all $B_1,B_2 \in \mathcal{B}$ and for each element $x\in B_1\cap B_2 $, there exists $ B_3 \in \mathcal{B}$ such that $x \in B_3\subseteq B_1 \cap B_2$.
\end{defn}
\indent \indent Given a basis $\mathcal{B}$  on a set $X$, let $\mathcal{T_B}=\{\bigcup \mathcal{C} \hspace{1ex} |\hspace{1ex}\mathcal{C} \in \mathcal{B} \}\bigcup \{\emptyset\}$. Then $\mathcal{T_B}$ is a topology on $X$ \cite{Mun2000}. We say that $\mathcal{T_B}$ is the topology generated by $\mathcal{B}$.
\\\indent  In the case of metric spaces the following  well-known technique is used \cite{Mat1994} to prove that a certain subset of $\mathcal{P}(X)$ is a basis.
\begin{lem}\label{L3-2} Given a set $X$ and a function $b:\mathbb{R}^{> 0} \times X \to \mathcal{P}(X)$ satisfying the  conditions:
\\$(b_1)$: For all $ x \in X$ and for each  $\alpha \in \mathbb{R}^{>0}$, $x \in b(\alpha,x)$.
\\$(b_2)$: For all $ x \in X$ and for all  positive real numbers $\alpha \le \beta$,  $b(\alpha,x)\subseteq b(\beta,x)$.
\\$(b_3)$: For all $x,y \in X$ and $\alpha \in \mathbb{R}^{>0}$ such that $y\in b(\alpha,x)$, there exists a positive real number $\epsilon$ such that $b(\epsilon, y) \subseteq b(\alpha ,x)$.
\\ Then $b(\mathbb{R}^{> 0} \times X)$ is a basis on $X$.
\end{lem}
\begin{defn}\label{D3-5} Let $(X,\mathcal{T})$ be a topological space and a set $A\subseteq X$. The \textbf{\uline{closure}} of $A$, denoted $\bar A$ or $Cl(A)$, is the intersection of all closed sets containing $A$.
\end{defn}
\begin{lem}\label{L3-3}Let $(X,\mathcal{T})$ be a topological space and a set $A\subseteq X$. $Cl(A)$ is the smallest closed set containing $A$.
\end{lem}
\begin{lem} \label{L3-4}Let $(X,\mathcal{T})$ be a topological space and a set $A\subseteq X$. $x\in Cl(A)$ if and only if every open neighborhood of $x$ intersects $A$.
\end{lem}
\indent To prove Banach-type fixed point theorems on a topological space, we need to be able to topologically distinguish distinct points.\begin{defn}\label{D3-6} Let $(X,\mathcal{T})$ be a topological space. We say that $(X,\mathcal{T})$ is $\uline{\mathbf{T_0}}$ if for every two distinct elements $x$ and $y$ of $X$, there exists an open set $ U$ in $\mathcal{T}$ such that
either$$[x \in U \text{ and } y \notin U]\text{ \textbf{\uline{or}} }[x \notin U \text{ and }y \in U].$$
\end{defn}
\begin{defn}\label{D3-7} Let $(X,\mathcal{T})$ be a topological space. We say that $(X,\mathcal{T})$ is $\uline{\mathbf{T_1}}$ if for every two distinct elements $x$ and $y$ of $X$, there exist two open sets $U$ and $V$ in $\mathcal{T}$ such that$$[x \in U \text{ and } y \notin U]\text{  \textbf{\uline{and}} }[x \notin V \text{ and } y \in V].$$
\end{defn}
\begin{defn}\label{D3-8} Let $(X,\mathcal{T})$ be a topological space. We say that $(X,\mathcal{T})$ is $\uline{\mathbf{T_2}}$ (or \textbf{\uline{Hausdorff}}) if for every two distinct elements $x$ and $y$ of $X$, there exists two open sets $U$ and $V$ in $\mathcal{T}$ such that
$$x\in U,\hspace{1ex} y \in V \text{and } U \cap V = \emptyset.$$
\end{defn}
\indent \indent Clearly if a topological space is  $T_2$ then it is  $T_1$ and if a topological space is $T_1$ then it is $T_0$ \cite{Mun2000}.
\begin{defn}\label{D3-9} Let $(X,\mathcal{T})$ be a topological space. We say that $(X,\mathcal{T})$ is \uline{\textbf{first countable}} if for every $x\in X$, there exists a countable collection $\mathcal{B}_x=\{B_i\}_{i\in \mathbb{N}}\subseteq\mathcal{T}$ such that for each open set $U$ containing $x$ there is an index $i\in\mathbb{N}$ with $x\in B_i\subseteq U$. We call $\mathcal{B}_x$ a countable local basis at $x$.
\end{defn}
\indent \indent Finally we introduce the notion of ``partial ordering" on topologies of on a particular domain.
\begin{defn}\label{D3-10} Let $(X,\mathcal{T})$ and $(X,\mathcal{T}')$ be two topological spaces.
\\We say $ \mathcal{T}$ is \uline{\textbf{coarser}} than $\mathcal{T}'$ or  $\mathcal{T}'$ is \uline{\textbf{finer}} than $\mathcal{T}$ if $\mathcal{T}\subseteq \mathcal{T'}$.
Alternatively, we may denote it as 
\\$(X,\mathcal{T})$ is coarser than $(X,\mathcal{T}')$ or $(X,\mathcal{T}')$ is finer than $(X,\mathcal{T})$. \end{defn} 
\indent \indent In the sections to follow, we will use the generalized metrics defined in  \Cref{C2} to define  topologies, thus, giving us the structure needed for the study of convergence.
\section{Metric Space}
\label{C3a}
\begin{defn}\label{D3a1} Let $d$ be a metric on a set $X$. For each element $x \in X$ and  positive real number $\epsilon$,  the \uline{\textbf{d-open ball}} around $x$ of radius $\epsilon$ is $$B^d_\epsilon(x)=\{y \in X \hspace{1ex}\vert \hspace{1ex} d(x,y)<\epsilon \}.$$
\end{defn}
\begin{lem} \label{L3a1} Let $d$ be a metric on a set $X$. The collection of all d-open balls on $X$, $\mathcal{B}^d =\{B^d_\epsilon (x)\}_{x\in X}^{\epsilon \in \mathbb{R}^{>0}}$ forms a basis on $X$.
\end{lem}
\begin{defn}\label{D3a2} Let $(X,\mathcal{T})$ be a topological space. We say that $(X,\mathcal{T})$ is a \textbf{\uline{metric space}} if there exists a metric $d$ on $X$ such that $$\mathcal{T} =\mathcal{T}_{\mathcal{B}^d}.$$
(See \Cref{D3-4}.)
\end{defn}
\begin{nott}\label{N3a1} Let  $d$ be a metric on a set $X$. We denote  the topological space $(X,\mathcal{T}_{\mathcal{B}^d})$ by $(X,d)$.
\end{nott}
\begin{lem}\label{L3a2} Every metric space $(X,\mathcal{T})$ is first countable.
\end{lem}
\begin{lem}\label{L3a3} Every metric space $(X,d)$ is $T_2$.
\end{lem}
\begin{defn}\label{D3a3} Let $(X,d)$ be a metric space. For each element $x \in X$ and  positive real number  $\epsilon$, the \textbf{\uline{d-gilded ball}} around $x$ of radius $\epsilon$ is $$\tilde B^d_\epsilon(x)=\{y \in X \hspace{1ex}\vert\hspace{1ex} d(x,y) \le \epsilon \}.$$
\end{defn}
\begin{lem}\label{L3a-4} If $(X,d)$ is a metric space then every d-gilded ball is closed in $(X,d)$.
\end{lem}
\begin{rmk}\label{R3a1} In the literature, a d-gilded ball is referred to as a d-closed ball. When discussing a partial metric $p$, a p-gilded ball  need not be a closed set in the relevant topology. Hence, we chose to change the usual name to avoid ambiguity.
\end{rmk}
\section{Partial Metric  and Strong Partial Metric Space}
\label{C3b}
\indent \indent Although other possible definitions were given for an open ball relative to a partial metric \cite{Mat1994}, we felt that O'Neill's definition \cite{One1996} is the most natural generalization of a d-open ball into the partial metric case.
\begin{defn}\label{D3b1} Let $p$ be  a partial metric on a set $X$. For each element $x \in X$ and  positive real number $\epsilon$, the \textbf{\uline{p-open ball}} around $x$ of radius $\epsilon$ is
$$B_\epsilon^{p}(x)=\{y \in X\hspace{1ex} \vert \hspace{1ex}p(x,y)-p(x,x)<\epsilon \}.$$
\end{defn}
\indent \indent Matthews (\cite{Mat2009}: Definition 13) noted that   
$$B_\epsilon^{p \star}(x)=\{y \in X \hspace{1ex}\vert\hspace{1ex} p(x,y)-p(y,y)<\epsilon \}$$
is another  possible definition which may give us a totally different topology.
 He also presented us with \Cref{L3b1} which is in a fact a special case of our \Cref{L3d1}.
\begin{lem}\label{L3b1} Let $p$ be a partial metric on a set $X$. The collection of all p-open balls on $X$, 
\\$\mathcal{B}^p =\{B_\epsilon ^p (x)\}_{x\in X}^{\epsilon \in \mathbb{R}^{>0}}$ forms a basis on $X$.
\end{lem}
\textbf{Proof:} Use \Cref{L3-2}.
\begin{defn}\label{D3b2} Let $(X,\mathcal{T})$ be a topological space. We say that $(X,\mathcal{T})$ is a  \uline{\textbf{(strong) partial metric space}} if there exists a (strong) partial metric $p$ on $X$ such that 
$$\mathcal{T}=\mathcal{T}_{\mathcal{B}^p}.$$
(See \Cref{D3-4}.)
\end{defn}
\begin{nott}
\label{N3b1} Let $p$ be a partial metric on a set $X$. We denote the topological space $(X,\mathcal{T}_{\mathcal{B}^p})$ by $(X,p)$.
\end{nott}
\begin{lem}\label{L3b2}Every partial metric space $(X,\mathcal{T})$ is first countable.
\end{lem}
\textbf{Proof:} Trivial.
\begin{lem}\label{L3b3}(Matthews \cite{Mat2009}) Every partial metric space $(X,\mathcal{T})$ is $T_0$.
\end{lem}
\textbf{Proof:} Let $(X,\mathcal{T})$ be a partial metric space. By \Cref{D3b2}, let $p$ be a partial metric on $X$ such that
$$\mathcal{T}=\mathcal{T}_{\mathcal{B}^p}.$$
Consider two distinct elements $x$ and $y$ in $X$. From (p-lbnd) and (p-sep) we get
$$p(x,y)>p(x,x) \text{ or }p(x,y)>p(y,y).$$
Then $$\epsilon _x = p(x,y)-p(x,x)>0 \text{ or } \epsilon _y = p(x,y)-p(y,y)>0.$$
Without any loss of generality we may assume that $\epsilon _x >0.$ Hence, $B^p_{\epsilon_x}$ exists. Additionally 
$$x \in B^p_{\epsilon_x} (x) \text{ and } y \notin B^p_{\epsilon_x}(x) \hspace{6ex} \square$$
\indent Note that a partial metric space need not be $T_1$.
\begin{eg}{\textbf{(A partial metric space that is not $T_1$):}}
\label{E3b1}
\\Let $p$ be a partial metric on $X=\mathbb{R}\cup\{a\}$ where $a\notin\mathbb{R}$ as defined in  \Cref{E2b3} by
$$p(a,a)=0,p(a,x)=\vert x \vert \text{ and }p(x,y)=\vert x-y \vert -1.$$
Then  $(X,p)$  is not $T_1$.
 \end{eg}
 \textbf{Proof:} We know that $a\ne 0$ and for every $\epsilon \in \mathbb{R}^{>0}$, $$p(a,0)-p(a,a)=0-0=0<\epsilon.$$ Hence, $0 \in B^p_\epsilon (a) . \hspace{6ex} \square$
\begin{lem}\label{L3b4} Every strong partial metric space $(X,s)$ is $T_1$.
\end{lem}
\textbf{Proof: }Let $(X,\mathcal{T})$ be a strong partial metric space. By \Cref{D3b2}, let $s$ be a strong partial metric on $X$ such that
$$\mathcal{T}=\mathcal{T}_{\mathcal{B}^s}.$$
Consider two distinct elements $x$ and $y$ in $X$. From (s-lbnd)  we get
$$s(x,y)>s(x,x) \text{ and } s(x,y)>s(y,y).$$
Let $$\epsilon _x = s(x,y)-s(x,x)>0 \text{ and } \epsilon _y = s(x,y)-s(y,y)>0.$$
Then,   $$x \in B^s_{\epsilon_x} (x)\text{ , }y \notin B^s_{\epsilon_x}(x) \text{ , } y \in B^s_{\epsilon_y} (y) \text{ , and } x \notin B^s_{\epsilon_y}(y).\hspace{6ex}\square$$
\\Note that a strong partial metric space need not be $T_2$.
\begin{eg}{\textbf{(A strong partial metric space that is not $T_2$):}}
\label{E3b2}
\\Let $s$ be the strong partial metric on $X=(0,+\infty)$  defined in  \Cref{E2c2} by 
$$  s(x,y)=\begin{cases}x & \text{ if } x=y. \\
x+y & \text{ if } x\ne y. \\
\end{cases}   $$
Then $(X,s)$  is not $T_2$.
 \end{eg}
 \textbf{Proof:} For each element $x \in X$ and positive real number $\epsilon$,
$$B^s_\epsilon(x)=\{y \in X \hspace{1ex}\vert\hspace{1ex} s(x,y)-s(x,x)<\epsilon\}=\{x\}\cup\{y\in X \hspace{1ex}\vert\hspace{1ex} y<\epsilon\}.$$
Consider two distinct elements $x$ and $y$ in $X$ and  two positive real numbers $\epsilon$ and $\delta$.
\\Let $z < \min\{\epsilon,\delta\}$ then $z \in B^s_\epsilon(x)\cap B^s_\delta(y). \hspace{6ex}\square$
\begin{defn}\label{D3b3} Let $(X,p)$ be a partial metric space. For each element $ x \in X$ and positive real number $\epsilon$, the \textbf{\uline{p-gilded ball}} around $x$ of radius $\epsilon$ is
$$ \tilde B^p_\epsilon(x)=\{y \in X \hspace{1ex}\vert\hspace{1ex} p(x,y)-p(x,x) \le \epsilon \}.$$
\end{defn}
\indent \indent The next example shows that if $p$ is a partial metric on a set $X$, a p-gilded ball need not be closed in $(X,p).$
\begin{eg}{\textbf{(p-gilded but not closed ball):}}
\label{E3b3}
\\Let $s$ be the strong partial metric on $X=(0,+\infty)$  defined in  \Cref{E2c2} by 
$$  s(x,y)=\begin{cases}x & \text{ if } x=y. \\
x+y & \text{ if } x\ne y. \\
\end{cases}   $$
Then, for each $ x \in X$ and positive real number $\epsilon$, $\tilde B^s_\epsilon(x)$ is not closed in $(X,s)$.
\end{eg}
\textbf{Proof:} Let $x\in X$ and $\epsilon$ a positive real number. Then
$$\tilde B^s_\epsilon(x)=\{y\in X \hspace{1ex} \vert \hspace{1ex} y \le \epsilon \}\cup\{x\}.$$
For each $y\in X$ and every positive real number $\delta$,
$$\emptyset \ne (0, \min \{\epsilon,\delta \})\subseteq B^s_\delta(y)\cap \tilde B^s_\epsilon(x).$$
Hence, $y$ is in the closure of $\tilde B^s_\epsilon(x)$. I.e. the closure of $\tilde B^s_\epsilon(x)$ is $X$ and, therefore, $\tilde B^s_\epsilon(x)$ is not closed in $(X,s)$.  \hspace{6ex} $\square$
\section{$G-$metric and $n-\mathfrak{M}$etric Space}
\label{C3c}
\indent \indent Mustapha and Sims \cite{Mus2006} defined a topology on a set $X$  having a  $G-$metric (see \Cref{D2d1}). Later Khan \cite{Kha2012} used the $K-$metric, a stronger version of our $n-\mathfrak{M}$etric (see \Cref{D2d2}), on a set $X$ to  define a topology on $X$. He called this topological space a $K-$metric space \cite{Ass20152}.
\\\indent In this section we will develop the  $n-\mathfrak{M}$etric space, which is in fact a generalization of the $G-$metric space and the  $K-$metric space. We will also prove that an $n-\mathfrak{M}$etric space  is simply a metric space.
Hence, we will be dropping $n-\mathfrak{M}$etrics as of  \Cref{C4}. 
\begin{defn}\label{D3c1} Let $M$ be an $n-\mathfrak{M}$etric on a set $X$. For each element $x \in X$ and  positive real number $\epsilon$, the \uline{\textbf{M-open ball} }around $x$ of radius $\epsilon$ is
$$B_\epsilon^{M}(x)=\{y \in X \hspace{1ex}\vert\hspace{1ex} M(\langle x\rangle^{n-1},y)<\epsilon \}.$$
\end{defn}
\indent \indent Note that this is not the only possible definition.  For example
$$B_\epsilon^{M \star}(x)=\{y \in X \hspace{1ex}\vert\hspace{1ex} M(\langle y\rangle^{n-1},x))<\epsilon \}$$
 is another possible definition but  which will still generate the same metric topology.
\begin{lem}\label{L3c1} Let $M$ be an $n-\mathfrak{M}$etric on a set $X$. The collection of all M-open balls on $X$, 
\\$\mathcal{B}^M =\{B_\epsilon ^M (x)\}_{x\in X}^{\epsilon \in \mathbb{R}^{>0}}$ forms a basis on $X$.
\end{lem}
\textbf{Proof:} Consider the function  $b:\mathbb{R}^{> 0} \times X \to \mathcal{P}(X)$ given by $b(\epsilon,x) = B^M_\epsilon (x)$. We now check that $b$ satisfies the conditions of \Cref{L3-2}.
\\$(b_1)$: For every  $x \in X$, by (M-sep) we know that $M(\langle x\rangle^{n-1},x)=0$. Hence, for each positive real number $\epsilon$,  $x \in B^M_\epsilon(x)$.
\\$(b_2)$: Clearly, if $0<\delta \le \epsilon$ then $B^M_\delta(x) \subseteq B^M_\epsilon (x)$ for each $x \in X$ by the transitivity of the order on $\mathbb{R}$. Hence, $b(\delta, x) \subseteq b(\epsilon,x)$. 
\\$(b_3)$: For every $x \in X$ and $y \in B^M_\epsilon (x)$, we know that $M(\langle x\rangle^{n-1},y)< \epsilon$. Let $\delta = \epsilon -M(\langle x\rangle^{n-1},y)$. For each element $ z \in B^M_\delta (y)$ (i.e $M(\langle y\rangle^{n-1},z) < \delta$), by (M-inq) we get
$$M(\langle x\rangle^{n-1},z)\le M(\langle x\rangle^{n-1},y)+M(\langle y\rangle^{n-1},z)<M(\langle x\rangle^{n-1},y)+ \delta = \epsilon$$
and, hence,  $$M(\langle x\rangle^{n-1},z)< \epsilon.$$
Thus,  $ z \in B^M_\epsilon (x)$ and $B^M_\delta (y)\subseteq B^M_\epsilon(x). \hspace{6ex}\square$
\begin{defn}\label{D3c2} Let $(X,\mathcal{T})$ be a topological space. We say that $(X,\mathcal{T})$ is an \textbf{\uline{$n-\mathfrak{M}$etric space}} if there exists an $n-\mathfrak{M}$etric $M$ on $X$ such that 
$$\mathcal{T}=\mathcal{T}_{\mathcal{B}^M}.$$
\end{defn}
\begin{nott}\label{N3c1} Let $M$ be an $n-\mathfrak{M}$etric on a set $X$. We denote the topological space $(X,\mathcal{T}_{\mathcal{B}^M})$ by $(X,M)$.
\end{nott}
\begin{thm}\label{L3c2} Every $n-\mathfrak{M}$etric space $(X,\mathcal{T})$ is a metric space.
\end{thm}
\textbf{Proof:} Let $(X,\mathcal{T})$ be an $n-\mathfrak{M}$etric space. By \Cref{D3d2}, let $M$ be an $n-\mathfrak{M}$etric on $X$ such that
$$\mathcal{T}=\mathcal{T}_{\mathcal{B}^M}.$$
\indent In \Cref{T2d2} we proved that the function $$d(x,y)=M(\langle x\rangle^{n-1},y)+M(\langle y\rangle^{n-1},x)$$ is a metric  on $X$.
\\\indent First we prove that $(X,d)$ is finer than $(X,M)$. For each element $x \in X$ and  $y \in B^d_\epsilon(x)$ we know that $$d(x,y) < \epsilon.$$
I.e.
$$M(\langle x\rangle^{n-1},y)+M(x,\langle y\rangle^{n-1}) < \epsilon.$$
From (M-lbnd) we know that $$M(x,\langle y\rangle^{n-1})\ge 0$$
and, hence, $$M(y,\langle x\rangle^{n-1})\le M(\langle x\rangle^{n-1},y)+M(x,\langle y\rangle^{n-1}) < \epsilon.$$
Thus, $y \in B^M_\epsilon(x)$ and, as a result, $B^d_\epsilon(x) \subseteq B^M_\epsilon(x).$
\\
\\\indent Similarly we prove that $(X,M)$ is finer than $(X,d)$. For each element $x \in X$ and  $y \in B^M_\epsilon(x)$ we know that $$M(\langle x\rangle^{n-1},y) < \epsilon.$$ Let $\delta = \frac{\epsilon}{n}$. For each element $y \in B^M_\delta(x)$ we know that $$M(\langle x\rangle^{n-1},y)<\delta.$$
From \Cref{C2d3} we get
$$M(x,\langle y\rangle^{n-1})\le (n-1)M(\langle x\rangle^{n-1},y).$$
Hence, $$d(x,y) \le M(\langle x\rangle^{n-1},y)+(n-1)M(\langle x\rangle^{n-1},y)$$ $$=nM(\langle x\rangle^{n-1},y)<n \delta = \epsilon.$$
Thus, $y \in B^d_\epsilon (x)$ and $B^M_\delta(x) \subseteq B^d_\epsilon(x). \hspace{6ex} \square$
\section{Partial $n-\mathfrak{M}$etric and Strong Partial $n-\mathfrak{M}$etric Space}
\label{C3d}
\indent \indent As mentioned in \Cref{C2}, a partial $n-\mathfrak{M}$etric is a generalization  of a partial metric and an $n-\mathfrak{M}$etric. Unlike in the chapters to come, no special techniques are needed other than the ones presented in  \Cref{C3b} and \Cref{C3c}. Hence, we will present our statements without proofs.
\begin{defn}\label{D3d1} Let $P$ be a partial $n-\mathfrak{M}$etric   on a set $X$. For each element $x \in X$ and positive real number $\epsilon$, the \textbf{\uline{P-open ball}} around $x$ of radius $\epsilon$ is
$$B_\epsilon^{P}(x)=\{y \in X\hspace{1ex} \vert\hspace{1ex} P(\langle x\rangle^{n-1},y)-P(\langle x\rangle^{n})<\epsilon \}.$$
\end{defn}
\indent \indent Note that this is not the only possible definition. For example 
$$B_\epsilon^{P \star}(x)=\{y \in X \hspace{1ex}\vert\hspace{1ex} P(\langle y\rangle^{n-1},x)-P(\langle y\rangle^{n})<\epsilon \}$$ 
is another possible definition which may give us a totally different topology.
\begin{lem}\label{L3d1} Let $P$ be a partial $n-\mathfrak{M}$etric on a set $X$. The collection of all P-open balls on $X$,
\\$\mathcal{B}^P =\{B_\epsilon ^P (x)\}_{x\in X}^{\epsilon \in \mathbb{R}^{>0}}$ forms a basis on $X$.
\end{lem}
\begin{defn}\label{D3d2} Let $(X,\mathcal{T})$ be a topological space. We say that $(X,\mathcal{T})$ is a \textbf{\uline{(strong) partial $n-\mathfrak{M}$etric space}} if there exists a (strong) partial $n-\mathfrak{M}$etric $P$ on $X$ such that $$\mathcal{T} =\mathcal{T}_{\mathcal{B}^P}.$$
\end{defn}
\begin{nott}\label{N3d1} Let $P$ be a partial $n-\mathfrak{M}$etric on a set $X$. We denote  the topological space $(X,\mathcal{T}_{\mathcal{B}^P})$ by $(X,P)$.
\end{nott}
\begin{lem}\label{L3d2} Every partial $n-\mathfrak{M}$etric space $(X,\mathcal{T})$ is first countable.
\end{lem}
\begin{lem}\label{L3d3} Every partial $n-\mathfrak{M}$etric space $(X,\mathcal{T})$ is $T_0$.
\end{lem}
\indent  Note that a partial $n-\mathfrak{M}$etric space need not be $T_1$.
\begin{eg}{\textbf{(A partial $n-\mathfrak{M}$etric space that is not $T_1$):}}
\label{E3d1}
\\\indent For $n=2$, a partial $n-\mathfrak{M}$etric space is simply a partial metric space. In \Cref{E3b1} we gave an example of a partial metric space that is not $T_1$.
\end{eg}
\begin{lem}\label{L3d4} Every strong partial $n-\mathfrak{M}$etric space $(X,\mathcal{T})$ is $T_1$.
\end{lem}
\indent Note that a strong partial $n-\mathfrak{M}$etric space need not be $T_2$.
\begin{eg}{\textbf{(A strong partial $n-\mathfrak{M}$etric space that is not $T_2$):}}
\label{E3d2}
\\\indent For $n=2$, a strong partial $n-\mathfrak{M}$etric space is a strong partial metric space. In \Cref{E3b2} we gave an example of a strong partial metric space that is not $T_2$.
\end{eg}
\begin{defn}\label{D3d3} Let $(X,P)$ be a partial $n-\mathfrak{M}$etric space. For each element $ x \in X$ and positive real number $\epsilon$, the \uline{\textbf{P-gilded ball}} around $x$ of radius $\epsilon$ is
$$\tilde B^P_\epsilon(x)=\{y \in X\hspace{1ex} \vert \hspace{1ex}P(\langle x\rangle^{n-1},y)-P(\langle x\rangle^{n}) \le \epsilon \}.$$
\end{defn}
\begin{eg}{\textbf{(P-gilded but not closed ball):}}
\label{E3d3}
\\\indent For $n=2$, a strong partial $n-\mathfrak{M}$etric space is a strong partial metric space. In \Cref{E3b3} we gave an example of an s-gilded ball that is not closed.
\end{eg}
\chapter{Sequences and Limits}
\label{C4}
\setcounter{thm}{0}
\setcounter{defn}{0}
\begin{defn}\label{D4-1} Let $\{x_i\}_{i \in \mathbb{N}}$ be a sequence in a topological space $(X,\mathcal{T})$. We say that $a$ is a \uline{\textbf{limit}} of $\{x_i\}_{i \in \mathbb{N}}$ if and only if for every open neighborhood $U$ of $a$, there exists a natural number $N$ such that for all   $i > N$, $x_i \in U$.
\end{defn}
\indent\Cref{C4} deals with limits, Cauchy sequences and Cauchy pairs. As mentioned in \Cref{C2}, metrics are a special case of partial metrics which in turn are a special case of partial $n-\mathfrak{M}$etrics.  That is why, in \Cref{C4a} we will be presenting the results without proofs. In \Cref{C4b}, we supply the proofs because they are considerably simpler than their counterparts in \Cref{C4c}. Additionally for many users, we expect the level of generality found in \Cref{C4b} to suffice.
\section{Metric Space}
\label{C4a}
The proofs for this section are readily available in \cite{Mun2000}.
\begin{defn}\label{D4a1} Let $\{x_i\}_{i \in \mathbb{N}}$ be a sequence in a metric space $(X,d)$. We say that the sequence $\{x_i\}_{i \in \mathbb{N}}$ is \textbf{\uline{Cauchy}} if and only if for each positive real number $\epsilon$ there exists a natural number $N$ such that for all $i\ge j >N$, 
$$d(x_i,x_j)<\epsilon.$$
\end{defn}
\indent \indent Merging  \Cref{D4-1} and \Cref{D4a1}, we can now identify a limit by analyzing distances.
\begin{lem}\label{L4a1} Let $\{x_i\}_{i \in \mathbb{N}}$ be a sequence in a metric space $(X,d)$. A point $a$ in $X$ is a limit of $\{x_i\}_{i \in \mathbb{N}}$ if and only if for each positive real number $\epsilon$ there exists a natural number $N$ such that for all  $i\ge j>N$,
$$d(a,x_i)< \epsilon.$$
\end{lem}
\begin{lem}\label{L4a2} Let $\{x_i\}_{i \in \mathbb{N}}$ be a Cauchy sequence in a metric space $(X,d)$. If $\{x_i\}_{i \in \mathbb{N}}$ has a limit in $X$ then that limit is unique.
\end{lem}
\begin{defn}\label{D4a2} Let $(X,d)$ be a metric space. We say that $(X,d)$ is \textbf{\uline{complete}} if and only if  every Cauchy sequence in $X$ has a limit in $X$.
\end{defn}
\indent \indent Suppose we now have two different sequences in a topological space. The study of their respective terms and whether or not they get, in some sense, closer to each other proves to be quite useful in common fixed point theorems.
\begin{defn}\label{D4a3} Let $\{x_i\}_{i \in \mathbb{N}}$ and $\{y_i\}_{i \in \mathbb{N}}$ be  sequences in a metric space $(X,d)$. We say that $\{x_i\}_{i \in \mathbb{N}}$ and $\{y_i\}_{i \in \mathbb{N}}$ form a \textbf{\uline{Cauchy pair}} if and only if for each positive real number $\epsilon$ there exists a natural number $N$ such that for all $ i,j>N$,
$$d(x_i,y_j)< \epsilon.$$
\end{defn}
\begin{thm}{\textbf{(Cauchy pair term comparison):}}
\label{T4a1}
\\Let $\{x_i\}_{i \in \mathbb{N}}$ and $\{y_i\}_{i \in \mathbb{N}}$ be two sequences in a  metric space $(X,d)$. The statements below are equivalent.
\\(a) $\{x_i\}_{i \in \mathbb{N}}$ and $\{y_i\}_{i \in \mathbb{N}}$ form a Cauchy pair.
\\(b) For every positive real number $\epsilon$ there exists a natural number $N$ such that for all  $i\ge j> N$,
$$d(x_i,y_j)< \epsilon.$$
(c) For every positive real number $\epsilon$ there exists a natural number $N$ such that for all  $i\ge j> N$,
$$d(x_j,y_i)< \epsilon.$$
\end{thm}
\indent When we introduce a Cauchy pair, not only are we considering a pair of sequences whose corresponding terms  get eventually  arbitrarily close, but also the sequences  are Cauchy.
\begin{lem}\label{L4a3} Let $\{x_i\}_{i \in \mathbb{N}}$ and $\{y_i\}_{i \in \mathbb{N}}$ be a Cauchy pair in a  metric space $(X,d)$. Then $\{x_i\}_{i \in \mathbb{N}}$ and $\{y_i\}_{i \in \mathbb{N}}$ are both Cauchy sequences. Additionally, if $a$ is a limit of $\{x_i\}_{i \in \mathbb{N}}$ then $a$ is also a limit of $\{y_i\}_{i \in \mathbb{N}}$.
\end{lem}
\section{Partial Metric Space}
\label{C4b}
\indent \indent Matthews \cite{Mat1994} generalized the notion of a Cauchy sequence to partial metrics. Even though Matthews' work was restricted to non-negative values in $\mathbb{R}$, we will use his generalization of Cauchy sequences  when considering O'Neill's definition \cite{One1996} of partial metrics.
\begin{defn}\label{D4b1} Let $\{x_i\}_{i \in \mathbb{N}}$ be a sequence in a partial metric space $(X,p)$. We say that the sequence $\{x_i\}_{i \in \mathbb{N}}$ is \textbf{\uline{Cauchy}} if and only if  there exists a real number $r$ such that for each positive real number $\epsilon$ there exists a natural number $N$ where for all  $i\ge j> N$,
$$-\epsilon <p(x_i,x_j)-r< \epsilon.$$
We call $r$  the \textbf{\uline{central distance}} of $\{x_i\}_{i \in \mathbb{N}}$.
\end{defn}
\indent Notice that from \Cref{D4b1}, $p(x_i,x_i)$ tends to $r$. Hence, in the metric case the central distance $r$ of a Cauchy sequence is $0$ and coincides with \Cref{D4a1}.
\\\indent From \Cref{D4-1} we get the  lemma below.
\begin{lem}\label{L4b1} Let $\{x_i\}_{i \in \mathbb{N}}$ be a sequence in a partial metric space $(X,p)$. A point $a$ in $X$ is a limit of $\{x_i\}_{i \in \mathbb{N}}$ if and only if for each positive real number $\epsilon$ there exists a natural number $N$ such that for all   $i> N$,
$$p(a,x_i)-p(a,a)< \epsilon.$$
\end{lem}
\textbf{Proof:} Let $\{x_i\}_{i \in \mathbb{N}}$ be a sequence in a partial metric space $(X,p)$.
\\\indent$(\Rightarrow)$ Let $a$ be  a limit of $\{x_i\}_{i \in \mathbb{N}}$. From \Cref{D3b2} we know that for each positive real number $\epsilon$, $B^p_\epsilon(a)$ is an open neighborhood of $a$. Since $a$ is a limit of $\{x_i\}_{i \in \mathbb{N}}$ then there exists a natural number $N$ such that for all $i>N$,
$$x_i \in B^p_\epsilon(a).$$
Therefore, from \Cref{D3b1}
$$p(a,x_i)-p(a,a)<\epsilon.$$
\indent$(\Leftarrow)$ Assume that for each positive real number $\epsilon$ there exists a natural number $N$ such that for all $i\ge j > N$,
$$p(a,x_i)-p(a,a)< \epsilon.$$
By \Cref{D3b1}
$$x_i \in B^p_\epsilon (a).$$
\indent Consider $U$ an open neighborhood of $a$. From \Cref{D3b2} we know that there exists a positive real number $\epsilon$ such that
$$B^p_\epsilon (a)\subseteq U.$$
Therefore, there exists a natural number $N$ such that for all $i>N$,
$$x_i\in B^p_\epsilon(a)\subseteq U.\hspace{6ex} \square$$ 

\begin{lem}\label{L4b2} Let $(X,p)$ be a partial metric space. Consider a Cauchy sequence $\{x_i\}_{i \in \mathbb{N}}$   in $X$ with  central distance $r$.
 If the point $a$ in $X$ is a limit of $\{x_i\}_{i \in \mathbb{N}}$, then
$$r \le p(a,a).$$
\end{lem}
\textbf{Proof:} Let  $\{x_i\}_{i \in \mathbb{N}}$ be a Cauchy sequence in $X$ with central distance $r$.
 Then, from \Cref{D4b1}, for every positive real number $\epsilon$ there exists a natural number $N_1$ such that for all $i \ge j >N_1$,
$$-\frac{\epsilon}{3}<p(x_i,x_j)-r<\frac{\epsilon}{3}$$
therefore, $$r<p(x_i,x_j)+\frac{\epsilon}{3}.$$
Since $a$ is a limit of $\{x_i\}_{i \in \mathbb{N}}$, from \Cref{L4b1}, there exists a natural number $N_2$ such that for all $i>N_2$,
$$p(a,x_i)-p(a,a)<\frac{\epsilon}{3}$$
therefore, $$p(a,x_i)<p(a,a)+\frac{\epsilon}{3}.$$
\\Now taking $N=\max\{N_1,N_2\}$ we get that for all $i \ge j>N$ 
$$r <p(x_i,x_j)+\frac{\epsilon}{3}$$
by (p-inq)$$ \le p(x_i,a) + p(a,x_j)-p(a,a) + \frac{\epsilon}{3}$$
by (p-sym)
$$p(a,x_i) + p(a,x_j)-p(a,a) + \frac{\epsilon}{3}$$
$$<(p(a,a)+\frac{\epsilon}{3})+(p(a,a)+\frac{\epsilon}{3})-p(a,a)+\frac{\epsilon}{3}$$
$$=p(a,a)+\epsilon.$$
Hence, for every positive real number $\epsilon$
$$r<p(a,a)+\epsilon$$
and, therefore,
$$r \le p(a,a).\hspace{6ex} \square$$
\indent In a metric space, a limit to a Cauchy sequence is unique. In a partial metric space though, that need not be the case.
\begin{eg}{\textbf{(Multiple Limits):}}
\label{E4b1}
\\Let $p$ be a partial metric on $X=\mathbb{R}\cup\{a\}$ where $a\notin\mathbb{R}$ as defined in  \Cref{E2b3} by:
\\For all $x,y\in \mathbb{R}$,
$$p(a,a)=0,p(a,x)=\vert x \vert \text{ and }p(x,y)=\vert x-y \vert -1.$$
Then the sequence $\{\frac{1}{2^n}\}_{n \in \mathbb{N}}$  in $X$ is Cauchy with a central distance $r=-1$. Moreover, $0$ and $a$ are both limits of $\{\frac{1}{2^n}\}_{n \in \mathbb{N}}$.
\end{eg}
\textbf{Proof:} Assuming that $i\ge j$, then 
$$p(\frac{1}{2^i},\frac{1}{2^j})=\vert \frac{1}{2^i}- \frac{1}{2^j}\vert-1$$
$$=\frac{1}{2^j} -\frac{1}{2^i}-1.$$
Therefore,
$$-1-\frac{1}{2^i}<p(\frac{1}{2^i},\frac{1}{2^j})<\frac{1}{2^j}-1.$$
Hence, for every positive real number $\epsilon$ there exists a natural number $N$ such that $\frac{1}{2^N}<\epsilon$. Thus, for all $i\ge j>N$,
$$-1-\epsilon<-1-\frac{1}{2^N}<-1-\frac{1}{2^i}<p(\frac{1}{2^i},\frac{1}{2^j})<\frac{1}{2^j}-1<\frac{1}{2^N}-1<\epsilon-1.$$
From \Cref{D4b1}, $\{\frac{1}{2^n}\}_{n \in \mathbb{N}}$ is a Cauchy sequence with a central distance $r=-1$.
\\\indent Additionally, for all natural numbers $i$,
$$p(0,\frac{1}{2^i})-p(0,0)=\vert0-\frac{1}{2^i}\vert-1-(-1)= \frac{1}{2^i}$$
and
$$p(a,\frac{1}{2^i})-p(a,a)=\vert \frac{1}{2^i}\vert-0=\frac{1}{2^i}.$$
Hence, for every positive real number $\epsilon$ there exists a natural number $N$ such that for
all $i>N$,
$$p(0,\frac{1}{2^i})-p(0,0)= \frac{1}{2^i}<\frac{1}{2^N}<\epsilon$$
and
$$p(a,\frac{1}{2^i})-p(a,a)= \frac{1}{2^i}<\frac{1}{2^N}<\epsilon.$$
Therefore, by \Cref{L4b1},  $0$ and $a$ are both limits of $\{\frac{1}{2^n}\}_{n \in \mathbb{N}}. \hspace{6ex}\square$
\\\indent In the above example both $0$ and $a$ are limits of $\{\frac{1}{2^n}\}_{n \in \mathbb{N}}$. However, $p(a,a)=0 \ne \lim \limits_{i,j\to+\infty} p(\frac{1}{2^i},\frac{1}{2^j})$ while $p(0,0)=-1=\lim \limits_{i,j\to+\infty} p(\frac{1}{2^i},\frac{1}{2^j})$.  Hence, in some sense, $0$ has more significance to the sequence $\{\frac{1}{2^n}\}_{n \in \mathbb{N}}$  than $a$ does. We will call limits like this: special limits.
\begin{defn}\label{D4b2} Let $(X,p)$ be a partial metric space. Consider a Cauchy sequence $\{x_i\}_{i \in \mathbb{N}}$ in $X$   with a central distance $r$. A point $a$ in $X$ is called a \textbf{\uline{special limit}} of  $\{x_i\}_{i \in \mathbb{N}}$ if and only if $a$ is a limit of $\{x_i\}_{i \in \mathbb{N}}$ and $r=p(a,a)$.
\end{defn}
\indent\indent In a partial metric space, the special limit is analogous to the limit in a metric space since, if it exists, it is unique.
\begin{lem}\label{L4b3} Let $\{x_i\}_{i \in \mathbb{N}}$ be a Cauchy sequence in a partial metric space $(X,p)$.
If  $\{x_i\}_{i \in \mathbb{N}}$ has a special limit in $X$ then that special limit is unique.
\end{lem}
\textbf{Proof:} Consider the Cauchy sequence $\{x_i\}_{i \in \mathbb{N}}$ in $X$   with a central distance $r$. If the points $a$ and $b$ of $X$ are both special limits of $\{x_i\}_{i \in \mathbb{N}}$ then, by \Cref{D4b2},
$$p(a,a)=r=p(b,b).$$
Furthermore, for every positive real number $\epsilon$ there exists a natural number $N_1$ such that for all $i>N_1$,
$$r-\frac{\epsilon}{3} < p(x_i,x_j)<r+\frac{\epsilon}{3} $$
i.e. $$-p(x_i,x_i)<-r+\frac{\epsilon}{3}.$$
By \Cref{D4b2}, both $a$ and $b$ are  limits of $\{x_i\}_{i \in \mathbb{N}}$. Then, by \Cref{L4b1}, there exists a natural number $N_2$ such that for all $i>N_2$,
$$p(a,x_i)-r=p(a,x_i)-p(a,a)<\frac{\epsilon}{3}$$
i.e. $$p(a,x_i)<r+\frac{\epsilon}{3}$$
and there exists a natural number $N_3$ such that for all $i>N_3$,
$$p(b,x_i)-r=p(b,x_i)-p(b,b)<\frac{\epsilon}{3}$$
i.e. $$p(b,x_i)<r+\frac{\epsilon}{3}.$$
Hence, using (p-lbnd) we get for every positive real number $\epsilon$ there exists a natural number $N=\max\{N_1,N_2,N_3\}$ such that for all $i>N$,
$$p(a,a)\le p(a,b)$$
by (p-inq) $$\le p(a,x_i)+p(x_i,b)-p(x_i,x_i)$$
by (p-sym) $$=p(a,x_i)+p(b,x_i)-p(x_i,x_i)$$
$$<r+\frac{\epsilon}{3}+r+\frac{\epsilon}{3}-r+\frac{\epsilon}{3}$$
$$=r+\epsilon=p(a,a)+\epsilon.$$
Therefore, $$p(a,a)\le p(a,b)<p(a,a)+\epsilon$$
i.e. $$p(a,a)=p(a,b).$$
Similarly $p(b,b)=p(b,a)$ and, hence, by (p-sep) $a=b. \hspace{6ex}\square$
\\\indent An additional property of a special limit is that it preserves a notion of sequential continuity.
\begin{lem}\label{L4b4} Let $(X,p)$ be a partial metric space. Let $a$ be the special limit of\ the Cauchy sequence $\{x_i\}_{i \in \mathbb{N}}$ in $X$   with a central distance $r$.
 For every $y$ in $X$ and positive real number $\epsilon$, there exists a natural number $N$ such that for all   $i>N$,
$$-\epsilon <p(y,x_i)-p(y,a)<\epsilon.$$
I.e. $$\lim\limits_{i\to +\infty}p(y,x_i)=p(y,a).$$
\end{lem}
\textbf{Proof:} Let $a$ be the special limit of\ the Cauchy sequence $\{x_i\}_{i \in \mathbb{N}}$ in $X$   with a central distance $r$. From \Cref{D4b2} we have $$p(a,a)=r.$$
For every $y$ in $X$ and positive real number $\epsilon$, there exists a natural number $N_1$ such that for all $i>N_1$,
$$p(a,x_i) -r=p(a,x_i)-p(a,a)< \frac{\epsilon}{2}$$
i.e. $$p(a,x_i)< r+\frac{\epsilon}{2}$$
additionally, there exists a natural number $N_2$ such that for all $i>N_2$,
$$r-\frac{\epsilon}{2}<p(x_i,x_j)<r+\frac{\epsilon}{2}.$$
In particular for $i=j$, $$-p(x_i,x_i)<-r+\frac{\epsilon}{2}.$$
 Hence, using (p-inq) we get that for every positive real number $\epsilon$, there exists a natural number $N=\max\{N_1,N_2\}$, such that for all $i>N$,
$$p(y,x_i)-p(y,a)\le p(y,a)+p(a,x_i)-p(a,a)-p(y,a)$$
$$=p(a,x_i)-p(a,a)<\frac{\epsilon}{2}<\epsilon.$$
Additionally, by (p-inq)
$$p(y,a)\le p(y,x_i)+p(x_i,a)-p(x_i,x_i)$$
by (p-sym) $$ =p(y,x_i)+p(a,x_i)-p(x_i,x_i)$$
$$<p(y,x_i)+r+\frac{\epsilon}{2}-r+\frac{\epsilon}{2}$$
$$=p(y,x_i)+\epsilon.$$
Therefore, $$-\epsilon <p(y,x_i)-p(y,a). \hspace{6ex}\square$$

\begin{defn}\label{D4b3} Let $(X,p)$ be a partial metric space. We say that $(X,p)$ is \textbf{\uline{complete}}  if and only if every Cauchy sequence has a special limit in $X$.
\end{defn}
\indent We conclude this section by extending the definition of a  Cauchy pair to the partial metric case.
\begin{defn}\label{Db4} Let $\{x_i\}_{i \in \mathbb{N}}$ and $\{y_i\}_{i \in \mathbb{N}}$ be two  sequences in a partial metric space $(X,p)$. We say that $\{x_i\}_{i \in \mathbb{N}}$ and $\{y_i\}_{i \in \mathbb{N}}$ form a \textbf{\uline{Cauchy pair}} if and only if there exists a real number $r$ such that for every positive real number $\epsilon$ there exists a natural number $N$ where for all $i,j>N$,
$$r-\epsilon<\min \{p(x_i,x_i),p(y_j,y_j)\}\le p(x_i,y_j)<r+\epsilon.$$
We call $r$ the \textbf{\uline{central distance}} of the Cauchy pair $\{x_i\}_{i \in \mathbb{N}}$ and $\{y_i\}_{i \in \mathbb{N}}$.
\end{defn}
\begin{thm}{\textbf{(Cauchy pair term comparison):}}
\label{T4b1}
\\Let $\{x_i\}_{i \in \mathbb{N}}$ and $\{y_i\}_{i \in \mathbb{N}}$ be two sequences in a partial metric space $(X,p)$. The statements below are equivalent.
\\(a) $\{x_i\}_{i \in \mathbb{N}}$ and $\{y_i\}_{i \in \mathbb{N}}$ form a Cauchy pair with central distance $r$.
\\(b) There exists a real number $r$ where for every positive real number $\epsilon$ there exists a natural number $N$ such that for all  $i\ge j> N$,
$$r-\epsilon<\min \{p(x_i,x_i),p(y_j,y_j)\}\le p(x_i,y_j)<r+\epsilon.$$
(c) There exists a real number $r$ where for every positive real number $\epsilon$ there exists a natural number $N$ such that for all  $i\ge j> N$,
$$r-\epsilon<\min \{p(x_j,x_j),p(y_i,y_i)\}\le p(x_j,y_i)<r+\epsilon.$$
\end{thm}
\textbf{Proof:} It is clear that (a) is true if and only if (b) and (c) are true.
\\\indent(b)$\Rightarrow$ (c): For every positive real number $\epsilon$ there exists a natural number $N$ such that for all $i\ge j>N$,
$$r-\frac{\epsilon}{5}<p(x_i,x_i) \text{ , } r-\frac{\epsilon}{5} <p(y_j,y_j)\text{ and } p(x_i,y_j)<r+\frac{\epsilon}{5}$$
i.e.
$$-p(x_i,x_i)<-r+\frac{\epsilon}{5} \text{ , }-p(y_j,y_j)<-r+\frac{\epsilon}{5}\text{ and } p(x_i,y_j)<r+\frac{\epsilon}{5}.\hspace{6ex}(\nabla)$$
Hence, for every positive real number $\epsilon$  for all  $i\ge j> N$, from $(\nabla)$ and using (p-lbnd)
$$r-\epsilon<r-\frac{\epsilon}{5}<r-\epsilon<\min \{p(x_j,x_j),p(y_i,y_i)\}\le p(x_j,x_j)\le p(x_j,y_i)$$
using (p-inq) twice we get
$$\le p(x_j,y_j)+p(y_j,y_i)-p(y_j,y_j)$$
$$\le p(x_j,y_j)+p(y_j,x_i)+p(x_i,y_i)-p(x_i,x_i)-p(y_j,y_j)$$
by (p-sym)
$$=p(x_j,y_j)+p(x_i,y_j)+p(x_i,y_i)-p(x_i,x_i)-p(y_j,y_j)$$
by $(\nabla)$
$$<r+\frac{\epsilon}{5}+r+\frac{\epsilon}{5}+r+\frac{\epsilon}{5}-r+\frac{\epsilon}{5}-r+\frac{\epsilon}{5}=r+\epsilon.$$
Therefore, for every positive real number $\epsilon$ there exists a natural number $N$ such that for all $i\ge j>N$,
$$r-\epsilon<\min \{p(x_j,x_j),p(y_i,y_i)\}\le p(x_j,y_i)<r+\epsilon.$$
Similarly we can prove that  (c)$\Rightarrow$ (b).$\hspace{6ex}\square$
\begin{lem}\label{L4b5} Let $\{x_i\}_{i \in \mathbb{N}}$ and $\{y_i\}_{i \in \mathbb{N}}$ be a Cauchy pair with a central distance $r$ in a partial metric space $(X,p)$. Then $\{x_i\}_{i \in \mathbb{N}}$ and $\{y_i\}_{i \in \mathbb{N}}$ are both Cauchy sequences with central distance $r$. If $a$ is a (special) limit of $\{x_i\}_{i \in \mathbb{N}}$ then a is also a (special) limit of $\{y_i\}_{i \in \mathbb{N}}$.
\end{lem}
\textbf{Proof:} Let $\{x_i\}_{i \in \mathbb{N}}$ and $\{y_i\}_{i \in \mathbb{N}}$ be a Cauchy pair in $X$ with a central distance $r$. Then for every positive real number $\epsilon$ there exists a natural number $N_1$ such that for all $i,j>N_1$,
$$r-\frac{\epsilon}{3}<p(x_i,x_i)\le p(x_i,y_j)<r+\frac{\epsilon}{3}$$
and $$r-\frac{\epsilon}{3}<p(y_j,y_j)\le p(x_i,y_j)<r+\frac{\epsilon}{3}.$$
In particular, since the above is true, then for all $i,j>N_1$,
$$p(x_j,y_j)<r+\frac{\epsilon}{3}$$
and
$$-p(y_j,y_j)<-r+\frac{\epsilon}{3}.$$
Hence, by (p-lbnd)$$r-\epsilon<r-\frac{\epsilon}{3}<p(x_i,x_i)\le p(x_i,x_j)$$
by (p-inq) $$\le p(x_i,y_j)+p(y_j,x_j)-p(y_j,y_j)$$
by (p-sym) $$=p(x_i,y_j)+p(x_j,y_j)-p(y_j,y_j)$$
$$<r+\frac{\epsilon}{3}+r+\frac{\epsilon}{3}-r+\frac{\epsilon}{3}=r+\epsilon.$$
Therefore, $\{x_i\}_{i \in \mathbb{N}}$ (and similarly $\{y_i\}_{i \in \mathbb{N}}$) is a Cauchy sequence with central distance $r$.
Therefore for every positive real number $\epsilon$ there exists a natural number $N_2$, such that for all  $i,j>N_2$,
$$r-\frac{\epsilon}{3}<p(x_i,x_j)<r+\frac{\epsilon}{3}$$
in particular, for $i=j$,$$-p(x_i,x_i)<-r+\frac{\epsilon}{3}.$$
\\\indent Now assume that $a$ is a  limit of $\{x_i\}_{i \in \mathbb{N}}$.  By \Cref{L4b1}, for every positive real number $\epsilon$ there exists a natural number $N_3$ such that for all $i>N_3$,
$$p(a,x_i)-p(a,a)<\frac{\epsilon}{3}$$
i.e. $$ p(a,x_i)<p(a,a)+\frac{\epsilon}{3}.$$
Therefore, for every positive real number $\epsilon$ there exists a natural number $N=\max\{N_1,N_2, N_3\}$ such that for all $i,j>N$,
by (p-inq) $$p(a,y_j)-p(a,a)\le p(a,x_i)+p(x_i,y_j)-p(x_i,x_i)-p(a,a)$$
$$<p(a,a)+\frac{\epsilon}{3}+r+\frac{\epsilon}{3}-r+\frac{\epsilon}{3}-p(a,a)=\epsilon.$$
\indent The special limit case follows from the fact that $\{x_i\}_{i \in \mathbb{N}}$ and $\{y_i\}_{i \in \mathbb{N}}$ have the same central distance $r$ as shown above.   $\hspace{1ex} \square$
\section{Partial $n-\mathfrak{M}$etric Space}
\label{C4c}
\indent \indent In \cite{Ass20152}, and to retain the feel of \Cref{D4b1}, we generalized a Cauchy sequence to the partial $n-\mathfrak{M}$etric space in the manner below.
\begin{defn}\label{D4c1}  Let $\{x_i\}_{i \in \mathbb{N}}$ be a sequence in a partial $n-\mathfrak{M}$etric space $(X,P)$.
We say that the sequence $\{x_i\}_{i \in \mathbb{N}}$ is \uline{\textbf{Cauchy}} if and only if there exists a real number $r$ where for each positive number $\epsilon$ there exists a natural number $N$ such that for all $i_1,i_2,...,i_n>N$, $$-\epsilon <P(\langle x_{i_t}\rangle_{t=1}^n)-r< \epsilon.$$
We call $r$ the \textbf{\uline{central distance}} of $\{x_i\}_{i \in \mathbb{N}}$.
\end{defn}
\indent Although the above definition is a natural generalization from the partial metric case, it may seem to the reader that it is a condition that is difficult to check in practice. \Cref{T4c1} makes it much simpler to check if a sequence in a partial $n-\mathfrak{M}$etric space is Cauchy.
\begin{thm}{\textbf{(Cauchy Sequence Two Term Comparison):}}
\label{T4c1}
\\Let $\{x_i\}_{i \in \mathbb{N}}$ be a sequence in a partial $n-\mathfrak{M}$etric space $(X,P)$. Then the statements below are equivalent.
\\(a) $\{x_i\}_{i \in \mathbb{N}}$ is a Cauchy sequence with a central distance $r$.
\\(b) There exists a positive real number $r$ where for every positive real number $\epsilon$ there exists a natural number $N$ such that for all  $i,j> N$,
$$-\epsilon <P(\langle x_{i}\rangle^{n-1},x_j)-r< \epsilon.$$
(c) There exists a positive real number $r$ where for every positive real number $\epsilon$ there exists a natural number $N$ such that for all $i \ge j > N$,
$$-\epsilon <P(\langle x_{i}\rangle^{n-1},x_j)-r< \epsilon.$$
(d) There exists a positive real number $r$ where for every positive real number $\epsilon$ there exists a natural number $N$ such that for all   $i \ge j > N$,
$$-\epsilon <P(\langle x_{j}\rangle^{n-1},x_i)-r< \epsilon.$$
\end{thm}
\textbf{Proof:} (a)$\Rightarrow$ (b) and (b)$\Rightarrow$ (c) are trivial.
\\(c) $\Rightarrow$ (d):  For every positive real number $\epsilon$,  let $\epsilon'= \frac{\epsilon}{2n-3} \le \epsilon$. Then, there exists a natural number $N$ such that for all $i\ge j>N$,
$$-\epsilon' <P(\langle x_{i}\rangle^{n-1},x_j)-r< \epsilon'.$$
In particular, for $i=j$ $$-\epsilon'<P(\langle x_j \rangle^n)-r$$
and $$-P(\langle x_i \rangle^n)<-r+\epsilon'.$$
By ($P_n$-lbnd) for all $i \ge j>N$,
$$-\epsilon \le-\epsilon'< P(\langle x_j \rangle^n)-r\le 
P(\langle x_{j}\rangle^{n-1},x_i)-r$$
by \Cref{C2e3},
$$\le (n-1)P(\langle x_{i}\rangle^{n-1},x_j)-(n-2)P(\langle x_{i}\rangle^{n})-r$$
$$= (n-1)P(\langle x_{i}\rangle^{n-1},x_j)+(n-2)[-P(\langle x_{i}\rangle^{n})]-r$$
$$<(n-1)(r+\epsilon')+(n-2)(-r+\epsilon')-r$$
$$=(2n-3)\epsilon'=\epsilon.$$
Hence,
$$-\epsilon <P(\langle x_{j}\rangle^{n-1},x_i)-r< \epsilon.$$
(d)$\Rightarrow$ (a): For every positive real number $\epsilon$,  let $\epsilon'= \frac{\epsilon}{2n-3}\le\epsilon$. Then, there exists a natural number $N$ such that for all $i \ge j >N$,
$$-\epsilon '<P(\langle x_{j}\rangle^{n-1},x_i)-r< \epsilon'.$$
For all $i_1,i_2,...,i_n>N$, we will prove that 
$$-\epsilon<P(\langle x_{i_t}\rangle_{t=1}^{n})-r< \epsilon'.$$
Without loss of generality by ($P_n$-sym),  assume that $i_n \ge i_{n-1} \ge ..... \ge i_2 \ge i_1 > N$. Hence, from (d), for all $t\ge k$ in $\{1,...,n\}$ we get $$(\star)\begin{cases}P(\langle x_{i_k} \rangle^{n-1},x_{i_{t}})<r+\epsilon' &  \\P(\langle x_{i_k} \rangle^{n})<r+\epsilon' & \\  
-\epsilon' <P(\langle x_{i_t} \rangle^{n})-r &
\\-P(\langle x_{i_k} \rangle^{n})<-r+\epsilon' &  \\
\end{cases}$$
In particular, for $t=n$ in $(\star)$,
$$-\epsilon' < P(\langle x_{i_n}\rangle^{n})-r$$
by ($P_n$-lbnd)
$$\le P(\langle x_{i_n} \rangle^{n-1},x_{i_1})-r$$
by ($P_n$-sym) 
$$=P(\langle x_{i_n} \rangle^{n-2},x_{i_1},x_{i_n})-r$$
by \Cref{T2e1},
$$\le P(\langle x_{i_k} \rangle_{k=2}^{n-1},x_{i_1},x_{i_n})+\sum\limits_{k=2}^{n-1}[P(\langle x_{i_k} \rangle^{n-1},x_{i_n})-P(\langle x_{i_k} \rangle^{n})]-r$$
by ($P_n$-sym)
$$=P(\langle x_{i_k} \rangle_{k=1}^{n})+\sum\limits_{k=2}^{n-1}[P(\langle x_{i_k} \rangle^{n-1},x_{i_n})-P(\langle x_{i_k} \rangle^{n})]-r$$
by $(\star)$ and for $t= n$,
$$<P(\langle x_{i_k} \rangle_{k=1}^{n})+\sum\limits_{k=2}^{n-1}[(r+\epsilon')+(-r+\epsilon')]-r$$
$$=P(\langle x_{i_k} \rangle_{k=1}^{n})+\sum\limits_{k=2}^{n-1}2\epsilon'-r$$
$$=P(\langle x_{i_k} \rangle_{k=1}^{n})+(n-2)(2\epsilon')-r$$
$$=P(\langle x_{i_k} \rangle_{k=1}^{n})-r +(2n-4)\epsilon'$$
and, hence, $$-\epsilon'<P(\langle x_{i_k} \rangle_{k=1}^{n})-r +(2n-4)\epsilon'$$
therefore, 
$$P(\langle x_{i_k} \rangle_{k=1}^{n})-r >-(2n-3)\epsilon'=-\epsilon.$$
On the other hand  by \Cref{C2e1},
 $$P(\langle x_{i_t} \rangle_{t=1}^{n})-r\le P(\langle x_{i_1} \rangle^n) +\sum\limits_{t=2}^{n-1}[P(\langle x_{i_1} \rangle^{n-1},x_{i_t})-P(\langle x_{i_1} \rangle^{n})]-r$$
by $(\star)$ and for $k= 1$,
$$<r+\epsilon' +\sum\limits_{t=2}^{n-1}[(r+\epsilon')+(-r+\epsilon')]-r$$

$$=\epsilon'+\sum\limits_{t=2}^{n-1}2\epsilon'=\epsilon'+(n-2)(2\epsilon')$$
$$=(2n-3)\epsilon'=\epsilon.$$
Therefore
$$-\epsilon < P(\langle x_{i_t}\rangle_{t=1}^n)-r<\epsilon. \hspace{6ex}\square$$

\begin{lem}\label{L4c1} Let $\{x_i\}_{i \in \mathbb{N}}$ be a sequence in a partial $n-\mathfrak{M}$etric space $(X,P)$. A point $a$ in $X$ is a limit of $\{x_i\}_{i \in \mathbb{N}}$ if and only if for every positive real number $\epsilon$ there exists a natural number $N$ such that for all $i>N$,
$$P(\langle a\rangle^{n-1},x_i)-P(\langle a \rangle^n)< \epsilon.$$
\end{lem}
\textbf{Proof: }Let $\{x_i\}_{i \in \mathbb{N}}$ be a sequence in a partial $n-\mathfrak{M}$etric space $(X,P)$.
\\\indent$(\Rightarrow)$ Let $a$ be a limit of $\{x_i\}_{i \in \mathbb{N}}$ then for every positive real number $\epsilon$, $B_\epsilon^P(a)$ is an open neighborhood of $a$. Hence, there exists a positive natural number $N$ such that for all $i>N$,
$$x_i\in B_\epsilon^P(a).$$
Therefore, from \Cref{D3c1} 
$$P(\langle a\rangle^{n-1},x_i)-P(\langle a \rangle^n)< \epsilon.$$
\indent$(\Leftarrow)$ Assume that for each positive real number $\epsilon$ there exists a natural number $N$ such that for all $i>N$,
$$P(\langle a\rangle^{n-1},x_i)-P(\langle a \rangle^n)< \epsilon$$
and, hence, from \Cref{D3c1}
$$x_i \in B_\epsilon^P(a).$$
\indent For every  open neighborhood $U$ of $a$, from \Cref{D3c2} we know that their exists a positive real number $\epsilon$ such that 
$$B_\epsilon^P(a)\subseteq U.$$
Therefore, there exists a natural number $N$ such that for all $i>N$,
$$x_i \in B_\epsilon^P(a)\subseteq U. \hspace{6ex} \square$$ 
\indent The next theorem shows that partial $n-\mathfrak{M}$etrics possess a kind of upper semi-continuity property.
\begin{thm}{\textbf{(Upper Semi Continuity):}}
\label{T4c2}
\\ Let $\{x_i\}_{i \in \mathbb{N}}$ be a Cauchy sequence in a partial $n-\mathfrak{M}$etric space $(X,P)$ with  limit $a$ in $X$. Then for every $0\le q\le n-1$ , $\{b_k\}_{k=1}^q \subseteq X$ and positive real number $\epsilon$, there exists a natural number $N$ such that  for all $i_1,i_2,...,i_{n-q}>N$,
$$P(\langle x_{i_t}\rangle_{t=1}^{n-q},\langle b_k\rangle_{k=1}^q) <P(\langle a\rangle^{n-q},\langle b_k\rangle_{k=1}^q) +\epsilon.$$
\end{thm}
\textbf{Proof:} Let $0\le q\le n-1$, $\{b_k\}_{k=1}^q \subseteq X$ and positive real number $\epsilon$. Let $\epsilon'=\frac{\epsilon}{n-q}$. Then, by \Cref{L4c1}, there exists a natural number $N$ such that  for all   $i_1,i_2,...,i_{n-q}>N$, 
$$\hspace{2ex}P(\langle a\rangle^{n-1},x_{i_t})-P(\langle a \rangle^n)< \epsilon'.$$
Hence,  without loss of generality by ($P_n$-sym) we may assume that $i_{n-q} \ge i_{n-q-1} \ge ...\ge i_2 \ge i_1 > N$. Then by \Cref{T2e1} , 
$$P(\langle x_{i_t}\rangle_{t=1}^{n-q},\langle b_k\rangle_{k=1}^q)<P(\langle a\rangle^{n-q},\langle b_k\rangle_{k=1}^q)+\sum\limits_{t=1}^{n-q}[P(\langle a\rangle^{n-1},x_{i_t})-P(\langle a\rangle^{n})]$$
$$<P(\langle a\rangle^{n-q},\langle b_k\rangle_{k=1}^q)+\sum\limits_{t=1}^{n-q}\epsilon'$$
$$=P(\langle a\rangle^{n-q},\langle b_k\rangle_{k=1}^q)+(n-q)\epsilon'=P(\langle a\rangle^{n-q},\langle b_k\rangle_{k=1}^q)+\epsilon. \hspace{6ex} \square$$
\begin{cor}\label{C4c1} Let $\{x_i\}_{i \in \mathbb{N}}$ be a Cauchy sequence in a partial $n-\mathfrak{M}$etric space $(X,P)$ with a limit $a $ in $X$. Then for every positive real number $\epsilon$ there exists a natural number $N$ such that for all $i_1,i_2,...,i_n>N$ and  all $0\le q \le n$ the statements below hold true.
\\(a) $P(\langle x_{i_t}\rangle_{t=1}^{n-q},\langle a\rangle^q) < P(\langle a\rangle^n)+\epsilon.$
\\(b) $P(\langle x_{i_t}\rangle_{t=1}^{n}) < P(\langle a\rangle^n)+\epsilon.$
\\(c)  $P(\langle x_{i_t}\rangle_{t=1}^{n-1},a) < P(\langle a\rangle^n)+\epsilon.$
\end{cor}
\indent \indent The above Corollary is trivial to prove using \Cref{T4c2} while varying $1\le q\le n-1$ and taking  for all $k$, $b_k=a$. The case where $q=n$ is trivial since
$$P(\langle a\rangle^n)<P(\langle a\rangle^n)+\epsilon$$ for any positive real number $\epsilon$. 
\\\indent As in \Cref{E4b1}, a limit of a Cauchy sequence need not be unique.
\begin{defn}\label{D4c2} Let $(X,p)$ be a partial $n-\mathfrak{M}$etric space. Consider a Cauchy sequence $\{x_i\}_{i \in \mathbb{N}}$ in $X$ with a central distance $r$.  A point $a$ in $X$ is called a \textbf{\uline{special limit}} of  $\{x_i\}_{i \in \mathbb{N}}$ if and only if $a$ is a limit of $\{x_i\}_{i \in \mathbb{N}}$ and $r=P(\langle a\rangle^n)$.
\end{defn}
\indent \indent As in the case of a partial metric space a special limit is unique.
\begin{thm}{\textbf{(Uniqueness of Special Limits):}}
\label{T4c3}
\\ Let $\{x_i\}_{i \in \mathbb{N}}$ be a Cauchy sequence in a partial $n-\mathfrak{M}$etric space $(X,P)$. If $\{x_i\}_{i \in \mathbb{N}}$ has a special limit in $X$ then that special limit is unique.
\end{thm}
\textbf{Proof:} Consider the Cauchy sequence $\{x_i\}_{i \in \mathbb{N}}$ in $X$ with a central distance $r$.
If $a$ and $b$ are both special limits of  $\{x_i\}_{i \in \mathbb{N}}$ , then by \Cref{D4c2}
$$P(\langle a\rangle^n)=r=P(\langle b\rangle^n).$$
Furthermore, from \Cref{D4c1}, for every positive real number $\epsilon$ there exists a natural number $N_1$ such that for all $i>N_1$,
$$-\frac{\epsilon}{3} <P(\langle x_i\rangle^n)-r,$$
i.e.
$$-P(\langle x_i\rangle^n)<-r+\frac{\epsilon}{3}.$$
The special limit $a$ of $\{x_i\}_{i \in \mathbb{N}}$ is also a limit
of $\{x_i\}_{i \in \mathbb{N}}$. Hence, by \Cref{C4c1} for every positive real number $\epsilon$ there exists a natural number $N_2$ such that for all $i>N_2$,
$$P(\langle a\rangle^{n-1},x_i)-P(\langle a\rangle^n)<\frac{\epsilon}{3}$$
i.e.
$$P(\langle a\rangle^{n-1},x_i)<r+\frac{\epsilon}{3}.$$
The special limit $b$ of $\{x_i\}_{i \in \mathbb{N}}$ is also a limit
of $\{x_i\}_{i \in \mathbb{N}}$. Hence, by \Cref{C4c2}  there exists a natural number $N_3$ such that for all $i>N_3$,
$$P(\langle x_i\rangle^{n-1},b) < P(\langle b\rangle^n)+\frac{\epsilon}{3},$$
i.e.
$$P(\langle x_{i}\rangle^{n-1},b) < r+\frac{\epsilon}{3}.$$
Hence, using ($P_n$-lbnd) we get that for every positive real number $\epsilon$ there exists a natural number 
\\$N=\max\{N_1,N_2,N_3\}$ such that for all $i>N$,
$$P(\langle a\rangle^n)\le P(\langle a\rangle^{n-1},b)$$
by ($P_n$-inq)
$$\le P(\langle a\rangle^{n-1},x_i)+P(\langle x_i\rangle^{n-1},b)-P(\langle x_i\rangle^{n})$$
$$<r+\frac{\epsilon}{3}+r+\frac{\epsilon}{3}-r+\frac{\epsilon}{3}$$
$$=r+\epsilon=P(\langle a\rangle^n)+\epsilon$$
therefore, $$P(\langle a\rangle^n)=P(\langle a\rangle^{n-1},b).$$
Similarly $P(\langle b\rangle^n)=P(\langle b\rangle^{n-1},a)$ and, hence, by ($P_n$-sep)
$$a=b.\hspace{6ex} \square$$
\begin{thm}{\textbf{(Continuity of Special Limits):}}
\label{T4c4}
\\ Let $\{x_i\}_{i \in \mathbb{N}}$ be a Cauchy sequence in a partial $n-\mathfrak{M}$etric space $(X,P)$ with special limit $a$ in $X$. Then for every $0\le q\le n-1$ ,  $\{b_k\}_{k=1}^q \subseteq X$ and positive real number $\epsilon$, there exists a natural number $N$ such that  for all $i_1,i_2,...,i_{n-q}>N$,
$$P(\langle a\rangle^{n-q},\langle b_k\rangle_{k=1}^q) -\epsilon <P(\langle x_{i_t}\rangle_{t=1}^{n-q},\langle b_k\rangle_{k=1}^q) <P(\langle a\rangle^{n-q},\langle b_k\rangle_{k=1}^q) +\epsilon.$$
\end{thm}
\textbf{Proof:} Let $r$ be the central distance of   $\{x_i\}_{i \in \mathbb{N}}$. By \Cref{D4c2}
$$P(\langle a \rangle ^n)=r.$$
Let $0\le q\le n-1$ ,   $\{b_k\}_{k=1}^q \subseteq X$, and  positive real number $\epsilon$. Let $\epsilon'=\frac{\epsilon}{2(n-q)}$.  By \Cref{D4c2}, $a$ is also a limit of $\{x_i\}_{i \in \mathbb{N}}$ and, hence, there exists a natural number $N_1$ such that for all $i_t>N_1$,
$$P(\langle x_{i_t}\rangle^{n-1},a)<P(\langle a \rangle ^n)+\epsilon',$$
i.e. $$P(\langle x_{i_t} \rangle ^{n-1},a)<r+\epsilon'.$$
By \Cref{D4c1}, there exists a natural number $N_2$ such that for all $i_t > N_2$,
$$-\epsilon'<P(\langle x_{i_t} \rangle^n)-r,$$
i.e. $$-P(\langle x_{i_t} \rangle^n)<-r+\epsilon'.$$
Hence, using \Cref{T2e1} there exists a natural number $N=\max\{N_1,N_2\}$ such that for all 
\\$i_1,i_2,...,i_{n-q}>N$,
$$P(\langle a\rangle^{n-q},\langle b_k\rangle_{k=1}^q) -\epsilon$$
$$\le P(\langle x_{i_t}\rangle_{t=1}^{n-q},\langle b_k\rangle_{k=1}^q)+\sum\limits_{t=1}^{n-q}[P(\langle x_{i_t}\rangle^{n-1},a)-P(\langle x_{i_t}\rangle^{n})]-\epsilon$$
$$<P(\langle x_{i_t}\rangle_{t=1}^{n-q},\langle b_k\rangle_{k=1}^q)+\sum\limits_{t=1}^{n-q}[r+\epsilon'-r+\epsilon']-\epsilon$$
$$=P(\langle x_{i_t}\rangle_{t=1}^{n-q},\langle b_k\rangle_{k=1}^q)+\sum\limits_{t=1}^{n-q}2\epsilon'-\epsilon$$
$$=P(\langle x_{i_t}\rangle_{t=1}^{n-q},\langle b_k\rangle_{k=1}^q)+(n-q)(2\epsilon')-\epsilon$$
$$=P(\langle x_{i_t}\rangle_{t=1}^{n-q},\langle b_k\rangle_{k=1}^q)+\epsilon-\epsilon=P(\langle x_{i_t}\rangle_{t=1}^{n-q},\langle b_k\rangle_{k=1}^q).$$
\indent \indent The right side of the inequality was already proved in \Cref{T4c2}. $\hspace{6ex} \square$
\begin{cor}\label{C4c2} Let $\{x_i\}_{i \in \mathbb{N}}$ be a Cauchy sequence in a partial $n-\mathfrak{M}$etric space $(X,P)$ with special limit $a$ in $X$. Then for every positive real number $\epsilon$ there exists a natural number $N$ such that for all $i_1,i_2,...,i_n>N$  and $0\le q \le n$ the statements below hold true.
\\(a) $-\epsilon<P(\langle x_{i_t}\rangle_{t=1}^{n-q},\langle a\rangle^q) -P(\langle a\rangle^n)<\epsilon.$
\\(b) $-\epsilon<P(\langle x_{i_t}\rangle_{t=1}^{n}) -P(\langle a\rangle^n)<\epsilon.$
\\(c)  $-\epsilon<P(\langle x_{i_t}\rangle_{t=1}^{n-1},a) - P(\langle a\rangle^n)<\epsilon.$
\end{cor}
\indent The above Corollary is trivial to prove using \Cref{T4c4} while varying $1\le q\le n-1$ and taking  for all $k$, $b_k=a$. The case where $q=n$ is trivial.
\begin{defn}\label{D4c3} Let $(X,P)$ be a partial $n-\mathfrak{M}$etric space. We say that $(X,P)$ is \textbf{\uline{complete}} if and only if  every Cauchy sequence has a special limit in $X$.
\end{defn}
\begin{defn}\label{D4c4} Let $\{x_i\}_{i \in \mathbb{N}}$ and $\{y_i\}_{i \in \mathbb{N}}$ be two sequences in a partial $n-\mathfrak{M}$etric space $(X,P)$. We say that $\{x_i\}_{i \in \mathbb{N}}$ and $\{y_i\}_{i \in \mathbb{N}}$ form a \textbf{\uline{Cauchy pair}} if and only if
there exists a real number $r$ such that for every positive real number $\epsilon$ there exists a natural number $N$ such that for all $i,j>N$,
$$r-\epsilon<\min \{P(\langle x_{i}\rangle^n),P(\langle y_{i}\rangle^n)\}\le P(\langle x_{i}\rangle^{n-1},y_j)<r+\epsilon.$$
We call $r$ the \textbf{\uline{central distance}} of the Cauchy pair $\{x_i\}_{i \in \mathbb{N}}$ and $\{y_i\}_{i \in \mathbb{N}}$.
\end{defn}
\begin{thm}{\textbf{(Cauchy pair term comparison):}}
\label{T4c5}
\\Let $\{x_i\}_{i \in \mathbb{N}}$ and $\{y_i\}_{i \in \mathbb{N}}$ be two sequences in a partial $n-\mathfrak{M}$etric space $(X,P)$. The statements below are equivalent.
\\(a) $\{x_i\}_{i \in \mathbb{N}}$ and $\{y_i\}_{i \in \mathbb{N}}$ form a Cauchy pair with central distance $r$.
\\(b) There exists a real number $r$ where for every positive real number $\epsilon$ there exists a natural number $N$ such that for all  $i, j> N$,
$$r-\epsilon<\min \{P(\langle x_{i}\rangle^n),P(\langle y_{i}\rangle^n)\}\le P(\langle y_{j}\rangle^{n-1},x_i)<r+\epsilon.$$
(c) There exists a real number $r$ where for every positive real number $\epsilon$ there exists a natural number $N$ such that for all  $i\ge j> N$,
$$r-\epsilon<\min \{P(\langle x_{i}\rangle^n),P(\langle y_{i}\rangle^n)\}\le P(\langle x_{i}\rangle^{n-1},y_j)<r+\epsilon.$$
(d) There exists a real number $r$ where for every positive real number $\epsilon$ there exists a natural number $N$ such that for all  $i\ge j> N$,
$$r-\epsilon<\min \{P(\langle x_{i}\rangle^n),P(\langle y_{i}\rangle^n)\}\le P(\langle x_{j}\rangle^{n-1},y_i)<r+\epsilon.$$
\end{thm}
\textbf{Proof:} It is clear that (a) is true if and only if (c) and (d) are true. Hence, it will be enough to prove that (a) is equivalent to (b) and (c) is equivalent to (d).
\\\indent (a)$\Rightarrow$(b): Assume that  $\{x_i\}_{i \in \mathbb{N}}$ and $\{y_i\}_{i \in \mathbb{N}}$ form a Cauchy pair with central distance $r$. For every positive real number $\epsilon$, let $\epsilon'=\frac{\epsilon}{2n-3}\le \epsilon$. Hence, there exists a natural number $N$ such that for all  $i, j> N$,
$$r-\epsilon'<\min \{P(\langle x_{i}\rangle^n),P(\langle y_{i}\rangle^n)\}\le P(\langle x_{i}\rangle^{n-1},y_j)<r+\epsilon'$$
in particular,
$$-P(\langle x_{i}\rangle^n)<-r+\epsilon'.$$
Therefore,
$$r-\epsilon\le r-\epsilon'\le \min \{P(\langle x_{i}\rangle^n),P(\langle y_{i}\rangle^n)\}$$
from ($P_n$-lbnd) 
$$\le P(\langle y_{i}\rangle^n)\le P(\langle y_{j}\rangle^{n-1},x_i)$$
by \Cref{C2e3}
$$\le (n-1)P(\langle x_{i}\rangle^{n-1},y_j)-(n-2)P(\langle x_{i}\rangle^{n})$$
$$<(n-1)(r+\epsilon') +(n-2)(-r+\epsilon')=r+(2n-3)\epsilon'=r+\epsilon.$$
The proof that (b)$\Rightarrow$ (b) is similar.
\\\indent (c)$\Rightarrow$ (d): For every positive real number $\epsilon$, let $\epsilon'=\frac{\epsilon}{2n+1}\le \epsilon$. Hence, there exists a natural number $N$ such that for all  $i\ge j> N$,
$$r-\epsilon'<\min \{P(\langle x_{i}\rangle^n),P(\langle y_{i}\rangle^n)\}\le P(\langle x_{i}\rangle^{n-1},y_j)<r+\epsilon'$$
in particular,
$$-P(\langle x_{i}\rangle^n)<-r+\epsilon' \text{ , } -P(\langle y_{j}\rangle^n)<-r+\epsilon',$$
$$P(\langle x_{i}\rangle^{n-1},y_i)<r+\epsilon' \text{ and } P(\langle x_{j}\rangle^{n-1},y_j)<r+\epsilon'.$$
Therefore,
$$r-\epsilon\le r-\epsilon'\le \min \{P(\langle x_{i}\rangle^n),P(\langle y_{i}\rangle^n)\}$$
from ($P_n$-lbnd)
$$\le P(\langle x_{j}\rangle^n) \le P(\langle x_{j}\rangle^{n-1},y_i)$$
by using ($P_n$-inq) twice we get
$$\le P(\langle x_{j}\rangle^{n-1},y_j)+P(\langle y_{j}\rangle^{n-1},y_i)-P(\langle y_{j}\rangle^{n})$$
$$\le P(\langle x_{j}\rangle^{n-1},y_j)+P(\langle y_{j}\rangle^{n-1},x_i)+P(\langle x_{i}\rangle^{n-1},y_i)-P(\langle x_{i}\rangle^n)-P(\langle y_{j}\rangle^{n})$$
by \Cref{C2e3}
$$\le P(\langle x_{j}\rangle^{n-1},y_j)+(n-1)P(\langle x_{i}\rangle^{n-1},y_j)-(n-2)P(\langle x_{i}\rangle^n)+P(\langle x_{i}\rangle^{n-1},y_i)-P(\langle x_{i}\rangle^n)-P(\langle y_{j}\rangle^{n})$$
$$< (r+\epsilon')+(n-1)(r+\epsilon')+(n-2)(-r+\epsilon')+(r+\epsilon')+(-r+\epsilon')+(-r+\epsilon')$$
$$=r+(2n+1)\epsilon'=r+\epsilon.$$  
Similarly we can prove that (d)$\Rightarrow$(c).$\hspace{6ex}\square$
\begin{lem}\label{L4c2} Let $\{x_i\}_{i \in \mathbb{N}}$ and $\{y_i\}_{i \in \mathbb{N}}$ be a Cauchy pair with a central distance $r$ in a  partial $n-\mathfrak{M}$etric space $(X,P)$. Then $\{x_i\}_{i \in \mathbb{N}}$ and $\{y_i\}_{i \in \mathbb{N}}$ are both Cauchy sequences with central distance $r$. If $a$ is a (special) limit of $\{x_i\}_{i \in \mathbb{N}}$ then $a$ is a (special) limit of $\{y_i\}_{i \in \mathbb{N}}$.
\end{lem}
\textbf{Proof:} Let $\{x_i\}_{i \in \mathbb{N}}$ and $\{y_i\}_{i \in \mathbb{N}}$ be a Cauchy pair in $X$ with a central distance $r$. Then   for every positive real number $\epsilon$ there exists a natural number $N_1$ such that for all $i,j>N_1$,
$$r-\frac{\epsilon}{3}<\min \{P(\langle x_{i}\rangle^n),P(\langle y_{i}\rangle^n)\}\le P(\langle x_{i}\rangle^{n-1},y_j)<r+\frac{\epsilon}{3}$$
and
$$r-\frac{\epsilon}{3}<\min \{P(\langle x_{i}\rangle^n),P(\langle y_{i}\rangle^n)\}\le P(\langle y_{i}\rangle^{n-1},x_j)<r+\frac{\epsilon}{3}.$$
Hence, $$-P(\langle x_{i}\rangle^n)<-r+\frac{\epsilon}{3}.$$
By ($P_n$-lbnd)
$$r-\epsilon<r-\frac{\epsilon}{3}<P(\langle y_{i}\rangle^n)\le P(\langle y_{i}\rangle^{n-1},y_j)$$
by ($P_n$-inq)
$$\le P(\langle y_{i}\rangle^{n-1},x_j)+P(\langle x_{i}\rangle^{n-1},y_j)-P(\langle x_{i}\rangle^n)$$
$$<r+\frac{\epsilon}{3}+r+\frac{\epsilon}{3}-r+\frac{\epsilon}{3}=r+\epsilon.$$
Hence, by \Cref{T4c1}, $\{y_i\}_{i \in \mathbb{N}}$ (and similarly $\{x_i\}_{i \in \mathbb{N}}$ ) is a Cauchy sequence with central distance $r$.
\indent Now assume that $a$ is a limit of $\{x_i\}_{i \in \mathbb{N}}$. By \Cref{L4c1}, for every positive real number $\epsilon$, there exists
a natural number $N_2$ such that for all $i > N_2$,
$$P(\langle a \rangle^{n-1},x_i)-P(\langle a \rangle^n)<\frac{\epsilon}{3}.$$
Therefore, for every positive number $\epsilon$, there exists a natural number $N=\max\{N_1,N_2\}$ such that for all $i>N$,
$$ P(\langle a \rangle^{n-1},y_i)-P(\langle a \rangle^n)$$
by ($P_n$-inq)
$$\le P(\langle a \rangle^{n-1},x_i)+P(\langle x_{i} \rangle^{n-1},y_i)-P(\langle x_{i} \rangle^n)-P(\langle a \rangle^n)$$
$$=P(\langle a \rangle^{n-1},x_i)-P(\langle a \rangle^n)+P(\langle x_{i} \rangle^{n-1},y_i)-P(\langle x_{i} \rangle^n)$$
$$<\frac{\epsilon}{3}+r+\frac{\epsilon}{3}-r+\frac{\epsilon}{3}=\epsilon.$$
\indent The special limit case follows from the fact that $\{x_i\}_{i \in \mathbb{N}}$ and $\{y_i\}_{i \in \mathbb{N}}$ have the same central distance $r$ as shown above. $\hspace{6ex} \square$
\chapter{Cauchy Functions}
\label{C5}
\setcounter{thm}{0}
\setcounter{defn}{0}
It all started in \textbf{1922} with  Banach \cite{Ban1922}. Given a metric space $(X,d)$ and a function $f:X \to X$, Banach gave contracting criteria on $f$ allowing him to generate from that function a Cauchy sequence. He then proved that the limit of this sequence  is a fixed point. His fixed point theorem was generalized in many ways, but all generalizations had the same flow.
\\Step 1: Give criteria for the function to generate a Cauchy sequence.
\\Step 2: Make sure the limit of that Cauchy sequence exists.
\\Step 3: Give criteria on that function that leads to sequential continuity on the limit of the Cauchy sequence.
\\ This ensures the existence of  a fixed point is found.
\\\indent The most notable generalization was given in \textbf{1962} by Edelstein \cite{Ede1961}, who restricted his attention to continuous functions contractive on an orbit. In \textbf{1977},  Alber and Guerre-Delabriere \cite{Alb1997}  generalized  contraction to what they called weak $\varphi$-contraction.
Their work was restricted to maps on Hilbert spaces. In \textbf{2001}, Rhoades \cite{Rho2001} showed that the results in \cite{Alb1997} still hold in any Banach space.
\\\indent In  \Cref{C5} we will investigate some contractive criteria on a function that enables it to generate a Cauchy sequence in its domain.
\begin{defn}\label{D5-1} Let $(X,\mathcal{T})$ be a topological space with $x_o$ in $X$. Let $f: X\to X$ be a function on $X$. Denote $f^0(x_o)=x_0$, $f^1(x_o)=f(x_o)$ and inductively $f^{n+1}(x_o)=f(f^n(x_o))$. The \uline{\textbf{orbit of $f$ at $x_o$}} is the sequence $\{f^i(x_o)\}_{i\in \mathbb{N}}$.
\end{defn}
 We will state the definitions below on a partial $n-\mathfrak{M}$etric space knowing that it includes all other cases discussed in this thesis.
\begin{defn}\label{D5-2} Let  $(X,P)$ be a partial $n-\mathfrak{M}$etric space with $x_o$ in $X$ and suppose $f:X\to X$ is a function on $X$. We say that $f$ is a \uline{\textbf{Cauchy function} at $x_o$} if and only if $ \{f^i(x_o)\}_{i \in \mathbb{N}}$ is a Cauchy sequence. We say that $f$ is a \uline{\textbf{Cauchy function}} if and only if for every $x$ in $X$, $\{f^i(x)\}_{i \in \mathbb{N}}$ is a Cauchy sequence.
\end{defn}
\begin{defn}\label{D5-3} Let $(X,P)$ be a partial $n-\mathfrak{M}$etric space with $x_o$ and $y_o$ in $X$ and suppose  $f:X\to X$ and $g:X \to X$ are two functions on $X$. We say that $f$ and $g$ form a  \uline{\textbf{Cauchy pair over $(x_o,y_o)$}} if and only if $\{f^i(x_o)\}_{i\in\mathbb{N}}$ and $\{g^i(y_o)\}_{i\in\mathbb{N}}$ form a Cauchy pair. I.e. there exists a real number $r$ such that
$$\lim_{i,j\to+\infty} P(\langle x_{i}\rangle^{n-1},y_j)=r.$$
\end{defn}
\section{Metric Space}
\label{C5a}
The material in this section
will be presented more completely in the more general  case in \Cref{C5b} as partial metric spaces.
\begin{defn}\label{D5a2} (Edelstein \cite{Ede1961}) Let $(X,d)$ be a metric space with $x_o$ in $X$. Let $f:X\to X$ be a function on $X$  and $0<c<1$ be a real number. We say that $f$ is an \uline{\textbf{orbital $c_0-$contraction at $x_o$}} (or $f$ is \textbf{\uline{orbitally $c_0-$contractive at $x_o$}})  if and only if  for all natural numbers  $i$,
$$d(f^{i+1}(x_o),f^{i}(x_0) \le c d(f^{i}(x_o),f^{i-1}(x_o)).$$
We say that $f$ is an \textbf{orbital $c_0-$contraction} (or $f$ is \textbf{\uline{orbitally $c_0-$contractive}}) if and only if for every $x$ in $X$, $f$ is an orbital contraction at $x$.
\end{defn}
In the partial metric case, the central distance $r$ of the Cauchy sequence obtained need not be $0$. That is why the notation \textbf{\uline{orbital $c_0-$contraction }}was presented to allow the use of the term \textbf{\uline{orbital $c_r-$contraction}} in the more general case.\begin{lem}\label{L5a1} (Edelstein \cite{Ede1961}) Let $(X,d)$ be a metric space with $x_o$ in $X$. Let $f:X\to X$ be a function on $X$. If $f$ is an orbital  contraction at $x_o$ then $f$ is a Cauchy function at $x_o$.
\end{lem}
We build upon the work of Rhoades \cite{Rho2001} and Edelstein \cite{Ede1961,Ede1962,Ede1964} to present the definition below.    
\begin{defn}\label{D5a3} Let $(X,d)$ be a metric space with $x_o$ in $X$. Let $f:X\to X$ be a function on $X$ and $\varphi:[0,+\infty )\subset \mathbb{R}\to [0, + \infty)$ be a non-decreasing function such that
$$\varphi(t)=0 \text{ if and only if } t=0.$$
 We say that $f$ is an \uline{\textbf{orbital $\varphi_0-$contraction at $x_o$}} (or $f$ is \uline{\textbf{orbitally $\varphi_0-$contractive at $x_o$}}) if and only if for all  $i$ and $j$,
$$d(f^{i+1}(x_o),f^{j+1}(x_o))\le d(f^i(x_o),f^j(x_o))-\varphi (d(f^i(x_o),f^j(x_o)).$$
We say that $f$ is an \uline{\textbf{orbital $\varphi_0-$contraction}} (or $f$ is \uline{\textbf{orbitally $\varphi_0-$contractive}}) if and only if for every $x$ in $X$, $f$ is an orbital $\varphi_0-$contraction at $x$.
\end{defn}
\begin{rmk}\label{R5a1} The reader should note that any orbital $c_0-$contraction is an orbital $\varphi_0-$contraction by taking
$$\varphi(t)=(1-c)t.$$
\end{rmk}
\begin{lem}\label{L5a2} Let $(X,d)$ be a metric space with $x_o$ in $X$. Let $f:X\to X$ be a function on $X$. If $f$ is an orbital $\varphi_0-$contraction at $x_o$ then $f$ is a Cauchy function at $x_o$.
\end{lem}
\indent While searching the literature for common fixed point theorems, those we found required the two functions to commute and the space was required to have certain conditions in addition to being complete \cite{Fis1979}. We wanted to present a common fixed point theorem that relies on a contraction criteria.  In \textbf{2009}, Zhanga and Song \cite{Zha2009} presented us with just that, however, their contraction is defined over the whole space rather than an orbit. 
\begin{defn}\label{D5a4} Let $(X,d)$ be a metric space with $x_o$ and $y_o$ in $X$. Let $f:X\to X$ and $g:X \to X$ be two functions on $X$ and $0<c<1$ be a real number. We say that $f$ and $g$ are \textbf{\uline{$f-$pairwise $c_0-$contractive over $(x_o,y_o)$}}  if and only if  for all natural numbers  $i$,
$$\begin{cases}d(f^{i+1}(x_0),g^{i}(y_o))\le cd(f^{i}(x_o),g^{i-1}(y_o))
\\d(f^{i}(x_0),g^{i}(y_o))\le cd(f^{i-1}(x_o),g^{i-1}(y_o))
\end{cases}$$
\end{defn}
In fact \Cref{D5a4} is more general and much easier to check than the following possible alternative definition,
$$d(f^{i+1}(x_0),g^{j+1}(y_o))\le cd(f^{i}(x_o),g^{j}(y_o)).$$
\begin{thm}{\textbf{(Cauchy $f-$Pairwise $c_0-$contractive):}}
\label{T5a1}
\\Let $(X,d)$ be a metric space with $x_o$ and $y_o$ in $X$. Let $f:X\to X$ and $g:X \to X$ be two functions on $X$. If $f$ and $g$ are $f-$pairwise $c_0-$contractive over $(x_o,y_o)$ then $f$ and $g$ form a Cauchy pair over $(x_o,y_o)$.
\end{thm}
\begin{defn}\label{D5a6} Let $(X,l)$ and $(Y,d)$ be two metric spaces. Let $f:X\to Y$ and $g:X \to Y$ two functions on $X$. Let $A\ge 0$ and  $0<c<1$ be two real numbers. We say that $f$ and $g$ are \textbf{\uline{mutually $c_0-$contractive}} if and only if  for each $x$ in $X$ we can find an element $z$ in $X$  such that 
$$d(f(z),g(z)) \le cd(f(x),g(x))$$
and
$$l(x,z)\le Ad(f(x),g(x)).$$
\end{defn}
\Cref{T5a2} is a special case of \Cref{T5b2}. We have chosen to include the proof because, as far as we have been able to determine, the result is new and the proof increases in complexity in the partial metric case.
\begin{thm}{\textbf{(Cauchy Mutually $c_0-$Contractive):}}
\label{T5a2}
\\Let $(X,l)$ and $(Y,d)$ be two metric spaces. Let $f:X\to Y$ and $g:X \to Y$ be two functions on $X$. If $f$ and $g$ are mutually $c_0-$contractive then there exists a Cauchy sequence $\{x_i\}_{i\in \mathbb{N}}$ in $X$ such that $\{f(x_i)\}_{n\in \mathbb{N}}$ and $\{g(x_i)\}_{n\in \mathbb{N}}$ form a Cauchy pair in $Y$.
\end{thm}
\textbf{Proof:} Since $f$ and $g$ are mutually contractive then there exist two real numbers $0<c<1$ and $A\ge0$ such that for each $x$ in $X$, there exists a $z$ in $X$ such that  
$$d(f(z),g(z)) \le cd(f(x),g(x))$$
and
$$l(x,z)\le Ad(f(x),g(x)).$$
Let us take an arbitrary element $x_1$ in $X$ then there exists an element $x_2$ in $X$ such that
$$d(f(x_2),g(x_2))\le cd(f(x_1),g(x_1)) \text{ and } l(x_1,x_2)\le Ad(f(x_1),g(x_1)).$$ There exists an element $x_3$ in $X$ such that 
$$d(f(x_3),g(x_3))\le cd(f(x_2),g(x_2)) \text{ and } l(x_2,x_3)\le Ad(f(x_2),g(x_2)).$$ We continue the above process to generate the sequence $\{x_i\}_{i\in \mathbb{N}}$ such that for all  $i$,
$$d(f(x_{i+1}),g(x_{i+1}))\le cd(f(x_i),g(x_i)) \text{ and } l(x_i,x_{i+1})\le Ad(f(x_i),g(x_i)).\hspace{6ex}(\ddot\diamond)$$
\indent We now prove that $\{x_i\}_{i\in \mathbb{N}}$ is a Cauchy sequence with two steps.
\\\indent Step 1: For all  $i$ let $t_i=  d(f(x_i),g(x_i))$. Then from $(\ddot\diamond)$
$$t_{i+1}=d(f(x_{i+1}),g(x_{i+1}))\le cd(f(x_i),g(x_i))\le c^2d(f(x_{i-1}),g(x_{i-1}))$$
and, hence, by induction
$$t_{i+1}\le  c^id(f(x_1),g(x_1))=c^i t_1. \hspace{6ex}(\check\diamond)$$
Since $0<c<1$ then for every positive real number $\epsilon$ there exists a natural number $N'$ such that $c^{N'-1}t_1<\epsilon$ and, hence, for all $i>N'$,
$$d(f(x_i),g(x_i))=t_i\le c^{i-1}t_1<c^{N'-1}t_1<\epsilon.$$
\\\indent Step 2: For all  $j>i$ and by repeatedly using (d-inq) we get,
$$l(x_i,x_j)\le l(x_i,x_{i+1})+l(x_{i+1},x_{i+2})+.......l(x_{j-1},x_j)$$
from \Cref{D5a6}
$$\le Ad(f(x_{i}),g(x_{i}))+Ad(f(x_{i+1}),g(x_{i+1}))+.....+Ad(f(x_{j-1}),g(x_{j-1}))$$
$$=At_i+At_{i+1}+......+At_{j-1}$$
by $(\check\diamond)$
$$\le Ac^{i-1}t_1+Ac^it_1+......+Ac^{j-2}t_1=At_1\sum\limits_{k=i-1}^{j-1}c^k=At_1c^{i-1}\sum\limits_{k=i-1}^{j-i-1}c^k.$$
Hence, knowing that $At_1 \ge 0$ we get
$$l(x_i,x_j)\le At_1c^{i-1}\sum\limits_{k=0}^{j-i-1}c^k \le At_1\sum\limits_{k=i-1}^{+\infty}c^k$$
by the geometric series formula
$$=At_1c^{i-1}\sum\limits_{k=0}^{+\infty}c^k=\frac{At_1}{1-c}c^{i-1}.$$
Therefore, for every positive number $\epsilon$ there exists a natural number $N$ where  $\frac{At_1}{1-c}c^{N-1}<\epsilon$ and, hence, for all  $j\ge i>N$,
$$l(x_j,x_i)\le\frac{At_1}{1-c}c^{i-1}< \frac{At_1}{1-c}c^{N-1}<\epsilon.\hspace{6ex}\square$$
\indent The corollary below is straight forward by taking $Y=X$.
\begin{cor}\label{C5a1} Let $(X,d)$ be a metric space. Let $f:X\to X$ and $g:X \to X$ be two functions on $X$. If $f$ and $g$ are mutually contractive then there exists a Cauchy sequence $\{x_i\}_{i\in \mathbb{N}}$ in $X$ such that $\{f(x_i)\}_{n\in \mathbb{N}}$ and $\{g(x_i)\}_{n\in \mathbb{N}}$ form a Cauchy pair.
\end{cor}
\section{Partial Metric Space}
\label{C5b}
\indent \indent The results in this section are an extension of the work of Matthews et al. \cite{Mat1994,Mat2009} and  Karapinar et al. \cite{Kar20111,Kar20131}.  Although their theorems were quite  elegant, we felt that by requiring their Cauchy sequences to have a central distance $r=0$,  the partial metric spaces were not used to their full potential. In \cite{Ass20151}, we give  contracting criteria on a function $f$ so that the central distance of the Cauchy sequence generated is not restricted to be $0$.
\begin{defn}\label{D5b2} Let  $(X,p)$ be a partial metric space with $x_o$ in $X$. Let $f:X\to X$ be a function on $X$. Let $r$ and $0<c<1$ be two real numbers. We say that $f$ is an \uline{\textbf{orbital $c_r-$contraction at $x_o$}} (or $f$ is \uline{\textbf{orbitally $c_r-$contractive at $x_o$}}) if and only if  for all natural numbers $i$,
$$r \le p(f^{i}(x_o),f^{i}(x_0))$$
and
$$p(f^{i+2}(x_o),f^{i+1}(x_0) )\le r +c^{i+1}\vert p(f(x_o),x_o)\vert.$$
\\We say that $f$ is an \textbf{\uline{orbital $c_r-$contraction}} (or $f$ is \textbf{\uline{orbitally $c_r-$contractive}}) if and only if for every $x$ in $X$, $f$ is a orbital $r-$contraction at $x$.
\end{defn}
\begin{lem}\label{L5b1} Let  $(X,p)$ be a partial metric space with $x_o$ in $X$. Let $f:X\to X$ be a function on $X$. If $f$ is an orbital $c_r-$contraction at $x_o$ then $f$ is a Cauchy function at $x_o$.
\end{lem}
\textbf{Proof:} Let $x_o\in X$ and suppose $f: X \to X$ is an orbital $c_r-$contraction at $x_o$. Denote  $x_i =f^i(x_o)$. We now move to prove that $\{x_i\}_{i\in\mathbb{N}}$ is a Cauchy sequence. Let $M= \vert p(x_1,x_o)\vert=\vert p(f(x_o),x_o)\vert $. Then from \Cref{D5b2} and by (p-lbnd)
we get for all $i$,
$$r\le p(x_i,x_i)\le p(x_{i+1},x_i) \le r+c^iM \hspace{6ex} (\triangle)$$
and, hence,
$$-p(x_i,x_i)\le -r.$$
For all $j \ge i$, and by  using (p-inq)
$$p(x_j,x_i)\le p(x_j,x_{i+1})+p(x_{i+1},x_i)-p(x_{i+1},x_{i+1})$$
by $(\triangle)$
$$\le p(x_j,x_{i+1})+r+c^iM-r=p(x_j,x_{i+1})+c^iM$$
by (p-inq)
$$\le p(x_j,x_{i+2})+p(x_{i+2},x_{i+1})-p(x_{i+1},x_{i+1})+c^iM$$
by $(\triangle)$
$$\le p(x_j,x_{i+2})+r+c^{i+1}M-r+c^iM$$
$$=p(x_j,x_{i+2})+c^iM+c^{i+1}M$$
by repeating this process
$$\le p(x_j,x_{j-1})+ \sum\limits_{t=i}^{j-2} c^tM $$
by $(\triangle)$
$$\le r+c^{j-1}M+ \sum\limits_{t=i}^{j-2} c^tM$$
$$=r+ \sum\limits_{t=i}^{j-1} c^tM=r+c^i M\sum\limits_{t=0}^{j-i-1} c^t.$$
We know that $M \ge0$ and, hence, from the geometric series formula
$$p(x_j,x_i)\le r+c^i M\sum\limits_{t=0}^{j-i-1} c^t$$
$$\le r+c^i M\sum\limits_{t=0}^{+\infty} c^t=r+c^i\frac{M}{1-c}.$$
Since $0<c<1$, then for every positive real number $\epsilon$ there exists a natural number $N$ such that
$$c^N\frac{M}{1-c}<\epsilon$$
and hence, by (p-lbnd) and $(\triangle)$ for all $j\ge i >N$,
$$r-\epsilon <r \le p(x_i,x_i)\le p(x_j,x_i)$$
$$\le r+c^i\frac{M}{1-c}<r+c^N\frac{M}{1-c}<r+\epsilon. \hspace{6ex}\square$$
\indent As we mentioned before, Karapinar et al. \cite{Kar20111} generalized what was defined as a $\varphi-$weak contraction on a metric space \cite{Rho2001} to what they called a weak $\varphi-$contraction on a partial metric space. Their proposed generalization contracted over the whole space and the central distance of the obtained Cauchy sequence was forced to be zero.
In \cite{Ass20151}, we relax the constraints on the function $f$ requiring it  only to be contracting  on an orbit. We also allow the central distance of the obtained Cauchy sequence to be any arbitrary real number $r$.
\begin{defn}\label{D5b3} Let  $(X,p)$ be a partial metric space with $x_o$ in $X$. Let $f:X\to X$ be a function on $X$. Let $r$ be a real number and  $\varphi:[r,+\infty )\subset\mathbb{R}\to [0, + \infty)$  be a non-decreasing function such that
$$\varphi(t)=0 \text{ if and only if }t=r.$$
We say that $f$ is an \uline{\textbf{orbital $\varphi_{r}-$Contraction at $x_o$}} (or $f$ is \textbf{\uline{orbitally $\varphi_r-$contractive at $x_o$}}) if and only if for all $i$ and $j$,
$$r\le p(f^i(x_o),f^i(x_o))$$ 
and
$$p(f^{i+1}(x_o),f^{j+1}(x_o))\le p(f^i(x_o),f^j(x_o))-\varphi (p(f^i(x_o),f^j(x_o)).$$
We say that $f$ is an \uline{\textbf{orbital $\varphi_r-$contraction}} (or $f$ is \textbf{\uline{orbitally $\varphi_r-$contractive}} ) if and only if for every $x$ in $X$, $f$ is an orbital $\varphi_r$-contraction at $x$.
\end{defn}
\begin{lem}\label{L5b2} Let  $(X,p)$ be a partial metric space with $x_o$ in $X$. Let $f:X\to X$ be a function on $X$. If $f$ is an orbital $\varphi_r$-contraction at $x_o$ then $f$ is a Cauchy function at $x_o$.
\end{lem}
\textbf{Proof:} Let $x_o \in X$ and suppose $f: X \to X$ is an orbital $\varphi_r-$contraction at $x_o$. Denote  $x_i =f^i(x_o)$.
\\\indent \uline{Step 1:} Let $t_i=p(x_{i+1},x_i)$. In this step we will show that in the topological space $\mathbb{R}$ endowed with the standard topology, $\{t_i\}_{i\in\mathbb{N}}$ is a Cauchy sequence that converges to $r$.
\\From (p-lbnd) and since $\varphi (t_i)\ge 0$, 
$$r\le p(x_{i+1},x_{i+1})\le p(x_{i+2},x_{i+1})=t_{i+1}$$
and
$$t_{i+1}=p(x_{i+2},x_{i+1})\le p(x_{i+1},x_i)-\varphi (p(x_{i+1},x_i))$$
$$=t_i-\varphi(t_i)\le t_i.$$
Hence, for all $i$,
$$r\le t_{i+1}\le t_i$$
i.e. $\{t_i\}_{i\in \mathbb{N}}$ is a non-increasing sequence in $\mathbb{R}$ bounded below by $r$ and, therefore, $\{t_i\}_{i\in \mathbb{N}}$ is a Cauchy sequence in $\mathbb{R}$. Since $\mathbb{R}$ with the standard topology is a complete metric space, $\{t_i\}_{i\in\mathbb{N}}$ has a limit $L$ such that for all $i$,
$$t_i\ge L \ge r$$
and, since $\varphi$ is a non-decreasing function,
$$\varphi(t_i) \ge \varphi(L) \ge \varphi(r)=0$$
i.e.
$$-\varphi(t_i) \le -\varphi(L)\le 0.$$
Hence, by \Cref{D5b3}
$$r\le t_{i+1}\le t_i-\varphi(t_i)\le t_i-\varphi(L)$$
$$\le t_{i-1}-\varphi(t_{i-1})-\varphi(L)\le t_{i-1}-2\varphi(L)$$
by induction
$$t_{i+1}\le t_1-i\varphi(L).$$
Assume that $L>r$ then by \Cref{D5b3} $\varphi(L)>0$. By taking $i>\frac{t_1-r}{\varphi(L)}$ we get 
$$t_{i+1} \le t_1-i\varphi(L)<t_1-\frac{t_1-r}{\varphi (L)}\varphi(L) =r$$ a contradiction since $t_i \ge r$. Therefore, $L=r$.
\\\indent \uline{Step 2:} We now show that $\{x_i\}_{i\in \mathbb{N}}$ is a Cauchy sequence with central distance $r$ by supposing that it is not (a contrapositive approach). From \Cref{D5b3} we know that
for all $i$ and $j$,
$$r\le p(x_i,x_i)\le p(x_i,x_j)$$
in particular, for $i=j$,
$$-p(x_i,x_i)\le -r.$$
Hence, if $\{x_i\}_{i\in \mathbb{N}}$ is not a Cauchy sequence with central distance $r$ then there exists a positive real number $\delta$ such that for every natural number $N$, there exists  $i,j>N$ where 
$$p(x_i,x_j)\ge r+\delta >r$$
and from step 1, by choosing $N$  big enough 
$$r\le p(x_i,x_i)\le p(x_{i},x_{i+1})<r+\delta.$$
Then there exist  $j_1>m_1>N$ such that  
$$p(x_{m_1},x_{j_1})\ge r+\delta>r.$$
Let $n_1$ be the smallest number with $n_1 > m_1 $ and 
$$p( x_{m_1},x_{n_1}) \ge r+\delta.$$
Note 
$$p(x_{m_1},x_{n_1-1})<r+\delta.$$
There exist $j_2>m_2>n_1$ such that $$p(x_{m_2},x_{j_2})\ge r+\delta>r.$$
Let $n_2$ be the smallest number with $n_2 > m_2 $ and  
$$p( x_{m_2},x_{n_2}) \ge r+\delta.$$
Then
$$p(x_{m_2},x_{n_2-1})<r+\delta.$$
\indent Continuing this process, we build two increasing sequences in $\mathbb{N}$, $\{m_k\}_{k\in\mathbb{N}}$ and $\{n_k\}_{k\in\mathbb{N}}$ such that for all  $k$,
$$p(x_{m_k},x_{n_k-1})<r+\delta\le p(x_{m_k},x_{n_k}).$$
For all  $k$, denote $s_k=p(x_{m_k},x_{n_k})$. By (p-inq) 
$$s_k=p(x_{m_k},x_{n_k})\le p(x_{m_k},x_{n_k-1})+p(x_{n_k-1},x_{n_k})-p(x_{n_k-1},x_{n_k-1})$$
by (p-sym) and Step 1
$$=p(x_{m_k},x_{n_k-1})+p(x_{n_k},x_{n_k-1})-p(x_{n_k-1},x_{n_k-1})$$
$$\le p(x_{m_k},x_{n_k-1})+t_{n_k-1}-r$$
and, hence,
$$ s_k \le p(x_{m_k},x_{n_k-1})+t_{n_k-1}-r.\hspace{6ex} (\triangle')$$
Additionally, since $\{t_i\}_{i\in\mathbb{N}}$ is a Cauchy sequence tending to $r$, for every positive real number $\epsilon$ there exists a natural number $N$ such that for all $n_k-1>N$,
$$r \le t_{n_k-1}<r+\epsilon.$$
Since $\{m_k\}_{k\in\mathbb{N}}$ and $\{n_k\}_{k\in\mathbb{N}}$ are increasing sequences, there exists a natural number $N'$ such that for all $k>N'$, $n_k-1>N$. Therefore for all $k>N'$ and since $ r+\delta \le s_k$,
$$0\le s_k-(r+\delta)$$
from $(\triangle')$
$$ \le p(x_{m_k},x_{n_k-1})+t_{n_k-1}-r-(r+\delta)$$
$$<(r+\delta )+(r+\epsilon) -r-(r+\delta)=\epsilon.$$
Therefore,  $\{s_k\}_{k\in\mathbb{N}}$ is a Cauchy sequence with $r+\delta$ as a limit. On the other hand by applying (p-inq) twice we get
$$s_k=p(x_{m_k},x_{n_k})\le p(x_{m_k},x_{n_k+1})+p(x_{n_k+1},x_{n_k})-p(x_{n_k+1},x_{n_k+1})$$
$$\le p(x_{m_k},x_{m_k+1})+p(x_{m_k+1},x_{n_k+1})-p(x_{m_k+1},x_{m_k+1})+p(x_{n_k+1},x_{n_k})-p(x_{n_k+1},x_{n_k+1})$$
from (p-sym) and Step 1
$$=p(x_{m_k+1},x_{m_k})+p(x_{m_k+1},x_{n_k+1})-p(x_{m_k+1},x_{m_k+1})+p(x_{n_k+1},x_{n_k})-p(x_{n_k+1},x_{n_k+1})$$
$$=t_{m_k}+p(x_{m_k+1},x_{n_k+1})-p(x_{m_k+1},x_{m_k+1})+t_{n_k}-p(x_{n_k+1},x_{n_k+1})$$
$$\le t_{m_k}+p(x_{m_k+1},x_{n_k+1})-r+t_{n_k}-r$$
from \Cref{D5b3}
$$\le t_{m_k}+p(x_{m_k},x_{n_k})-\varphi (p(x_{m_k},x_{n_k}))-r +t_{n_k}-r$$
and, hence,
$$ s_k \le t_{m_k}+s_k-\varphi (s_k)+t_{n_k} -2r$$
i.e.
$$\varphi(s_k)\le t_{m_k}+t_{n_k}-2r.$$
Since $r<r+\delta \le s_k$ and from \Cref{D5b3} we get
$$0 <\varphi(r+\delta)\le \varphi(s_k)\le t_{m_k}+t_{n_k}-2r.$$
Since $\{t_i\}_{i\in\mathbb{N}}$ is a Cauchy sequence that tends to $r$ then for every positive real number $\epsilon$, there exists a natural number $N$ such that for all $m_k,n_k>N$, 
$$0 < \varphi(r+\delta)\le  t_{m_k}+t_{n_k}-2r<(r+\frac{\epsilon}{2})+(r+\frac{\epsilon}{2})-2r=\epsilon$$
and, hence,
$$0<\varphi(r+\delta)\le 0$$ a clear contradiction. Therefore, the assumption considered at the beginning of Step 2 is incorrect proving that $\{x_i\}_{i\in \mathbb{N}}$ is a Cauchy sequence with central distance $r.\hspace{6ex}\square$
\\\indent We now move to generalizing pairwise contractive functions.
\begin{defn}\label{D5b4} Let $(X,p)$ be a partial metric space with $x_o$ and $y_o$ in $X$. Let $f:X\to X$ and $g:X \to X$ be two functions on $X$. Let $r$ and $0<c<1$ be two real numbers. We say that $f$ and $g$ are \uline{\textbf{$f-$pairwise $c_r-$contractive over $(x_o,y_o)$}} if and only if  for all natural numbers $i$,
$$\begin{cases}r\le \min\{p(f^i(x_o),f^i(x_o)),p(g^i(y_o),g^i(y_o))\}
\\ p(f^{i+1}(x_0),g^{i}(y_o))\le r+ c^iM
\\p(f^{i}(x_0),g^{i}(y_o))\le r+c^iM
\end{cases}$$
\\where $M=\max\{\vert p(f(x_o),y_o)\vert,\vert p(x_o,y_o)\vert\}.$
\end{defn}
\begin{thm}{\textbf{(Cauchy $f-$Pairwise $c_r-$Contractive):}}
\label{T5b1}
\\Let $(X,p)$ be a partial metric space with $x_o$ and $y_o$ in $X$. Let $f:X\to X$ and $g:X \to X$ be two functions on $X$. If  $f$ and $g$ are $f-$pairwise $c_r-$contractive over $(x_o,y_o)$ then $f$ and $g$ form a Cauchy pair over $(x_o,y_o)$.
\end{thm}
\textbf{Proof:} For all  $i$, let $x_i=f^i(x_o)$  and $y_i=g^i(y_o)$. Let $M=\max\{\vert p(f(x_o),y_o)\vert,\vert p(x_o,y_o)\vert\}$.
\\\indent First we show that for all $i\ge j$,
$$r\le \min\{ p(x_i,x_i), p(y_j,y_j)\}\le p(x_i,y_j)\le r+ c^j\frac{2M}{1-c}.$$
 From \Cref{D5b4} we know that
$$\begin{cases}r\le \min\{p(x_i,x_i),p(y_i,y_i)\}
\\ p(x_{i+1},y_i)\le r+ c^iM
\\p(x_i,y_i)\le r+c^iM
\end{cases}\hspace{6ex} (\tilde\triangle)$$
In particular, for all $i$,
$$-p(x_i,x_i)\le -r \text{ and }-p(y_i,y_i)\le-r.$$
\indent By \Cref{T4b1} it suffices to bound $p(x_i,y_j)$ for $i\ge j$. Let us first investigate what happens for the specific values of $i=6$ and $j=3$. By repeatedly using (p-inq) we get
$$p(x_6,y_3)\le p(x_6,y_5)-p(y_5,y_5)+p(y_5,x_5)-p(x_5,x_5)+p(x_5,y_4)-p(y_4,y_4)+p(y_4,x_4)-p(x_4,x_4)+p(x_4,y_3)$$
by (p-sym)
$$=p(x_6,y_5)-p(y_5,y_5)+p(x_5,y_5)-p(x_5,x_5)+p(x_5,y_4)-p(y_4,y_4)+p(x_4,y_4)-p(x_4,x_4)+p(x_4,y_3)$$
by $(\tilde\triangle)$
$$\le (r+c^5M)-r+(r+c^5M)-r+(r+c^4M)-r+(r+c^4M)-r+r+c^3M$$
since $c^jM\ge 0$
$$\le r+2[c^5M+c^4M+c^3M].$$
We now move to derive an upper bound for   $p(x_i,y_j)$ where $i\ge j$. 
\\\textit{Case 1:} If $i=j$, by $(\tilde\triangle)$ and since $0<c<1$ and $c^jM\ge 0$,
$$r\le \min\{ p(x_j,x_j), p(y_j,y_j)\}\le p(x_j,y_j)\le r+c^jM\le r+c^j\frac{2M}{1-c}.$$
\textit{Case 2:} If $i>j$,  by $(\tilde\triangle)$ and (p-lbnd) we get
$$r \le p(x_i,x_i)\le p(x_i,y_j)$$
by repeatedly using (p-inq)
$$\le \sum\limits_{t=j+1}^{i-1}[p(x_{t+1},y_t)-p(y_t,y_t)+p(y_t,x_t)-p(x_t,x_t)]+p(x_{j+1},y_j)$$
by (p-sym)
$$\le \sum\limits_{t=j+1}^{i-1}[p(x_{t+1},y_t)-p(y_t,y_t)+p(x_t,y_t)-p(x_t,x_t)]+p(x_{j+1},y_j)$$
by $(\tilde\triangle)$
$$\le \sum\limits_{t=j+1}^{i-1}[(r+c^tM-r+(r+c^tM)-r]+r+c^jM=2\sum\limits_{t=j+1}^{i-1}[c^tM]+r+c^jM$$
since $c^tM\ge 0$
$$\le r+2\sum\limits_{t=j}^{i-1}c^tM=r+2c^j\sum\limits_{t=0}^{i-j-1}c^tM\le r+2c^j\sum\limits_{t=0}^{+\infty}c^tM$$
finally, using the geometric series formula
$$=r+c^j\frac{2M}{1-c}.$$
Hence, for all $i\ge j$,
$$r\le \min\{ p(x_i,x_i), p(y_j,y_j)\}\le p(x_i,y_j)\le r+ c^j\frac{2M}{1-c}.$$
\indent  For all positive real numbers $\epsilon$ there exists a natural number $N$ such that $c^N\frac{2M}{1-c}<\epsilon.$ Hence, for all $i\ge j>N$,
$$r-\epsilon <r\le \min\{ p(x_i,x_i), p(y_j,y_j)\}\le p(x_i,y_j)\le r+ c^j\frac{2M}{1-c}<r+\frac{M}{1-c}c^N<r+\epsilon.$$
Therefore, by \Cref{T4b1}, $f$ and $g$ are Cauchy pairs over $(x_o,y_o)$.$\hspace{6ex}\square$
\begin{defn}\label{D5b6} Let  $(X,p)$ and $(Y,h)$ be two partial metric spaces and suppose $f:X\to Y$ and $g:X \to Y$ are two functions on $X$. Let $r$, $ A\ge 0$ and $0<c<1$ be three real numbers. We say that  $f$ and $g$ are \textbf{\uline{$f-$mutually $c_r-$contractive}} if and only  for each $x$ in $X$ we can find an element $z$ of $X$ such that 
$$h(f(z),g(z)) -h(f(z),f(z))\le c[ h(f(x),g(x))-h(f(x),f(x))]$$
and
$$r \le p(z,z) \le p(x,z)\le r+A[ h(f(x),g(x))-h(f(x),f(x))].$$
\end{defn}
\indent \indent In the above definition, putting a heavier emphasis on one function  makes it much easier to apply the theorem below on a pair of functions when one is much more complex than the other. The above contraction is enough to generate a Cauchy sequence in $X$ as shown in \Cref{T5b2}. In the case of $(Y,h)$ being a strong partial metric space, $f-$mutually $c_r-$contraction is used to obtain a coincidence point. However, when $(Y,h)$ is a partial metric space, we need a stronger version.
\begin{defn}\label{D5b6-} Let  $(X,p)$ and $(Y,h)$ be two partial metric spaces and suppose $f:X\to Y$ and $g:X \to Y$ are two functions on $X$. Let $r$, $ A\ge 0$ and $0<c<1$ be three real numbers. We say that  $f$ and $g$ are \textbf{\uline{$(f,g)-$mutually $c_r-$contractive}} if and only  for each $x$ in $X$ we can find an element $z$ of $X$ such that 
$$h(f(z),g(z)) -h(f(z),f(z))\le c[ h(f(x),g(x))-h(f(x),f(x))]$$
$$h(f(z),g(z)) -h(g(z),g(z))\le c[ h(f(x),g(x))-h(g(x),g(x))]$$
$$r \le p(z,z) \le p(x,z)\le r+A[ h(f(x),g(x))-h(f(x),f(x))]$$
and
$$r \le p(z,z) \le p(x,z)\le r+A[ h(f(x),g(x))-h(g(x),g(x))].$$
\end{defn}
It is clear that if two functions $f$ and $g$ are $(f,g)-$mutually $c_r-$contractive then they are $f-$mutually $c_r-$contractive and $g-$mutually $c_r-$contractive.
 \begin{thm}{\textbf{(Cauchy $f-$Mutually $c_r$-Contractive):}}
\label{T5b2}
\\Let $(X,p)$ and $(Y,h)$ be two partial metric spaces. Let $f:X\to Y$ and $g:X \to Y$ be two functions on $X$. If $f$ and $g$ are $f-$mutually $c_r-$contractive then there exists a Cauchy sequence $\{x_i\}_{i\in \mathbb{N}}$ in $X$ with central distance  $r$ such that for all natural numbers $i$, $r\le p(x_i,x_i)$. Additionally for every positive real number $\epsilon$ there exists a natural number $N$ such that for all $i>N$,
$$h(f(x_i),g(x_i))-h(f(x_i),f(x_i))<\epsilon.$$
\end{thm}
\textbf{Proof:} Since $f$ and $g$ are $f-$mutually $c_r-$contractive then there exist two real numbers $0<c<1$ and $A\ge0$ such that for each $x$ in $X$, there exists a $z$ in $X$ such that  
$$h(f(z),g(z)) -h(f(z),f(z))\le c[ h(f(x),g(x))-h(f(x),f(x))]$$
and
$$r \le p(z,z) \le p(x,z)\le r+A[ h(f(x),g(x))-h(f(x),f(x))].$$
Let us take an arbitrary element $x_1$ in $X$ then there exists an element $x_2$ in $X$ such that
$$h(f(x_2),g(x_2))-h(f(x_2),f(x_2))\le c[h(f(x_1),g(x_1))-h(f(x_1),f(x_1))]$$
and $$r\le p(x_2,x_2)\le p(x_1,x_2)\le r+A[h(f(x_1),g(x_1))-h(f(x_1),f(x_1))].$$ There exists an element $x_3$ in $X$ such that 
$$h(f(x_3),g(x_3))-h(f(x_3),f(x_3))\le c[h(f(x_2),g(x_2))-h(f(x_2),f(x_2))]$$
and $$r\le p(x_3,x_3)\le p(x_2,x_3)\le r+A[h(f(x_2),g(x_2))-h(f(x_2),f(x_2))].$$ We continue the above process to generate the sequence $\{x_i\}_{i\in \mathbb{N}}$ such that for all  $i$,
$$h(f(x_{i+1}),g(x_{i+1}))-h(f(x_{i+1}),f(x_{i+1}))\le c[h(f(x_i),g(x_i))-h(f(x_i),f(x_i))]$$
and $$r\le p(x_{i+1},x_{i+1})\le p(x_i,x_{i+1})\le r+A[h(f(x_i),g(x_i))-h(f(x_i),f(x_i))]\hspace{6ex}(\ddot\triangle)$$
in particular, 
$$r\le p(x_i,x_i)$$
i.e.
$$- p(x_i,x_i)\le r.$$
\indent We now prove that $\{x_i\}_{i\in \mathbb{N}}$ is a Cauchy sequence with central distance $r$.
\\Step 1: For all  $i$, let $t_i=  h(f(x_i),g(x_i))-h(f(x_i),f(x_i))$. Then from $(\ddot\triangle)$
$$t_{i+1}=h(f(x_{i+1}),g(x_{i+1}))-h(f(x_{i+1}),f(x_{i+1}))$$
$$\le c[h(f(x_i),g(x_i))-h(f(x_i),f(x_i))]\le c^2[h(f(x_{i-1}),g(x_{i-1}))-h(f(x_{i-1}),f(x_{i-1}))]$$
and, hence, by induction
$$t_{i+1}\le  c^i[h(f(x_1),g(x_1))-h(f(x_1),g(x_1))]=c^i t_1. \hspace{6ex}(\check\triangle)$$
Therefore, and since $0<c<1$ for every positive real number $\epsilon$ there exists  natural number $N'$ such that for all $i>N'$,
$$h(f(x_i),g(x_i))-h(f(x_i),f(x_i))<\epsilon.$$
Step 2: For all  $j>i$ , by $(\ddot\triangle)$ and  repeatedly using (p-inq) we get,
$$r\le p(x_i,x_i)\le p(x_i,x_j)$$
$$\le [p(x_i,x_{i+1})-p(x_{i+1},x_{i+1})]+[p(x_{i+1},x_{i+2})-p(x_{i+2},x_{i+2})]+...+[p(x_{j-2},x_{j-1})-p(x_{j-1},x_{j-1})]+p(x_{j-1},x_{j})$$
by $(\ddot\triangle)$
$$\le [r+At_i-r]+[r+At_{i+1}-r]+...+[r+At_{j-2}-r]+r+At_{j-1}=r+\sum\limits_{k=i}^{j-1}At_k$$
by $(\check\triangle)$
$$\le r+\sum\limits_{k=i}^{j-1}Ac^{k-1}t_1=r+At_{1}c^{i-1}\sum\limits_{k=0}^{j-i-1}c^k$$
knowing that $At_1 \ge 0$ and by the geometric series formula
$$\le r+ At_1\sum\limits_{k=0}^{+\infty}c^k$$
by the geometric series formula
$$=r+At_1c^{i-1}\sum\limits_{k=0}^{+\infty}c^k=r+\frac{At_1}{c(1-c)}c^i.$$
Therefore, for every positive number $\epsilon$ there exists a natural number $N$ where  $\frac{At_1}{c(1-c)}c^N<\epsilon$ and, hence, for all  $j\ge i>N$,
$$r-\epsilon<r\le p(x_j,x_i)\le r+\frac{At_1}{c(1-c)}c^i<r+ \frac{At_1}{c(1-c)}c^N<r+\epsilon.\hspace{6ex}\square$$
\indent The corollary below is straight forward by taking $Y=X$.
\begin{cor}\label{C5b1} Let $(X,p)$  be a partial metric space. Let $f:X\to X$ and $g:X \to X$ be two functions on $X$. If $f$ and $g$ are $f-$mutually $c_r-$contractive then there exists a Cauchy sequence $\{x_i\}_{i\in \mathbb{N}}$ in $X$ with central distance $r$ such that for all natural numbers $i$, $r\le p(x_i,x_i)$. Additionally for every positive real number $\epsilon$ there exists a natural number $N$ such that for all $i>N$,
$$p(f(x_i),g(x_i))-p(f(x_i),f(x_i))<\epsilon.$$
\end{cor}
\section{Partial $n-\mathfrak{M}$etric Space}
\label{C5c}
\indent \indent The inspiration for this section mainly came from  Ayadi et al. \cite{Kar20121}. As explained in \Cref{C2}, we have relaxed the axioms of a $G_p-$metric (see \cite{Zan2011}) to obtain the partial $3-\mathfrak{M}$etric. We also have a much less restrictive condition on the contracting functions and the central distance of a Cauchy sequence. To this end, we will use \Cref{T4c1} to check whether  a sequence is Cauchy or not by  comparing the elements pairwise rather than comparing $n-$tuples.
\\\indent Computations with partial $n-\mathfrak{M}$etrics when $n>2$ require more attention than their partial metric counterparts because of the following: Let $P$ be a partial $n-\mathfrak{M}$etric on $X$. If $n=2$ then $P$ is a partial metric. Hence, for any two elements $a$ and $b$ in a set $X$ we have by (p-sym) 
$$ P(\langle a\rangle^{2-1},b)=P(\langle b\rangle^{2-1},a).$$
In the more general case with $n>2$, we have by \Cref{C2e3},
$$P(\langle a\rangle^{n-1},b) \le (n-1)P(\langle b\rangle^{n-1},a)-(n-2)P(\langle b \rangle^n).$$
\indent The proofs of \Cref{C5b} generalize to the proofs of \Cref{C5c} by adding steps and considering a different $\epsilon$ to compensate for the that fact. Otherwise, the proofs are very similar to the proofs in \Cref{C5b} and, hence, may be skipped if the reader so desires. 

\begin{defn}\label{D5c2}
Let  $(X,P)$ be a partial $n-\mathfrak{M}$etric space with $x_o$ in $X$ and suppose $f:X\to X$ is a function on $X$. Let $r$ and  $0<c<1$  be two real numbers. We say that $f$ is an \uline{\textbf{orbital $c_r-$contraction at $x_o$}} (or $f$ is \uline{\textbf{orbitally $c_r-$contractive at $x_o$}}) if and only if  for all natural numbers $i$,
$$r \le P(\langle f^{i+1}(x_o)\rangle^n)$$
and
$$ P(\langle f^{i}(x_o)\rangle^{n-1},f^{i+1}(x_0) )\le r +c^{i}\vert P(\langle x_o\rangle ^{n-1},f(x_o))\vert.$$
We say that $f$ is an \uline{\textbf{orbital $c_r-$contraction }} (or $f$ is \uline{\textbf{orbitally $c_r-$contractive}}) if and only if for every $x$ in $X$, $f$ is an orbital $c_r-$contraction at $x$.
\end{defn}
\begin{lem}\label{L5c1} Let  $(X,P)$ be a partial $n-\mathfrak{M}$etric space with $x_o$ in $X$. Let $f:X\to X$ be a function on $X$. If $f$ is an orbital $c_r-$contraction at $x_o$ then $f$ is a Cauchy function at $x_o$.
\end{lem}
\textbf{Proof:} Let $x_o\in X$ and suppose $f: X \to X$ is an orbital $c_r-$contraction at $x_o$. Denote  $x_i =f^i(x_o)$. We now move to prove that $\{x_i\}_{i\in\mathbb{N}}$ is a Cauchy sequence. Let $$M= \vert P(\langle x_o\rangle ^{n-1},x_1)\vert=\vert P(\langle x_o\rangle ^{n-1},f(x_o))\vert .$$ Then, from \Cref{D5c2} and by ($P_n$-lbnd)
we get for all $i$,
$$r \le P(\langle x_i\rangle^n) \le P(\langle x_i\rangle^{n-1},x_{i+1}) \le r +c^{i}M \hspace{6ex}(\otimes)$$
and, hence,
$$-P(\langle x_i\rangle^n)\le -r.$$\
For all $j \ge i$, and by using ($P_n$-inq)
$$P(\langle x_i\rangle^{n-1},x_{j})\le P(\langle x_i\rangle^{n-1},x_{i+1})+P(\langle x_{i+1}\rangle^{n-1},x_{j})-P(\langle x_{i+1}\rangle^{n})$$
by $(\otimes)$
$$\le r+c^iM +P(\langle x_{i+1}\rangle^{n-1},x_{j})-r=P(\langle x_{i+1}\rangle^{n-1},x_{j})+c^iM$$
by ($P_n$-inq)
$$\le P(\langle x_{i+1}\rangle^{n-1},x_{i+2})+P(\langle x_{i+2}\rangle^{n-1},x_{j})-P(\langle x_{i+2}\rangle^{n})+c^iM$$
by $(\otimes)$
$$\le r+c^{i+1}M+P(\langle x_{i+2}\rangle^{n-1},x_{j})-r+c^iM$$
$$=P(\langle x_{i+2}\rangle^{n-1},x_{j})+c^iM+c^{i+1}M$$
by repeating this process
$$\le P(\langle x_{j-1}\rangle^{n-1},x_{j})+\sum\limits_{k=i}^{j-2}c^kM$$
by $(\otimes)$
$$\le r +c^{j-1}+\sum\limits_{k=i}^{j-2}c^kM=r+M\sum\limits_{k=i}^{j-1}c^k=r+Mc^i\sum\limits_{k=0}^{j-i-1}c^k.$$
We know that $M\ge 0$and, hence, from the geometric series formula
$$ P(\langle x_i\rangle^{n-1},x_{j})\le r+Mc^i\sum\limits_{k=0}^{j-i-1}c^k$$
$$\le r+Mc^i\sum\limits_{k=0}^{+\infty}c^k=r+c^i\frac{M}{1-c}.$$
Since $0<c<1$, then for every positive real number $\epsilon$ there exists a natural number $N$ such that
$$c^N\frac{M}{1-c}<\epsilon$$
and, hence, for all $j\ge i>N$,
$$r-\epsilon < r \le P(\langle x_i\rangle^{n-1},x_{j}) \le r+c^i\frac{M}{1-c}<r+c^N\frac{M}{1-c}<r+\epsilon.$$
Therefore, by \Cref{T4c1}, $\{x_i\}_{i\in\mathbb{N}}$ is a Cauchy sequence with central distance $r$.$\hspace{6ex} \square$
\\\indent Bilgili et al. \cite{Kar20132} generalized the idea of a weak contraction into a $G_p-$metric space. We build on his work and generalize it to a partial $n-\mathfrak{M}$etric case.
\begin{defn}\label{D5c3}
Let  $(X,P)$ be a partial $n-\mathfrak{M}$etric space with $x_o$ in $X$ and suppose $f:X\to X$ is a function on $X$. Let $r$ be a real number and  $\varphi:[r,+\infty )\subset\mathbb{R}\to [0, + \infty)$ be a non-decreasing function such that
$$\varphi(t)= 0 \text{ if and only if } t=r.$$
We say that $f$ is an \uline{\textbf{orbital $\varphi_{r}-$contraction at $x_o$}} (or $f$ is \textbf{\uline{orbitally $\varphi_r-$contractive at $x_o$}}) if and only if 
for all $i$ and $j$,
$$ r \le P(\langle f^{i}(x_o)\rangle^{n})$$
and
$$P(\langle f^{i+1}(x_o)\rangle^{n-1}, f^{j+1}(x_o))\le P(\langle f^{i}(x_o)\rangle^{n-1}, f^{j}(x_o))-\varphi(P(\langle f^{i}(x_o)\rangle^{n-1}, f^{j}(x_o))).$$
We say that $f$ is an \uline{\textbf{orbital $\varphi_r-$contraction}} (or $f$ is \textbf{\uline{orbitally $\varphi_r-$contractive}}) if and only if for every $x$ in $X$, $f$ is an orbital $\varphi_r$-contraction at $x$.\end{defn}
\begin{lem}\label{L5c2} Let  $(X,P)$ be a partial $n-\mathfrak{M}$etric space with $x_o$ in $X$. Let $f:X\to X$ be a function on $X$. If $f$ is an orbital $\varphi_r$-contraction at $x_o$ then $f$ is a Cauchy function at $x_o$.

\end{lem}
\textbf{Proof:} Let $x_o \in X$ and suppose $f: X \to X$ is an orbital $\varphi_r-$contraction at $x_o$. Denote  $x_i =f^i(x_o)$.
\\\indent \uline{Step 1:} Let $t_i=P(\langle x_{i}\rangle^{n-1}, x_{i+1})$. In this step we will show that in the topological space $\mathbb{R}$ endowed with the standard topology, $\{t_i\}_{i\in\mathbb{N}}$ is a Cauchy sequence that converges to $r$.
\\From ($P_n$-lbnd) and \Cref{D5c3},
$$r\le P(\langle x_{i+1}\rangle^{n})\le P(\langle x_{i+1}\rangle^{n-1}, x_{i+2})=t_{i+1}$$
$$\le P(\langle x_{i}\rangle^{n-1}, x_{i+1})-\varphi(P(\langle x_{i}\rangle^{n-1}, x_{i+1}))$$
since $\varphi(t_i)\ge 0$
$$=t_i-\varphi(t_i)\le t_i.$$
Hence, for all $n$,
$$r\le t_{i+1}\le t_i$$
i.e. $\{t_i\}_{i\in\mathbb{N}}$ is a non-increasing sequence in $\mathbb{R}$ bounded below by $r$ and, therefore, $\{t_i\}_{i\in\mathbb{N}}$ is a Cauchy sequence in $\mathbb{R}$. Since $\mathbb{R}$ with the standard topology is a complete metric space, $\{t_i\}_{i\in\mathbb{N}}$ has a limit $L$ such that for all $i$,
$$\varphi(t_i)\ge \varphi(L) \ge \varphi(r)=0$$
i.e.
$$-\varphi(t_i)\le-\varphi(L)\le 0.$$
Hence, by \Cref{D5c3}
$$r\le t_{i+1}\le t_i-\varphi(t_i)\le t_i-\varphi(L)$$
$$\le t_{i-1}-\varphi(t_{i-1})-\varphi(L)\le t_{i-1}-2\varphi(L)$$
by induction
$$t_{i+1}\le t_1-i\varphi(L).$$
Assume that $L>r$ then by \Cref{D5c3} $\varphi(L)>0$. By taking $i>\frac{t_1-r}{\varphi(L)}$ we get 
$$t_{i+1} \le t_1-i\varphi(L)<t_1-\frac{t_1-r}{\varphi (L)}\varphi(L) =r$$ a contradiction since $t_i \ge r$. Therefore, $L=r$.
\\\indent \uline{Step 2:} We now show that $\{x_i\}_{i\in \mathbb{N}}$ is a Cauchy sequence with central distance $r$ by supposing that it is not (a contrapositive approach). To do that we refer the reader back to \Cref{T4c1} (c).  From \Cref{D5c3} we know that
for all $i$ and $j$,
$$r\le P(\langle x_{i}\rangle^{n})\le P(\langle x_{i}\rangle^{n-1}, x_{i+1})$$
in particular, for $i=j$
$$-P(\langle x_{i}\rangle^{n})\le -r.$$
Hence, if $\{x_i\}_{i\in \mathbb{N}}$ is not a Cauchy sequence with central distance $r$ then by \Cref{T4c1} (d) there exists a positive real number $\delta$ such that for every natural number $N$, there exists $j\ge i>N$ where 
$$P(\langle x_{i}\rangle^{n-1}, x_{j})\ge r+\delta >r$$
and from step 1, by choosing $N$  big enough 
$$r\le P(\langle x_{i}\rangle^{n})\le P(\langle x_{i}\rangle^{n-1}, x_{i+1})<r+\delta.$$
Then there exist  $j_1>m_1>N$ such that  
$$P(\langle x_{m_1}\rangle^{n-1}, x_{j_1}))\ge r+\delta>r.$$
Let $n_1$ be the smallest number with $n_1 > m_1 $ and 
$$P(\langle x_{m_1}\rangle^{n-1}, x_{n_1}) \ge r+\delta.$$
Note
$$P(\langle x_{m_1}\rangle^{n-1}, x_{n_1-1})<r+\delta.$$
There exist $j_2>m_2>n_1$ such that $$P(\langle x_{m_2}\rangle^{n-1}, x_{j_2}) \ge r+\delta\ge r+\delta>r.$$
Let $n_2$ be the smallest number with $n_2 > m_2 $ and 
$$P(\langle x_{m_2}\rangle^{n-1}, x_{n_2}) \ge r+\delta \ge r+\delta.$$
Then 
$$P(\langle x_{m_2}\rangle^{n-1}, x_{n_2-1})<r+\delta.$$
\indent Continuing this process, we build two increasing sequences in $\mathbb{N}$, $\{m_k\}_{k\in\mathbb{N}}$ and $\{n_k\}_{k\in\mathbb{N}}$ such that for all  $k$,
$$P(\langle x_{m_k}\rangle^{n-1}, x_{n_k-1})<r+\delta\le P(\langle x_{m_k}\rangle^{n-1}, x_{n_k}).$$
For all  $k$, denote $s_k=P(\langle x_{m_k}\rangle^{n-1}, x_{n_k})$. By ($P_n$-inq) 
$$s_k=P(\langle x_{m_k}\rangle^{n-1}, x_{n_k})\le P(\langle x_{m_k}\rangle^{n-1}, x_{n_k-1})+P(\langle x_{n_k-1}\rangle^{n-1}, x_{n_k})-P(\langle x_{n_k-1}\rangle^{n})$$
by Step 1
$$s_k \le  P(\langle x_{m_k}\rangle^{n-1}, x_{n_k-1}) +t_{n_k-1}-r \hspace{6ex}(\tilde\otimes)$$
Additionally, since $\{t_i\}_{i\in\mathbb{N}}$ is a Cauchy sequence tending to $r$, for every positive real number $\epsilon$ there exists a natural number $N$ such that for all $n_k-1>N$,
$$r \le t_{n_k-1}<r+\epsilon.$$
Since $\{m_k\}_{k\in\mathbb{N}}$ and $\{n_k\}_{k\in\mathbb{N}}$ are increasing sequences, there exists a natural number $N'$ such that for all $k>N'$, $n_k-1>N$. Therefore for all $k>N'$ and since $ r+\delta \le s_k$,
$$0\le s_k-(r+\delta)$$
from $(\tilde\otimes)$
$$ \le P(\langle x_{m_k}\rangle^{n-1}, x_{n_k-1})+t_{n_k-1}-r-(r+\delta)$$
$$<(r+\delta )+(r+\epsilon) -r-(r+\delta)=\epsilon.$$
Therefore,  $\{s_k\}_{k\in\mathbb{N}}$ is a Cauchy sequence with $r+\delta$ as a limit. On the other hand by applying ($P_n$-inq)  we get
$$s_k=P(\langle x_{m_k}\rangle^{n-1}, x_{n_k})\le P(\langle x_{m_k}\rangle^{n-1}, x_{m_k+1})+P(\langle x_{m_k+1}\rangle^{n-1}, x_{n_k})-P(\langle x_{m_k+1}\rangle^{n})$$
by Step 1
$$\le t_{m_k}+P(\langle x_{m_k+1}\rangle^{n-1}, x_{n_k})-r$$
by ($P_n$-inq)
$$\le t_{m_k}+P(\langle x_{m_k+1}\rangle^{n-1}, x_{n_k+1})+P(\langle x_{n_k+1}\rangle^{n-1}, x_{n_k})-P(\langle x_{n_k+1}\rangle^{n})-r$$
by \Cref{C2e3}
$$\le t_{m_k}+ P(\langle x_{m_k+1}\rangle^{n-1}, x_{n_k+1})+(n-1)P(\langle x_{n_k}\rangle^{n-1}, x_{n_k+1})-(n-2)P(\langle x_{n_k}\rangle^{n})-P(\langle x_{n_k+1}\rangle^{n})-r$$
by Step 1
$$\le t_{m_k} +P(\langle x_{m_k+1}\rangle^{n-1}, x_{n_k+1})+(n-1)t_{n_k}-(n-2)r-r-r$$
by \Cref{D5c3}
$$\le t_{m_k} +P(\langle x_{m_k}\rangle^{n-1}, x_{n_k})-\varphi(P(\langle x_{m_k}\rangle^{n-1}, x_{n_k}))+(n-1)t_{n_k}-(n-2)r-r-r $$
and, hence,
$$s_k \le t_{m_k}+s_k-\varphi(s_k)+(n-1)t_{n_k}-nr$$
therefore, 
$$\varphi(s_k) \le t_{m_k}+(n-1)t_{n_k}-nr.$$
Since $r<r+\delta\le s_k$ and from \Cref{D5c3} we get
$$0<\varphi(r+\delta)\le\varphi(s_k)\le t_{m_k}+(n-1)t_{n_k}-nr$$
Additionally, since $\{t_i\}_{i\in\mathbb{N}}$ is a Cauchy sequence tending to $r$, for every positive real number $\epsilon$ there exists a natural number $N$ such that for all $n_k>m_k>N$,
$$r\le t_{n_k}\le t_{m_k} <r+\frac{\epsilon}{n}$$
therefore, for every positive real number $\epsilon$ there exists a natural number $N$ such that for all $n_k>m_k>N$,

$$ 0< \varphi(r+\delta) < r+\frac{\epsilon}{n}+(n-1)(r+\frac{\epsilon}{n})-nr=\epsilon$$
and, hence,
$$0<\varphi (r+\delta)\le 0$$
a clear contradiction. Therefore, the assumption considered at the beginning of Step 2 is incorrect proving that $\{x_i\}_{i\in \mathbb{N}}$ is a Cauchy sequence with central distance $r.\hspace{6ex}\square$
\begin{defn}\label{D5c4}Let $(X,P)$ be a partial $n-\mathfrak{M}$etric space with $x_o$ and $y_o$ in $X$ and suppose $f:X\to X$ and $g:X \to X$ are two functions on $X$. Let $r$ and  $0<c<1$ be two real numbers. We say that $f$ and $g$ are \uline{\textbf{$f-$pairwise $c_r-$contractive over $(x_o,y_o)$}} if and only if for all  natural numbers $i$,
$$\begin{cases}r\le \min\{P(\langle f^i(x_o)\rangle^n),P(\langle g^i(x_o)\rangle^n\}
\\ P(\langle f^{i+1}(x_o)\rangle^{n-1},g^{i}(y_o))\le r+ c^iM
\\P(\langle f^{i}(x_o)\rangle^{n-1},g^{i}(y_o))\le r+ c^iM
\end{cases}$$
\\where $M=\max\{\vert P(\langle f(x_o)\rangle^{n-1},y_o)\vert,\vert P(\langle x_o\rangle^{n-1},y_o)\vert\}.$
\end{defn}

\begin{thm}{\textbf{(Cauchy $f-$Pairwise $c_r-$Contractive):}}
\label{T5c1}
\\Let $(X,P)$ be a partial $n-\mathfrak{M}$etric space with $x_o$ and $y_o$ in $X$. Let $f:X\to X$ and $g:X \to X$ be two functions on $X$. If  $f$ and $g$ are $f-$pairwise $c_r-$contractive over $(x_o,y_o)$ then $f$ and $g$ form a Cauchy pair over $(x_o,y_o)$.
\end{thm}
\textbf{Proof:} For all  $i$, let $x_i=f^i(x_o)$  and $y_i=g^i(y_o)$.
Let $$M=\max\{\vert P(\langle f(x_o)\rangle^{n-1},y_o)\vert,\vert P(\langle x_o\rangle^{n-1},y_o)\vert\}.$$ \indent First we prove that for all $i \ge j$,
$$r\le \min\{P(\langle x_i\rangle^n),P(\langle y_i\rangle^n)\}\le P(\langle x_{i}\rangle^{n-1},y_j)\le r+c^j\frac{nM}{1-c}.$$
From \Cref{D5c4} we know that
$$\begin{cases}r\le \min\{P(\langle x_i\rangle^n),P(\langle y_i\rangle^n)\}
\\ P(\langle x_{i+1}\rangle^{n-1},y_{i})\le r+ c^iM
\\P(\langle x_{i}\rangle^{n-1},y_i)\le r+ c^iM
\end{cases}\hspace{6ex} (\ddot\otimes)$$
In particular, for all $i$,
$$-P(\langle x_{i}\rangle^{n}) \le -r \text{ and }-P(\langle y_{i}\rangle^{n})\le -r.$$
By \Cref{T4c5} it suffices to bound  $P(\langle x_i \rangle^{n-1},y_j)$ for $i\ge j$. Let us first investigate what happens for the specific values of $i=6$ and $j=3$. By  repeatedly using ($P_n$-inq) we get
$$P(\langle x_{6}\rangle^{n-1},y_3)$$
$$\le P(\langle x_{6}\rangle^{n-1},y_5)-P(\langle y_{5}\rangle^{n}) +P(\langle y_{5}\rangle^{n-1},x_5)-P(\langle x_{5}\rangle^{n})+P(\langle x_{5}\rangle^{n-1},y_4)-P(\langle y_{4}\rangle^{n})$$
$$ +P(\langle y_{4}\rangle^{n-1},x_4)-P(\langle x_{4}\rangle^{n})+P(\langle x_{4}\rangle^{n-1},y_3)$$
by \Cref{C2e3}
$$\le P(\langle x_{6}\rangle^{n-1},y_5)-P(\langle y_{5}\rangle^{n}) +[(n-1)P(\langle x_{5}\rangle^{n-1},y_5)-(n-2)P(\langle x_{5}\rangle^{n})]-P(\langle x_{5}\rangle^{n})+P(\langle x_{5}\rangle^{n-1},y_4)-P(\langle y_{4}\rangle^{n})$$
$$ +[(n-1)P(\langle x_{4}\rangle^{n-1},y_4)-(n-2)P(\langle x_{4}\rangle^{n})]-P(\langle x_{4}\rangle^{n})+P(\langle x_{4}\rangle^{n-1},y_3)$$
by $(\ddot\otimes)$
$$\le (r+c^5M)-r+(n-1)(r+c^5M)+(n-2)(-r)-r+(r+c^4M)-r+(n-1)(r+c^4M)+(n-2)(-r)-r+(r+c^3M)$$
$$=r+n[c^5M+c^4M
]+c^3M$$
since $c^jM\ge0$
$$\le r+n[c^5M+c^4M+c^3M].$$
We now move to derive an upper bound for $P(\langle x_{i}\rangle^{n-1},y_j)$ where $i \ge j$.
\\\textit{Case 1:} If $i=j$, by $(\ddot\otimes)$ and since $0<c<1$ and $c^jM\ge 0$,
$$r \le \min\{P(\langle x_j\rangle^n),P(\langle y_j\rangle^n)\}\le P(\langle x_{j}\rangle^{n-1},y_j)\le r+ c^jM \le r+c^j\frac{nM}{1-c}.$$
\textit{Case 2:} If $i>j$, by $(\ddot\otimes)$ and ($P_n$-lbnd) we get
$$r\le \min\{P(\langle x_i\rangle^n),P(\langle y_j\rangle^n)\}\le P(\langle x_{i}\rangle^{n-1},y_j)$$
by repeatedly using ($P_n$-inq)
$$\le \sum\limits_{t=j+1}^{i-1}[P(\langle x_{t+1}\rangle^{n-1},y_t)-P(\langle y_{t}\rangle^{n}) +P(\langle y_{t}\rangle^{n-1},x_t)-P(\langle x_{t}\rangle^{n})]+P(\langle x_{j+1}\rangle^{n-1},y_j)$$
by \Cref{C2e3}
$$\le \sum\limits_{t=j+1}^{i-1}[P(\langle x_{t+1}\rangle^{n-1},y_t)-P(\langle y_{t}\rangle^{n}) +(n-1)P(\langle y_{t}\rangle^{n-1},x_t)-(n-2)P(\langle x_{t}\rangle^{n})-P(\langle x_{t}\rangle^{n})]+P(\langle x_{j+1}\rangle^{n-1},y_j)$$
by $(\ddot\otimes)$
$$\le \sum\limits_{t=j+1}^{i-1}[(r+c^tM)-r +(n-1)(r+c^tM)+(n-2)(-r)-r]+(r+c^jM)$$
$$= \sum\limits_{t=j+1}^{i-1}[nc^tM]+(r+c^jM)$$
since $c^jM\ge0$
$$\le r+nM\sum\limits_{t=j+1}^{i-1}[c^t]+nMc^j=r+nM\sum\limits_{t=j}^{i-1}c^t=r+nMc^j\sum\limits_{t=0}^{i-j-1}c^t\le r+nMc^j\sum\limits_{t=0}^{+\infty}c^t$$
finally, using the geometric series formula 
$$=r+c^j\frac{nM}{1-c}.$$
Hence, for all $i \ge j$,
$$r\le \min\{P(\langle x_i\rangle^n),P(\langle y_i\rangle^n)\}\le P(\langle x_{i}\rangle^{n-1},y_j)\le r+c^j\frac{nM}{1-c}.$$
\indent  For all positive real number $\epsilon$ there exists a natural number  $N$ such that $c^N\frac{nM}{1-c}<\epsilon$. Hence, for all $i \ge j >N$,\
$$r-\epsilon<r\le \min\{P(\langle x_i\rangle^n),P(\langle y_i\rangle^n)\}\le P(\langle x_{i}\rangle^{n-1},y_j)\le r+c^j\frac{nM}{1-c}<r+c^N\frac{nM}{1-c}<r+\epsilon.$$
Therefore, by \Cref{T4c5}, $f$ and $g$ are Cauchy pairs over $(x_o,y_o)$.$\hspace{6ex}\square$
\begin{defn}\label{D5c6} Let $(X,P)$ and $(Y,H)$ be two partial $n-\mathfrak{M}$etric spaces and suppose $f:X\to Y$ and $g:X \to Y$ are two functions on $X$. Let $r$, $A\ge 0$ and  $0<c<1$ be three real numbers. We say that $f$ and $g$ are \uline{\textbf{$f-$mutually $c_r-$contractive}} if and only if  for each $x$ in $X$ we can find an element $z$ in $X$ such that
$$H(\langle f(z)\rangle^{n-1},g(z))-H(\langle f(z) \rangle^n)\le c[H(\langle f(x)\rangle^{n-1},g(x))-H(\langle f(x) \rangle^n)]$$
and
$$r\le P(\langle z\rangle^n)\le P(\langle z\rangle^{n-1},x)\le r+A[H(\langle f(x)\rangle^{n-1},g(x))-H(\langle f(x) \rangle^n)].$$
\end{defn}
As in the case of partial metrics, the above definition  is used in coincidence point theorems of strong partial $n-\mathfrak{M}$etric. In the partial $n-\mathfrak{M}$etric case, a stronger contraction is needed.
\begin{defn}\label{D5c6-} Let $(X,P)$ and $(Y,H)$ be two partial $n-\mathfrak{M}$etric spaces and suppose $f:X\to Y$ and $g:X \to Y$ are two functions on $X$. Let $r$, $A\ge 0$ and  $0<c<1$ be three real numbers. We say that $f$ and $g$ are \uline{\textbf{$(f,g)-$mutually $c_r-$contractive}} if and only if  for each $x$ in $X$ we can find an element $z$ in $X$ such that
$$H(\langle f(z)\rangle^{n-1},g(z))-H(\langle f(z) \rangle^n)\le c[H(\langle f(x)\rangle^{n-1},g(x))-H(\langle f(x) \rangle^n)]$$
$$H(\langle f(z)\rangle^{n-1},g(z))-H(\langle g(z) \rangle^n)\le c[H(\langle f(x)\rangle^{n-1},g(x))-H(\langle g(x) \rangle^n)]$$
$$r\le P(\langle z\rangle^n)\le P(\langle z\rangle^{n-1},x)\le r+A[H(\langle f(x)\rangle^{n-1},g(x))-H(\langle f(x) \rangle^n)]$$
and
$$r\le P(\langle z\rangle^n)\le P(\langle z\rangle^{n-1},x)\le r+A[H(\langle f(x)\rangle^{n-1},g(x))-H(\langle g(x) \rangle^n)].$$
\end{defn} 
\begin{thm}{\textbf{(Cauchy $f-$Mutually $c_r-$Contractive):}}
\label{T5c2}
\\Let $(X,P)$ and $(Y,H)$ be two partial $n-\mathfrak{M}$etric spaces. Let $f:X\to Y$ and $g:X \to Y$ be two functions on $X$. If $f$ and $g$ are $f-$mutually $c_r-$contractive then there exists a Cauchy sequence $\{x_i\}_{i\in \mathbb{N}}$ in $X$ with central distance  $r$ such that for all natural numbers $i$, $r\le P(\langle x_i\rangle^n)$. Additionally for every positive real number $\epsilon$ there exists a natural number $N$ such that for all $i>N$,
$$H(\langle f(x_i)\rangle^{n-1},g(x_i))-H(\langle f(x_i) \rangle^n)<\epsilon.$$
\end{thm}
\textbf{Proof:} Since $f$ and $g$ are $f-$mutually $c_r-$contractive then there exists two real numbers $0<c<1$ and $A\ge 0$ such that for each $x$ in $X$, there exists a $z$ in $X$ where
$$H(\langle f(z)\rangle^{n-1},g(z))-H(\langle f(z) \rangle^n)\le c[H(\langle f(x)\rangle^{n-1},g(x))-H(\langle f(x) \rangle^n)]$$
and
$$r\le P(\langle z\rangle^n)\le P(\langle z\rangle^{n-1},x)\le r+A[H(\langle f(x)\rangle^{n-1},g(x))-H(\langle f(x) \rangle^n)].$$
Let us take an arbitrary element $x_1$ in $X$ then there exists an element $x_2$ in $X$ such that
$$H(\langle f(x_2)\rangle^{n-1},g(x_2))-H(\langle f(x_2) \rangle^n)\le c[H(\langle f(x_1)\rangle^{n-1},g(x_1))-H(\langle f(x_1) \rangle^n)]$$
and
$$r\le P(\langle x_2\rangle^n)\le P(\langle x_2\rangle^{n-1},x_1)\le r+A[H(\langle f(x_1)\rangle^{n-1},g(x_1))-H(\langle f(x_1) \rangle^n)].$$
There exists an element $x_3$ in $X$ such that
$$H(\langle f(x_3)\rangle^{n-1},g(x_3))-H(\langle f(x_3) \rangle^n)\le c[H(\langle f(x_2)\rangle^{n-1},g(x_2))-H(\langle f(x_2) \rangle^n)]$$
and
$$r\le P(\langle x_3\rangle^n)\le P(\langle x_3\rangle^{n-1},x_2)\le r+A[H(\langle f(x_2)\rangle^{n-1},g(x_2))-H(\langle f(x_2) \rangle^n)].$$
We continue the above process to generate a sequence $\{x_i\}_{i\in\mathbb{N}}$ such that for all $i$,
$$H(\langle f(x_{i+1})\rangle^{n-1},g(x_{i+1}))-H(\langle f(x_{i+1}) \rangle^n)\le c[H(\langle f(x_i)\rangle^{n-1},g(x_i))-H(\langle f(x_i) \rangle^n)]$$
and
$$r\le P(\langle x_{i+1}\rangle^n)\le P(\langle x_{i+1}\rangle^{n-1},x_{n})\le r+A[H(\langle f(x_{i})\rangle^{n-1},g(x_{i}))-H(\langle f(x_{i}) \rangle^n)]\hspace{6ex}(\check\otimes)$$
in particular, it is easy to see that for all $i$,
$$r\le P(\langle x_i\rangle^n)$$
i.e.
$$-P(\langle x_i\rangle^{n})\le-r.$$
\indent \uline{Step 1:} For all $i$ let $t_i=H(\langle f(x_i)\rangle^{n-1},g(x_i))-H(\langle f(x_i) \rangle^n)$. Then, from $(\check\otimes)$
$$t_{i+1}=H(\langle f(x_{i+1})\rangle^{n-1},g(x_{i+1}))-H(\langle f(x_{i+1}) \rangle^n)$$
$$\le c[H(\langle f(x_i)\rangle^{n-1},g(x_i))-H(\langle f(x_i) \rangle^n)]$$
and, hence, by induction
$$t_{i+1}\le c^iH(\langle f(x_1)\rangle^{n-1},g(x_1))-H(\langle f(x_1) \rangle^n).\hspace{6ex} (\hat\otimes)$$
Therefore, and since $0<c<1$ for every positive real number $\epsilon$ there exists  natural number $N'$ such that for all $i>N'$,
$$H(\langle f(x_i)\rangle^{n-1},g(x_i))-H(\langle f(x_i) \rangle^n)<\epsilon.$$
\indent \uline{Step 2:} For all $i>j$, by $(\check\otimes)$ and repeatedly using ($P_n$-inq) we get,
$$r\le P(\langle x_i\rangle^n)\le P(\langle x_i\rangle^{n-1},x_j)$$
$$\le P(\langle x_i\rangle^{n-1},x_{i-1})-P(\langle x_{i-1}\rangle^n)+P(\langle x_{i-1}\rangle^{n-1},x_{i-2})-P(\langle x_{i-2}\rangle^n)$$
$$+...+P(\langle x_{j+2}\rangle^{n-1},x_{j+1})-P(\langle x_{j+1}\rangle^n)+P(\langle x_{j+1}\rangle^{n-1},x_j)$$
$$\le \sum\limits_{k=j+1}^{i-1}[P(\langle x_{k+1}\rangle^{n-1},x_k)-P(\langle x_{k}\rangle^{n}) ]+P(\langle x_{j+1}\rangle^{n-1},y_j)$$
$$\le \sum\limits_{k=j+1}^{i-1}[r+A[H(\langle f(x_{k})\rangle^{n-1},g(x_{k}))-H(\langle f(x_{k}) \rangle^n)]-r ]+r+A[H(\langle f(x_{j})\rangle^{n-1},g(x_{j}))-H(\langle f(x_{j}) \rangle^n)]$$
$$=\sum\limits_{k=j+1}^{i-1}[At_k ]+r+A[t_j]=r+A\sum\limits_{k=j}^{i-1}t_k$$
by $(\hat\otimes)$
$$\le r+A\sum\limits_{k=j}^{i-1}c^{k-1}t_1=r+c^{j-1}At_1\sum\limits_{k=0}^{i-j-1}c^{k}$$
since $Ac^jt_k\ge 0$ and by the geometric series formula
$$\le r+c^{j-1}At_1\sum\limits_{k=0}^{+\infty}c^{k}=r+c^{j-1}\frac{At_1}{1-c}.$$
Therefore, for every positive real number $\epsilon$ there exists a natural number $N$ where $c^{N-1}\frac{At_1}{1-c}<\epsilon$ and, hence, for all $i \ge j>N$,
$$r-\epsilon <r\le P(\langle x_i\rangle^n)\le P(\langle x_i\rangle^{n-1},x_j)\le r+c^{j-1}\frac{At_1}{1-c} < r+c^{N-1}\frac{At_1}{1-c}<r+\epsilon.\hspace{6ex}\square$$
\indent The corollary below is straightforward by taking $X=Y$.
\begin{cor}\label{C5c1} Let $(X,P)$  be a partial $n-\mathfrak{M}$etric space. Let $f:X\to X$ and $g:X \to X$ be two functions on $X$. If $f$ and $g$ are $f-$mutually $c_r-$contractive then there exists a Cauchy sequence $\{x_i\}_{i\in \mathbb{N}}$ in $X$ with central distance $r$ such that for all natural numbers $i$, $r\le P(\langle x_i\rangle^n)$. Additionally for every positive real number $\epsilon$ there exists a natural number $N$ such that for all $i>N$,
$$P(\langle f(x_i)\rangle^{n-1},g(x_i))-P(\langle f(x_i) \rangle^n)<\epsilon.$$
\end{cor}
\chapter{Continuity and Non-Expansiveness}
\label{C6}
\setcounter{thm}{0}
\setcounter{defn}{0}
In  \Cref{C5}, we established some criteria on functions that are sufficient to generate  Cauchy sequences and Cauchy pairs.  Given that the limits (or special limits) of these sequences exist, we will need extra criteria on the functions for them to have a fixed point, common fixed point or coincidence point.
\begin{defn}\label{D6-1} Let $(X,\mathcal{T}_X)$ and $(Y,\mathcal{T}_Y)$ be two topological spaces. Let $f:X \to Y$ be a function on $X$. We say that $f$ is \uline{\textbf{continuous}} if and only if for every set $U$ open in $Y$, $f^{-1}(U)$ is open in $X$.
\end{defn}
\begin{defn}\label{D6-2} Let $(X,\mathcal{T}_X)$ and $(Y,\mathcal{T}_Y)$ be two topological spaces. Let $f:X \to Y$ be a function on $X$. We say that $f$ is \uline{\textbf{sequentially continuous}} if and only if for every sequence $\{x_i\}_{i \in \mathbb{N}}$ in $X$ having a limit $a$ in $X$, $f(a)$ is a limit of $\{f(x_i)\}_{i \in \mathbb{N}}$.
\end{defn}
\begin{thm}{(\textbf{Continuity vs. Sequential Continuity):}}
\label{T6-1}
\\Let $(X,\mathcal{T}_X)$ and $(Y,\mathcal{T}_Y)$ be two topological spaces. Let $f:X \to Y$ be a function on $X$. If $(X,\mathcal{T}_X)$ is first countable, then the the two statements below are equivalent.
\\(a) $f$ is continuous.
\\(b) $f$ is sequentially continuous. 
\end{thm}
\indent \indent The proof of \Cref{T6-1} is found in \cite{Mun2000}: Theorem 21.3. We mention this theorem since all topologies discussed in this thesis are first countable. In most cases we only need sequential continuity on the orbit, even less,  we only need sequential continuity for the special limit rather than for all limits on that orbit.
\begin{defn}\label{D6-3} Let $(X,\mathcal{T}_X)$ be a topological space with $x_o$ in $X$. Let $f:X \to X$ be a function on $X$. We say that $f$ is \uline{\textbf{orbitally continuous at} $x_o$} if and only if $a$ is a limit of $\{f^i(x_o)\}_{i \in \mathbb{N}}$ implies that $f( a)$ is a limit of $\{f^i(x_o)\}_{i \in \mathbb{N}}$.
\end{defn}
\indent \indent The definition of a special limit is not a topological definition, but rather a definition deduced from our generalized metrics. However, we will define weakly orbitally continuous functions now to avoid repeating the definition in each section. We will state the definition on a partial $n-\mathfrak{M}$etric space knowing that it includes all other cases discussed in this thesis.
\begin{defn}\label{D6-4} Let $(X,P)$ a partial $n-\mathfrak{M}$etric space with $x_o$ in $X$. Let  $f:X \to X$ be a  function on $X$. We say that $f$ is \uline{\textbf{weakly orbitally continuous at} $x_o$} if and only if $a$ is a special limit of $\{f^i(x_o)\}_{i \in \mathbb{N}}$ implies that $f( a)$ is a limit of $\{f^i(x_o)\}_{i \in \mathbb{N}}$.
\end{defn}
Notice that the difference between orbitally continuous and weakly orbitally continuous is that in the latter we can only guarantee that $f(a)$ is a limit of $\{f^i(x_o)\}_{i \in \mathbb{N}}$ if $a$ is a special limit of $\{f^i(x_o)\}_{i \in \mathbb{N}}$. \Cref{T6-1} shows that continuity and sequential continuity become equivalent notions in first countable spaces. Hence, in first countable spaces, 
$$\text{Continuous}\iff \text{Sequentially Continuous}\Rightarrow\text{Orbitally Continuous at }x_o\Rightarrow \text{ Weakly Orbitally Continuous at }x_o.$$
The two notions left for us to define are non-expansiveness and consistency. We will present the definition in the partial $n-\mathfrak{M}$etric case as it is our most general case.
\begin{defn}\label{D6c1}  Let $(X,P)$ a partial $n-\mathfrak{M}$etric  space. Let  $f:X \to X$ be a function on $X$. We say that $f$ is \uline{\textbf{non-expansive}} if and only if for every two elements $x$ and $y$ in $X$,
$$P(\langle f(x)\rangle^{n-1},f(y))\le P(\langle x\rangle^{n-1},y).$$
\end{defn}
\begin{defn}\label{D6c2} Let $(X,P)$ and $(Y,H)$ two partial $n-\mathfrak{M}$etric spaces. Let $f:X \to Y$ be a function on $X$. We say that $f$ is \uline{\textbf{consistent}} if and only if for every two elements $x$ and $z$ in $X$,
$$P(\langle x\rangle^{n})\le P(\langle z\rangle^{n})$$
implies
$$H(\langle f( x)\rangle^{n}) \le H(\langle f( z)\rangle^{n}).$$
\end{defn}
\section{Metric Space}
\label{C6a}
\indent \indent In the case of a metric space, the definition of non-expansive functions given in \Cref{D6c1} is written as
$$d(f(x),f(y))\le d(x,y).$$  The Lemmas presented in this section are folklore.  Hence, we will state them while providing a minimal proof when needed.
\begin{lem}\label{L6a1} Let $(X,d)$ be a metric  space. If $f:X \to X$ is a non-expansive function on $X$ then $f$ is continuous.
\end{lem}
\begin{lem}\label{L6a2} Let $(X,d)$ be a metric  space with $x_o$ in $X$. Let  $f:X \to X$ be a non-expansive function on $X$. If $a$ is a limit of $\{f^i(x_o)\}_{i \in \mathbb{N}}$, then $f(a)=a.$ $\hspace{6ex}\square$
\end{lem}
\textbf{Proof:} A non-expansive function on a metric space is continuous and, hence,  by \Cref{T6-1} sequentially continuous.
Therefore, by \Cref{D6-2} $f(a)=a$. $\hspace{6ex}\square$
\\\indent In a metric space the definitions of orbitally continuous and weakly orbitally continuous coincide since the definitions of special limits and limits coincide.
\begin{lem}\label{L6a3} Let $(X,d)$ be a metric  space with $x_o$ in $X$ . Let  $f:X \to X$ be a weakly orbitally continuous function at $x_o$. If $a$ is a limit of $\{f^i(x_o)\}_{i \in \mathbb{N}}$, then $f(a)=a.$
\end{lem}
\textbf{Proof:}  In a metric space the limit is unique and, hence, by \Cref{D6-2} $f(a)=a$.$\hspace{6ex}\square$
\section{Partial Metric Space}
\label{C6b}
\indent \indent As in \Cref{C5a} the results in this section are special cases of results in \Cref{C6c}. We include the proofs because they are much simpler than those of the more general results.
 In the case of a partial metric space, the definition of non-expansive functions given in \Cref{D6c1} is written as
$$p(f(x),f(y))\le p(x,y).$$
\begin{rmk} In a metric space, a non-expansive function is continuous and, hence, weakly orbitally continuous. Additionally in the metric case, as pointed out in \Cref{C6a}, the notions of orbital continuity and weak orbital continuity coincide. On the other hand, in a partial metric space, a non-expansive function need not be continuous or even weakly orbitally continuous. Moreover, a weakly orbitally continuous function need not be orbitally continuous. We show these important differences using the three examples below.
\end{rmk} 
\begin{eg}{\textbf{(Non-Expansiveness vs Continuity):}}
\label{E6b1}
\\Let $p$ be a partial metric on $X=\mathbb{R}\cup\{a\}$ where $a\notin\mathbb{R}$ as defined in  \Cref{E2b3} by:
\\For all $x,y\in \mathbb{R}$,
$$p(a,a)=0,p(a,x)=\vert x \vert \text{ and }p(x,y)=\vert x-y \vert -1.$$
Let 
$$f(x)=\begin{cases}x & \text{ if } x \in \mathbb{R}.
\\1 & \text{ if }x=a. \\
\end{cases}$$
\indent The function $f$ is non-expansive over $\mathbb{R}$ since $$ p(f(a),f(a))=p(1,1)=-1\le 0=p(a,a).$$
Additionally, for each $x\in\mathbb{R}$,
$$p(f(a),f(x))=p(1,x)=\vert x-1\vert-1\le \vert x \vert =p(a,x).$$
  On the other hand, $a$ is a limit of $\{\frac{1}{2^n}\}_{n \in \mathbb{N}}$ but $f(a)=1$ is not. Hence, $f$ is not sequentially continuous and by \Cref{T6-1} is not continuous.
\end{eg}
\newpage
\begin{eg}{\textbf{(Non-Expansiveness vs. Weak Orbital Continuity):}}
\label{E6b2}
\\Let $p:X \times X\to \mathbb{R}$ be a partial metric on $X=[-1,1]$ as defined in \Cref{E2b2} by:
\\For all $x,y\in \mathbb{R}$,
$$p(x,y)=\max\{x,y\}.$$
Let
$$f(x)=\begin{cases}\frac{x}{2} & \text{ if } x \ne 0.
\\-1 & \text{ if }x=0. \\
\end{cases}$$ 
For each $x\in[-1,0)\cup(0,1]$
$$p(f(0),f(x))=p(-1,\frac{x}{2})=\max\{-1,\frac{x}{2}\}\le \max\{0,x\}=p(0,x).$$
Hence, $f$ is non-expansive. Showing that $0$ is a special limit of $\{f^n(1)\}_{n\in \mathbb{N}}=\{\frac{1}{2^n}\}_{n \in \mathbb{N}}$ is left as an exercise to the reader. On the other hand, $f(0)=-1$ is not a limit of $\{\frac{1}{2^n}\}_{n \in \mathbb{N}}$ since for each $n\in \mathbb{N}$ and $\epsilon<1$,
$$p(-1,\frac{1}{2^n})-p(-1,-1)=\frac{1}{2^n}+1>1>\epsilon.$$
\end{eg} 
\begin{eg}{\textbf{(Weak vs Usual Orbital Continuity):}}
\label{E6-1}
\\Let $p$ be a partial metric on $X=\mathbb{R}\cup\{a\}$ where $a\notin\mathbb{R}$ as defined in  \Cref{E2b3} by
$$p(x,y)=\begin{cases}0 & \text{ if } x=y=a. \\
\vert y \vert &\text{ if } x=a \text{ and } y\in\mathbb{R}.\\
\vert x-y \vert -1& \text{ if } \{x,y\}\subseteq\mathbb{R}. \\
\end{cases}$$
Let 
$$f(x)=\begin{cases}\frac{x}{2} & \text{ if } x \in \mathbb{R}-\{0\}. \\
a &\text{ if } x=0.\\
5& \text{ if }x=a. \\
\end{cases}$$
 Then, the sequence $\{f^n(1)\}_{n\in \mathbb{N}}=\{\frac{1}{2^n}\}_{n \in \mathbb{N}}$. As shown in \Cref{E2b3}, $0$ is a special limit and $a$ is a limit of the sequence $\{\frac{1}{2^n}\}_{n \in \mathbb{N}}$. Moreover, $f(0)=a$ is a limit (not a special limit though) of $\{\frac{1}{2^n}\}_{n \in \mathbb{N}}$ whereas $f(a)=5$ is not. Hence, $f$ is weakly orbitally continuous (but not orbitally continuous) at $x_o=1$.
\end{eg}
\begin{lem}\label{L6b1} Let $(X,p)$ a partial metric  space with $x_o$ in $X$. Let  $f:X \to X$ be a non-expansive function on $X$. If $a$ is a special limit of $\{f^i(x_o)\}_{i \in \mathbb{N}}$ then $p(a,f(a))=p(a,a)$.
\end{lem}
\textbf{Proof:}  Since $a$ is a special limit of $\{f^i(x_o)\}_{i \in \mathbb{N}}$ (see \Cref{D4b2}) this sequence is a Cauchy sequence with central distance $r=p(a,a)$. From (p-lbnd) we know that
$$p(a,a)\le p(a,f(a))$$
by (p-inq)
$$\le p(f(a), f^{i+1}(x_o))+p(f^{i+1}(x_o),a)-p(f^{i+1}(x_o),f^{i+1}(x_o))$$
by (p-sym)
$$=p(f(a), f^{i+1}(x_o))+p(a,f^{i+1}(x_o))-p(f^{i+1}(x_o),f^{i+1}(x_o))$$
since $f$ is non-expansive
$$\le p(a, f^{i}(x_o))+p(a,f^{i+1}(x_o))-p(f^{i+1}(x_o),f^{i+1}(x_o)).$$
For every positive real number $\epsilon$ by \Cref{D4b1} there exists a natural number $N_1$ such that for all $i>N_1$
$$-p(f^{i+1}(x_o),f^{i+1}(x_o))<-r+\frac{\epsilon}{3}=-p(a,a)+\frac{\epsilon}{3}$$
and, since a special limit is a limit, by \Cref{L4b1} there exists a natural number $N_2$ such that for all $i>N_2$,
$$p(a,f^i(x_o))<p(a,a)+\frac{\epsilon}{3}.$$
Therefore, for every positive real number $\epsilon$ there exists a natural number $N=\max\{N_1,N_2\}$ such that for all $i>N$,
$$p(a,a)\le p(a,f(a))\le  p(a, f^{i}(x_o))+p(a,f^{i+1}(x_o))-p(f^{i+1}(x_o),f^{i+1}(x_o))$$
$$<p(a,a)+ \frac{\epsilon}{3}+p(a,a)+\frac{\epsilon}{3}-p(a,a)+\frac{\epsilon}{3}=p(a,a)+\epsilon.$$
Hence, $p(a,a)=p(a,f(a))$.$\hspace{6ex}\square$
\begin{lem}\label{L6b2} Let $(X,p)$ a partial metric  space with $x_o \in X$. Let  $f:X \to X$ be a weakly orbitally continuous function at $x_o$. If $a$ is a special limit of $\{f^i(x_o)\}_{i \in \mathbb{N}}$ then $p(a,f(a))= p(f(a),f(a)).$
\end{lem}
\textbf{Proof:} From \Cref{D4b2}, since $a$ is a special limit of $\{f^i(x_o)\}_{i \in \mathbb{N}}$ then that sequence is a Cauchy sequence with central distance $r=p(a,a)$. From (p-lbnd) we know that
$$p(f(a),f(a))\le p(a,f(a))$$
by (p-inq)
$$\le p(f(a), f^{i+1}(x_o))+p(f^{i+1}(x_o),a)-p(f^{i+1}(x_o),f^{i+1}(x_o))$$
by (p-sym)
$$=p(f(a), f^{i+1}(x_o))+p(a,f^{i+1}(x_o))-p(f^{i+1}(x_o),f^{i+1}(x_o)).$$
For every positive real number $\epsilon$ by \Cref{D4b1} there exists a natural number $N_1$ such that for all $i>N_1$
$$-p(f^{i+1}(x_o),f^{i+1}(x_o))<-r+\frac{\epsilon}{3}=-p(a,a)+\frac{\epsilon}{3}.$$
Since $a$ is a limit of $\{f^i(x_o)\}_{i \in \mathbb{N}}$ by \Cref{L4b1} there exists a natural number $N_2$ such that for all $i>N_2$,
$$p(a,f^i(x_o)<p(a,a)+\frac{\epsilon}{3}.$$
\indent Furthermore, $f$ is weakly orbitally continuous at $x_o$ and $a$ is a special limit of $\{f^i(x_o)\}_{i \in \mathbb{N}}$ hence, $f(a)$ is a limit of $\{f^i(x_o)\}_{i \in \mathbb{N}}$. Therefore, from \Cref{L4b1} there exists a natural number $N_3$ such that for all $i>N_3$,
$$p(f(a),f^i(x_o))<p(f(a),f(a))+\frac{\epsilon}{3}.$$
Hence, for every positive real number $\epsilon$ there exists a natural number $N=\max\{N_1,N_2,N_3\}$ such that for all $i>N$,
$$p(f(a),f(a))\le p(a,f(a)\le p(f(a), f^{i+1}(x_o))+p(f^{i+1}(x_o),a)-p(f^{i+1}(x_o),f^{i+1}(x_o))$$
$$<p(f(a),f(a))+\frac{\epsilon}{3}+p(a,a)+\frac{\epsilon}{3}-p(a,a)=p(f(a),f(a))+\epsilon.$$
Therefore, $p(f(a),f(a))=p(a,f(a))$.$\hspace{6ex}\square$
\section{Partial $n-\mathfrak{M}$etric space}
\label{C6c}
\indent \indent As in \Cref{C5c}, the proofs in \Cref{C6c} are quite similar to those in \Cref{C6b} aside form the need to change $\epsilon$ to fit our needs.  We still present the proofs for completeness.
\begin{lem}\label{L6c1} Let $(X,P)$ a partial $n-\mathfrak{M}$etric  space with $x_o$ in $X$. Let  $f:X \to X$ be a non-expansive function on $X$. If $a$ is a special limit of $\{f^i(x_o)\}_{i \in \mathbb{N}}$ then 
$$P(\langle a\rangle^{n-1},f(a))=P(\langle a\rangle^{n})$$
and
$$P(\langle f(a)\rangle^{n-1},a)\le P(\langle a\rangle^{n}).$$
\end{lem}
\textbf{Proof:} From \Cref{D4c2}, since $a$ is a special limit of $\{f^i(x_o)\}_{i \in \mathbb{N}}$ then this sequence is a Cauchy sequence with central distance $r=P(\langle a\rangle^{n})$.  From ($P_n$-lbnd) we know that
$$P(\langle a\rangle^{n})\le P(\langle  a\rangle^{n-1},f(a))$$
by ($P_n$-inq)
$$\le P(\langle a\rangle^{n-1},f^{i+1}(x_o))+ P(\langle f^{i+1}(x_o)\rangle^{n-1},f(a))-P(\langle f^{i+1}(x_o)\rangle^{n})$$
since $f$ is non-expansive
$$\le P(\langle a\rangle^{n-1},f^{i+1}(x_o))+ P(\langle f^{i}(x_o)\rangle^{n-1},a)-P(\langle f^{i+1}(x_o)\rangle^{n})$$
and by \Cref{C2e3}
$$\le P(\langle a\rangle^{n-1},f^{i+1}(x_o))+ (n-1)P(\langle a\rangle^{n-1},f^{i}(x_o))-(n-2)P(\langle a\rangle^{n})-P(\langle f^{i+1}(x_o)\rangle^{n}).$$
By \Cref{D4c1}, for every positive real number $\epsilon$ there exists a natural number $N_1$ such that for all $i>N_1$,
$$-P(\langle f^{i}(x_o)\rangle^{n})<-r+\frac{\epsilon}{n+1}=-P(\langle a\rangle^{n})+ \frac{\epsilon}{n+1}.$$
Since a special limit is a limit, by \Cref{L4c1} there exists a natural number $N_2$ such that for all $i>N_2$,
$$P(\langle a\rangle^{n-1},f^{i}(x_o))<P(\langle a\rangle^{n})+ \frac{\epsilon}{n+1}.$$
Therefore, for every positive real number $\epsilon$ there exists a natural number $N=\max\{N_1,N_2,N_3\}$ such that for all $i>N$,
$$P(\langle a\rangle^{n})\le P(\langle a\rangle^{n-1},f(a))$$
$$\le P(\langle a\rangle^{n-1},f^{i+1}(x_o))+ (n-1)P(\langle a\rangle^{n-1},f^{i}(x_o))-(n-2)P(\langle a\rangle^{n})-P(\langle f^{i+1}(x_o)\rangle^{n})$$
$$<P(\langle a\rangle^{n})+ \frac{\epsilon}{n+1}+(n-1)(P(\langle a\rangle^{n})+ \frac{\epsilon}{n+1})+(n-2)P(\langle a\rangle^{n})-P(\langle a\rangle^{n})+ \frac{\epsilon}{n+1}$$
$$= P(\langle a\rangle^{n})+(n+1)\frac{\epsilon}{n+1}=P(\langle  a\rangle^{n})+\epsilon.$$
Hence, $P(\langle  a\rangle^{n})= P(\langle  a\rangle^{n-1},f(a)).$
\\Similarly by ($P_n$-inq)
$$P(\langle f(a)\rangle^{n-1},a)\le P(\langle f(a)\rangle^{n-1},f^{i+1}(x_o))+P(\langle f^{i+1}(x_o)\rangle^{n-1},a)-P(\langle f^{i+1}(x_o)\rangle^{n})$$
since $f$ is non-expansive
$$\le P(\langle a\rangle^{n-1},f^{i}(x_o))+P(\langle f^{i+1}(x_o)\rangle^{n-1},a)-P(\langle f^{i+1}(x_o)\rangle^{n})$$
by \Cref{C2e3}
$$\le P(\langle a\rangle^{n-1},f^{i}(x_o))+(n-1)P(\langle  a\rangle^{n-1},f^{i+1}(x_o))-(n-2)P(\langle  a\rangle^{n})-P(\langle f^{i+1}(x_o)\rangle^{n}).$$
Hence, for every positive real number $\epsilon$ there exists a natural number $N$ such that for all $i>N$,
$$P(\langle f( a)\rangle^{n-1},a)$$
$$\le P(\langle a\rangle^{n-1},f^{i}(x_o))+(n-1)P(\langle  a\rangle^{n-1},f^{i+1}(x_o))-(n-2)P(\langle  a\rangle^{n})-P(\langle f^{i+1}(x_o)\rangle^{n})$$
$$< P(\langle a\rangle^{n})+ \frac{\epsilon}{n+1}+(n-1)(P(\langle  a\rangle^{n})+ \frac{\epsilon}{n+1})-(n-2)P(\langle a\rangle^{n})-P(\langle a\rangle^{n})+ \frac{\epsilon}{n+1}$$
$$=P(\langle  a\rangle^{n})+(n+1) \frac{\epsilon}{n+1}=P(\langle  a\rangle^{n})+ \epsilon.$$
Therefore, $P(\langle f(a)\rangle^{n-1},a)\le P(\langle  a\rangle^{n}).\hspace{6ex}\square$
\begin{lem}\label{L6c2} Let $(X,P)$ a partial $n-\mathfrak{M}$etric  space with $x_o$ in $X$. Let  $f:X \to X$ be weakly orbitally continuous function at $x_o$. If $a$ is a special limit of $\{f^i(x_o)\}_{i \in \mathbb{N}}$ then 
$$P(\langle f(a)\rangle^{n-1},a)=P(\langle f( a)\rangle^{n})$$
and
$$P(\langle a\rangle^{n-1},f(a))\le P(\langle f( a)\rangle^{n}).$$
\end{lem}
\textbf{Proof:} From \Cref{D4c2}, since $a$ is a special limit of $\{f^i(x_o)\}_{i \in \mathbb{N}}$ then that sequence is a Cauchy sequence with central distance $r=P(\langle a\rangle^{n})$.  From ($P_n$-lbnd) we know that
$$P(\langle f(a)\rangle^{n})\le P(\langle f( a)\rangle^{n-1},a)$$
by ($P_n$-inq)
$$\le P(\langle f(a)\rangle^{n-1},f^{i+1}(x_o))+ P(\langle f^{i+1}(x_o)\rangle^{n-1},a)-P(\langle f^{i+1}(x_o)\rangle^{n})$$
by \Cref{C2e3}
$$\le P(\langle f(a)\rangle^{n-1},f^{i+1}(x_o))+ (n-1)P(\langle a\rangle^{n-1},f^{i+1}(x_o))-(n-2)P(\langle a\rangle^{n})-P(\langle f^{i+1}(x_o)\rangle^{n}).$$
By \Cref{D4c1}, for every positive real number $\epsilon$ there exists a natural number $N_1$ such that for all $i>N_1$,
$$-P(\langle f^{i}(x_o)\rangle^{n})<-r+\frac{\epsilon}{2n-1}=-P(\langle a\rangle^{n})+ \frac{\epsilon}{n+1}.$$
Since a special limit is a limit, by \Cref{L4c1} there exists a natural number $N_2$ such that for all $i>N_2$,
$$P(\langle a\rangle^{n-1},f^{i}(x_o))<P(\langle a\rangle^{n})+ \frac{\epsilon}{n+1}.$$
Since $f$ is weakly orbitally continuous at $x_o$ then $f(a)$ is also a limit of $\{f^i(x_o)\}_{i \in \mathbb{N}}$ then there exists a natural number $N_3$ such that for all $i>N_3$,
$$P(\langle f(a)\rangle^{n-1},f^{i}(x_o))<P(\langle f(a)\rangle^{n})+ \frac{\epsilon}{n+1}.$$
Therefore, $N=\max\{N_1,N_2,N_3\}$ is a natural number such that for all $i>N$,
$$P(\langle f(a)\rangle^{n})\le P(\langle f(a)\rangle^{n-1},a)$$
$$\le P(\langle f(a)\rangle^{n-1},f^{i+1}(x_o))+ (n-1)P(\langle a\rangle^{n-1},f^{i+1}(x_o))-(n-2)P(\langle a\rangle^{n})-P(\langle f^{i+1}(x_o)\rangle^{n})$$
$$<P(\langle f(a)\rangle^{n})+ \frac{\epsilon}{n+1}+(n-1)(P(\langle a\rangle^{n})+ \frac{\epsilon}{n+1})+(n-2)P(\langle a\rangle^{n})-P(\langle a\rangle^{n})+ \frac{\epsilon}{n+1}$$
$$= P(\langle f(a)\rangle^{n})+(n+1)\frac{\epsilon}{n+1}=P(\langle f( a)\rangle^{n})+\epsilon.$$
Hence, $P(\langle f( a)\rangle^{n})= P(\langle f( a)\rangle^{n-1},a).$
\\Similarly by ($P_n$-inq)
$$P(\langle a\rangle^{n-1},f(a))\le P(\langle a\rangle^{n-1},f^{i+1}(x_o))+P(\langle f^{i+1}(x_o)\rangle^{n-1},f(a))-P(\langle f^{i+1}(x_o)\rangle^{n})$$
by \Cref{C2e3}
$$\le P(\langle a\rangle^{n-1},f^{i+1}(x_o))+(n-1)P(\langle f( a)\rangle^{n-1},f^{i+1}(x_o))-(n-2)P(\langle f( a)\rangle^{n})-P(\langle f^{i+1}(x_o)\rangle^{n}).$$
Hence, for every positive real number $\epsilon$ there exists a natural number $N$ such that for all $i>N$,
$$P(\langle a\rangle^{n-1},f(a))$$
$$\le P(\langle a\rangle^{n-1},f^{i+1}(x_o))+(n-1)P(\langle f( a)\rangle^{n-1},f^{i+1}(x_o))-(n-2)P(\langle f( a)\rangle^{n})-P(\langle f^{i+1}(x_o)\rangle^{n})$$
$$< P(\langle a\rangle^{n})+ \frac{\epsilon}{n+1}+(n-1)(P(\langle f( a)\rangle^{n})+ \frac{\epsilon}{n+1})-(n-2)P(\langle f( a)\rangle^{n})-P(\langle a\rangle^{n})+ \frac{\epsilon}{n+1}$$
$$=P(\langle f( a)\rangle^{n})+(n+1) \frac{\epsilon}{n+1}=P(\langle f( a)\rangle^{n})+ \epsilon.$$
Therefore, $P(\langle a\rangle^{n-1},f(a))\le P(\langle f( a)\rangle^{n}).\hspace{6ex}\square$
\chapter{Applications to Fixed point and Coincidence Point Theory.}
\label{C7}
\setcounter{thm}{0}
\setcounter{defn}{0}
We have reached the end of the rainbow to find our pot of gold. In this chapter we state fixed, common fixed and coincidence point theorems whose sole constraint on the generalized metric spaces is that they be complete. This thesis was intentionally written in a way that minimizes the proofs in this section. The theorems and lemmas in previous chapters are building blocks for the theorems ahead. As previously stated, a (strong) partial $n-\mathfrak{M}$etric is a generalization  of a (strong) partial metric. Any special technique needed for the (strong) partial $n-\mathfrak{M}$etric case has already been presented in previous chapters. That is why we will be omitting the proofs of \Cref{C7d} and \Cref{C7e} to spare the reader any redundancy. We start with some basic definitions.
\begin{defn}\label{D7-1} Let $X$ be a non-empty set. Let $f:X \to X$ be a function on $X$. We say that $x$ in $X$ is a \textbf{\uline{fixed point of $f$}} if and only if $f(x)=x$.
\end{defn}
\begin{defn}\label{D7-2} Let $X$ be a non-empty set. Let $f:X \to X$ and $g:X \to X$ be two functions on $X$. We say that $x$ in $X$ is a \textbf{\uline{common fixed point of $f$ and $g$}} if and only if $f(x)=x=g(x)$.
\end{defn}
\begin{defn}\label{D7-3} Let $X$ and $Y$ be two non-empty sets. Let $f:X \to Y$ and $g:X \to Y$ be two functions on $X$. We say that $x$ in $X$ is a \textbf{\uline{coincidence point of $f$ and $g$}} if and only if $f(x)=g(x)$.
\end{defn}
\section{Metric Space}
\label{C7a}
\indent \indent Depending on the type of the contractive function used, \Cref{T7a1} and \Cref{T7a2} can be attributed to either Edelstein \cite{Ede1962,Ede1964,Ede1961} or Alber and  Guerre-Delabriere \cite{Alb1997}.
\begin{thm}{\textbf{(Fixed point and Non-expansive):}}
\label{T7a1}
\\Let $(X,d)$ be a complete metric space with $x_o$ in $X$. Let $f:X \to X$ be a Cauchy function at $x_o$. If $f$ is non-expansive then $f$ has a fixed point in $X$.
\end{thm}
\textbf{Proof:} Since $f$ is Cauchy at $x_o$ then by \Cref{D5-2}, $\{f^i(x_o)\}_{i \in \mathbb{N}}$ is a Cauchy sequence.
Since $(X,d)$ is a complete metric space then by \Cref{D4a2} $\{f^i(x_o)\}_{i \in \mathbb{N}}$ has a  limit $a$ in $X$. Finally $f$ is non-expansive, then by \Cref{L6a1} $f(a)=a$ and, hence, by \Cref{D7-3} $a$ is a fixed point of $f.\hspace{6ex}\square$
\\\indent In fact, \Cref{T7a1} can be considered a corollary of \Cref{T7a2} since any non-expansive function in a metric space is continuous. \begin{thm}{\textbf{(Fixed point and Weak orbital continuity):}}
\label{T7a2}
\\Let $(X,d)$ be a complete metric space with $x_o$ in $X$. Let $f:X \to X$ a Cauchy function at $x_o$. If $f$ is  weakly orbitally continuous  at $x_o$ then $f$ has a fixed point in $X$.
\end{thm}
\textbf{Proof:} Since $f$ is Cauchy at $x_o$ then by \Cref{D5-2}, $\{f^i(x_o)\}_{i \in \mathbb{N}}$ is a Cauchy sequence.
Since $(X,d)$ is a complete metric space then by \Cref{D4a2} $\{f^i(x_o)\}_{i \in \mathbb{N}}$ has a  limit $a$ in $X$. Finally $f$ is weakly orbitally continuous at $x_o$, then by \Cref{L6a2} $f(a)=a$ and, hence, by \Cref{D7-3} $a$ is a fixed point of $f.\hspace{6ex}\square$
\begin{rmk}\label{R7a1} From \Cref{L5a1} and \Cref{L5a2}, if $f$ is orbitally $c_0-$contractive or orbitally $\varphi_0-$contractive at $x_o$ then $f$ is Cauchy at $x_o$.
\end{rmk}
\begin{thm}{\textbf{(Common fixed point and Non-expansive):}}
\label{T7a3}
\\Let $(X,d)$ be a complete metric space with $x_o$ and $y_o$ in $X$. Let  $f:X \to X$ and $g: X \to X$ be two functions that form a Cauchy pair over $(x_o,y_o)$. If $f$ and $g$ are non-expansive  then $f$ and $g$ have a common fixed point.
\end{thm}
\textbf{Proof:} Since $f$ and $g$  form a Cauchy pair at $(x_o,y_o)$ then by \Cref{D5-3} $\{f^i(x_o)\}_{i\in\mathbb{N}}$ and $\{g^i(y_o)\}_{i\in\mathbb{N}}$ form a Cauchy pair. By \Cref{L4a3} $\{f^i(x_o)\}_{i\in\mathbb{N}}$ and $\{g^i(y_o)\}_{i\in\mathbb{N}}$ are both Cauchy sequences. Since $(X,d)$ is a complete metric space then by \Cref{D4a2}  and \Cref{L4a3} $\{f^i(x_o)\}_{i\in\mathbb{N}}$ and $\{g^i(y_o)\}_{i\in\mathbb{N}}$ both have the same limit  $a$ in $X$. Finally $f$ and $g$ are both non-expansive, then by \Cref{L6a1} $f(a)=a=g(a)$ and, hence, by \Cref{D7-2} $a$ is a common fixed point of $f$ and $g.\hspace{6ex}\square$
\begin{thm}{\textbf{(Common fixed point and Weak orbital continuity):}}
\label{T7a4}
\\Let $(X,d)$ be a complete metric space with $x_o$ and $y_o$ in $X$. Let  $f:X \to X$ and $g: X \to X$ be two functions that form a Cauchy pair over $(x_o,y_o)$. If $f$ and $g$ are    weakly orbitally continuous at $x_o$ and $y_o$ respectively then $f$ and $g$ have a common fixed point.
\end{thm}
\textbf{Proof:} Since $f$ and $g$  form a Cauchy pair at $(x_o,y_o)$ then by \Cref{D5-3} $\{f^i(x_o)\}_{i\in\mathbb{N}}$ and $\{g^i(y_o)\}_{i\in\mathbb{N}}$ form a Cauchy pair. By \Cref{L4a3} $\{f^i(x_o)\}_{i\in\mathbb{N}}$ and $\{g^i(y_o)\}_{i\in\mathbb{N}}$ are both Cauchy sequences. Since $(X,d)$ is a complete metric space then by \Cref{D4a2}  and \Cref{L4a3} $\{f^i(x_o)\}_{i\in\mathbb{N}}$ and $\{g^i(y_o)\}_{i\in\mathbb{N}}$ both have the same limit  $a$ in $X$. Finally $f$ and $g$ are weakly orbitally continuous at $x_o$ and $y_o$ respectively, then by \Cref{L6a2} $f(a)=a=g(a)$ and, hence, by \Cref{D7-2} $a$ is a common fixed point of $f$ and $g.\hspace{6ex}\square$
\begin{thm}{\textbf{(Common fixed point and Mixed criteria):}}
\label{T7a5}
\\Let $(X,d)$ be a complete metric space with $x_o$ and $y_o$ in $X$. Let  $f:X \to X$ and $g: X \to X$ be two functions that form a Cauchy pair over $(x_o,y_o)$. If $f$ is non-expansive and $g$ is weakly orbitally continuous at $y_o$ then $f$ and $g$ have a common fixed point.
\end{thm}
\textbf{Proof:} Since $f$ and $g$  form a Cauchy pair at $(x_o,y_o)$ then by \Cref{D5-3} $\{f^i(x_o)\}_{i\in\mathbb{N}}$ and $\{g^i(y_o)\}_{i\in\mathbb{N}}$ form a Cauchy pair. By \Cref{L4a3} $\{f^i(x_o)\}_{i\in\mathbb{N}}$ and $\{g^i(y_o)\}_{i\in\mathbb{N}}$ are both Cauchy sequences. Since $(X,d)$ is a complete metric space then by \Cref{D4a2} and \Cref{L4a3} $\{f^i(x_o)\}_{i\in\mathbb{N}}$ and $\{g^i(y_o)\}_{i\in\mathbb{N}}$ both have the same limit  $a$ in $X$. Finally $f$ is non-expansive then by \Cref{L6a1} $f(a)=a$ and $g$ is weakly orbitally continuous at  $y_o$  then by \Cref{L6a2} $a=g(a)$ and, hence, by \Cref{D7-2} $a$ is a common fixed point of $f$ and $g.\hspace{6ex}\square$
\begin{rmk}\label{R7a2} From \Cref{T5a1}, if $f$ and $g$ are $f-$pairwise $c_0-$contractive (similarly $g-$pairwise $c_0-$contractive) over $(x_o,y_o)$ then $f$ and $g$ form a Cauchy pair over $(x_o,y_o)$.
\end{rmk}
\begin{thm}{\textbf{(Coincidence Point Theorem):}}
\label{T7a6}
\\Let $(X,l)$ be a complete metric space and let $(Y,d)$ be a metric space. Let $f:X \to Y$ and $g: X \to Y$ be two sequentially continuous functions on $X$. If $f$ and $g$ are mutually $c_0-$contractive then $f$ and $g$ have a coincidence point.
\end{thm}
\textbf{Proof:} Since $f$ and $g$ are mutually $c_0-$contractive then by \Cref{T5a2} there exists a Cauchy sequence $\{x_i\}_{i \in \mathbb{N}}$ in $X$ where $\{f(x_i)\}_{i \in \mathbb{N}}$ and $\{g(x_i)\}_{i \in \mathbb{N}}$ form a Cauchy pair in $Y$. Hence, by \Cref{D4a3} for every positive real number $\epsilon$ there exists a natural number $N_1$ such that for all $i>N_1$,
$$d(f(x_i),g(x_i))<\frac{\epsilon}{3}.$$
Since $(X,l)$ is complete  then $\{x_i\}_{i \in \mathbb{N}}$ has a limit $a$ in $X$. Since $f$ and $g$   are sequentially continuous then by \Cref{D6-2} $f(a)$ and $g(a)$ are limits of $\{f(x_i)\}_{i \in \mathbb{N}}$ and $\{g(x_i)\}_{i \in \mathbb{N}}$ respectively.
\\$f(a)$ is the  limit of  $\{f(x_i)\}_{i \in \mathbb{N}}$ hence, by \Cref{L4a1} there exists a natural number $N_2$ such that for all $i>N_2$,
$$d(f(a),f(x_i))<\frac{\epsilon}{3}$$
$g(a)$ is the  limit of  $\{g(x_i)\}_{i \in \mathbb{N}}$ hence, by \Cref{L4a1} there exists a natural number $N_3$ such that for all $i>N_3$,
$$d(f(a),f(x_i))<\frac{\epsilon}{3}.$$
Hence, for every positive real number $\epsilon$ there exists a natural number $N=\max\{N_1,N_2,N_3\}$ such that for all $i>N$, by (m-lbnd)
$$0\le d(f(a),g(a))$$
using (m-inq) twice we get
$$\le d(f(a),f(x_i))+d(f(x_i),g(x_i))+d(g(x_i),g(a))$$
by (m-sym)
$$d(f(a),f(x_i))+d(f(x_i),g(x_i))+d(g(a),g(x_i))$$
$$<\frac{\epsilon}{3}+\frac{\epsilon}{3}+\frac{\epsilon}{3}=\epsilon.$$
Therefore, $d(f(a),g(a))=0$ and by (d-sep) $f(a)=g(a)$ and, hence, by \Cref{D7-3} $a$ is a coincidence point of $f$ and  $g. \hspace{6ex}\square$
\section{Partial Metric Space}
\label{C7b}
\begin{thm}{\textbf{(Fixed point and Partial metrics\cite{Ass20151}):}}
\label{T7b1}
\\Let $(X,p)$ be a complete partial metric space with $x_o$ in $X$. Let $f:X \to X$ be a Cauchy function at $x_o$. If $f$ is non-expansive and  weakly orbitally continuous at $x_o$ then $f$ has a fixed point.
\end{thm}
\textbf{Proof:} Since $f$ is Cauchy at $x_o$ then by \Cref{D5-2}, $\{f^i(x_o)\}_{i \in \mathbb{N}}$ is a Cauchy sequence.
Since $(X,p)$ is a complete partial metric space then by \Cref{D4b3} $\{f^i(x_o)\}_{i \in \mathbb{N}}$ has a  special limit $a$ in $X$. 
\\Since $f$ is non-expansive then by \Cref{L6b1} 
$$p(a,f(a))=p(a,a).$$
Since $f$ is weakly orbitally continuous at $x_o$ then by \Cref{L6b2}
$$p(a,f(a))=p(f(a),f(a)).$$
Hence, by (p-sep) $f(a)=a$. Therefore, by \Cref{D7-1} $a$ is a fixed point of $f.\hspace{6ex}\square$
\begin{cor}\label{C7b1}
Let $(X,p)$ be a complete partial metric space with $x_o$ in $X$. Let $f:X \to X$ be a non-expansive weakly orbitally continuous function at $x_o$. If either one of the below criteria holds true: 
\\a) $f$ is orbitally $c_r-$contractive at $x_o$.
\\b) $f$ is orbitally  $\varphi_r-$contractive at $x_o$.
\\then $f$ has a fixed point.
\end{cor}
\textbf{Proof:} From \Cref{L5b2} and \Cref{T5b1}, if $f$ is orbitally $c_r-$contractive or orbitally  $\varphi_r-$contractive at $x_o$ then $f$ is Cauchy at $x_o$.
\begin{thm}{\textbf{(Common fixed point and Partial metrics):}}
\label{T7b2}
\\Let $(X,p)$ be a complete partial metric space with $x_o$ and $y_o$ in $X$. Let  $f:X \to X$ and $g: X \to X$ be two functions that form a Cauchy pair over $(x_o,y_o)$. If $f$ and $g$ are non-expansive and $f$ and $g$ are weakly orbitally continuous at $x_o$ and $y_o$ respectively. then $f$ and $g$ have a common fixed point.
\end{thm}
\textbf{Proof: }Since $f$ and $g$  form a Cauchy pair at $(x_o,y_o)$ then by \Cref{D5-3} $\{f^i(x_o)\}_{i\in\mathbb{N}}$ and $\{g^i(y_o)\}_{i\in\mathbb{N}}$ form a Cauchy pair. By \Cref{L4b5} $\{f^i(x_o)\}_{i\in\mathbb{N}}$ and $\{g^i(y_o)\}_{i\in\mathbb{N}}$ are both Cauchy sequences. Since $(X,p)$ is a complete partial metric space then by \Cref{D4b3}  and \Cref{L4b5} $\{f^i(x_o)\}_{i\in\mathbb{N}}$ and $\{g^i(y_o)\}_{i\in\mathbb{N}}$ both have the same special limit  $a$ in $X$. 
\\Since $f$ and $g$ are both non-expansive then by \Cref{L6b1}
$$p(a,a)=p(a,f(a))\text{ and }p(a,a)=p(a,g(a)).$$
Since $f$ and $g$ are weakly orbitally continuous on $x_o$ and $y_o$ respectively then by \Cref{L6b2}
$$p(f(a),f(a))=p(a,f(a))\text{ and }p(g(a),g(a))=p(a,g(a)).$$
Hence, by (p-sep) $f(a)=a=g(a)$. Therefore, by \Cref{D7-2} $a$ is a common fixed point of $f$ and $g.\hspace{6ex}\square$
\begin{cor}\label{C7b2}
Let $(X,p)$ be a complete partial metric space with $x_o$ and $y_o$ in $X$. Let  $f:X \to X$ and $g: X \to X$ be two  non-expansive  functions with $f$ and $g$  weakly orbitally continuous at $x_o$ and $y_o$ respectively. If either one of the below criteria holds true: 
\\a) $f$ and $g$ are $f-$pairwise $c_r-$contractive over $(x_o,y_o)$.
\\b) $f$ and $g$ are $g-$pairwise $c_r-$contractive over $(x_o,y_o)$.
\\then $f$ and $g$ have a common fixed point.
\end{cor}
\textbf{Proof:} From \Cref{T5b1}, if $f$ and $g$ are $f-$pairwise $c_r-$contractive (similarly $g-$pairwise $c_r-$contractive) over $(x_o,y_o)$ then $f$ and $g$ form a Cauchy pair over $(x_o,y_o)$.
\begin{thm}{\textbf{(Coincidence Point Theorem):}}
\label{T7b3}
\\Let $(X,p)$ be a complete partial metric space and let $(Y,h)$ be a partial metric space. Let $f:X \to Y$ and $g: X \to Y$ be two sequentially continuous and consistent functions on $X$. If $f$ and $g$ are $(f,g)-$mutually $c_r$-contractive   then $f$ and $g$ have a coincidence point.
\end{thm}
\textbf{Proof:} Since $f$ and $g$ are  $(f,g)-$mutually $c_r$-contractive  then they are   $f-$mutually $c_r$-contractive . Hence, by \Cref{T5b2} there exists a Cauchy sequence $\{x_i\}_{i \in \mathbb{N}}$ in $X$. Let $r$ be the central distance of $\{x_i\}_{i \in \mathbb{N}}$ then again by \Cref{T5b2} for all natural numbers $i$,
$$r\le p(x_i,x_i).$$ 
Since $(X,p)$ is complete partial metric space then $\{x_i\}_{i \in \mathbb{N}}$ has a special limit $a$ in $X$. From \Cref{D4b2} for all natural numbers $i$,
$$r=p(a,a)\le p(x_i,x_i).$$
$g$ is consistent then
$$h(g(a),g(a))\le h(g(x_i),g(x_i))$$
and, hence,
$$-h(g(x_i),g(x_i))\le -h(g(a),g(a)).$$
For every positive real number $\epsilon$ we know that:
\\Since $g$ is sequentially continuous, by \Cref{D6-2} $g(a)$ is a limit of  $\{g(x_i)\}_{i \in \mathbb{N}}$. Therefore, there exists a natural number $N_1$ such that for all $i>N_1$,
$$h(g(a),g(x_i))<h(g(a),g(a))+\frac{\epsilon}{3}.$$
Similarly, $f$ is sequentially continuous hence, by \Cref{D6-2} $f(a)$ is limit of $\{f(x_i)\}_{i \in \mathbb{N}}$. Therefore, there exists a natural number $N_2$ such that for all $i>N_2$,
$$h(f(a),f(x_i))<h(f(a),f(a))+\frac{\epsilon}{3}.$$ 
By \Cref{T5b2}, there exists a natural number $N_3$ such that for all $i>N_3$,
$$h(f(x_i),g(x_i))-h(f(x_i),f(x_i))<\frac{\epsilon}{3}.$$
Hence, for every positive real number $\epsilon$ there exists a natural number $N=\max\{N_1,N_2,N_3\}$ such that for all $i>N$,by (p-lbnd)
$$h(f(a),f(a))\le h(f(a),g(a))$$
using (p-inq) twice we get for all $i$
$$\le h(f(a),f(x_i))-h(f(x_i),f(x_i))+h(f(x_i),g(x_i))-h(g(x_i),g(x_i))+h(g(x_i),g(a))$$
by (p-sym)
$$=h(f(a),f(x_i))+h(f(x_i),g(x_i))-h(f(x_i),f(x_i))-h(g(x_i),g(x_i))+h(g(a),g(x_i))$$
$$<h(f(a),f(a))+\frac{\epsilon}{3}+\frac{\epsilon}{3}-h(g(a),g(a))+h(g(a),g(a))+\frac{\epsilon}{3}$$
$$=h(f(a),f(a))+\epsilon.$$
Therefore, $h(f(a),f(a))=h(f(a),g(a)).$ Repeating the above process with $f$ and $g$ being  $g-$mutually $c_r$-contractive  and $f$ being consistent we get $h(g(a),g(a))= h(f(a),g(a))$ and, hence, by (p-sep) $f(a)=g(a)$. Therefore, by \Cref{D7-3} $a$ is a coincidence point of $f$ and $g.\hspace{6ex}\square$
\section{Strong Partial Metric Space}
\label{C7c}
\indent \indent We remind our reader that (s-lbnd) is a stronger version of (p-sep). Hence, if $(X,s)$ is a strong partial metric space, for any two element $x$ and $z$ in $X$, it is enough to have $s(x,z)\le s(x,x)$ to deduce that $x=z$. Therefore, we are able to relax the requirements on the functions studied to assert the existence of the fixed point, common fixed point or coincidence point in question. 
\begin{thm}{\textbf{(Fixed point and Non-expansive \cite{Ass20151}):}}
\label{T7c1}
\\Let $(X,s)$ be a complete strong partial metric space with $x_o$ in $X$. Let $f:X \to X$ be a Cauchy function at $x_o$. If $f$ is non-expansive  then $f$ has a fixed point.
\end{thm}
\textbf{Proof:} Since $f$ is Cauchy at $x_o$ then by \Cref{D5-2}, $\{f^i(x_o)\}_{i \in \mathbb{N}}$ is a Cauchy sequence.
Since $(X,s)$ is a complete strong partial metric space then by \Cref{D4b3} $\{f^i(x_o)\}_{i \in \mathbb{N}}$ has a  special limit $a$ in $X$. 
\\Since $f$ is non-expansive then by \Cref{L6b1} 
$$s(a,f(a))=s(a,a)$$
and, hence by (s-lbnd) $f(a)=a$. Therefore, by \Cref{D7-1} $a$ is a fixed point of $f.\hspace{6ex}\square$
\begin{thm}{\textbf{(Fixed point and Weak orbital continuity \cite{Ass20151}):}}
\label{T7c2}
\\Let $(X,s)$ be a complete strong partial metric space with $x_o$ in $X$. Let $f:X \to X$ be a Cauchy function at $x_o$. If $f$ is  weakly orbitally continuous at $x_o$ then $f$ has a fixed point.
\end{thm}
\textbf{Proof:} Since $f$ is Cauchy at $x_o$ then by \Cref{D5-2}, $\{f^i(x_o)\}_{i \in \mathbb{N}}$ is a Cauchy sequence.
Since $(X,p)$ is a complete strong partial metric space then by \Cref{D4b3} $\{f^i(x_o)\}_{i \in \mathbb{N}}$ has a  special limit $a$ in $X$. 
Since $f$ is weakly orbitally continuous at $x_o$ then by \Cref{L6b2}
$$s(a,f(a))=s(f(a),f(a))$$
and, hence, by (s-lbnd) $f(a)=a$. Therefore, by \Cref{D7-1} $a$ is a fixed point of $f.\hspace{6ex}\square$
\begin{cor}\label{C7c1}
Let $(X,s)$ be a complete strong partial metric space with $x_o$ in $X$. Let $f:X \to X$ be a non-expansive function or a weakly orbitally continuous function at $x_o$. If either one of the below criteria holds true: 
\\a) $f$ is orbitally $c_r-$contractive at $x_o$.
\\b) $f$ is orbitally  $\varphi_r-$contractive at $x_o$.
\\then $f$ has a fixed point.
\end{cor}
\textbf{Proof:} From \Cref{L5b2}  and \Cref{L5b1}, if $f$ is orbitally $c_r-$contractive or orbitally $\varphi_r-$contractive at $x_o$ then $f$ is Cauchy at $x_o$.
\begin{thm}{\textbf{(Common fixed point and Non-expansive):}}
\label{T7c3}
\\Let $(X,s)$ be a complete strong partial metric space with $x_o$ and $y_o$ in $X$. Let  $f:X \to X$ and $g: X \to X$ be two functions that form a Cauchy pair over $(x_o,y_o)$. If $f$ and $g$ are non-expansive  then $f$ and $g$ have a common fixed point.
\end{thm}
\textbf{Proof:} Since $f$ and $g$  form a Cauchy pair at $(x_o,y_o)$ then by \Cref{D5-3} $\{f^i(x_o)\}_{i\in\mathbb{N}}$ and $\{g^i(y_o)\}_{i\in\mathbb{N}}$ form a Cauchy pair. By \Cref{L4b5} $\{f^i(x_o)\}_{i\in\mathbb{N}}$ and $\{g^i(y_o)\}_{i\in\mathbb{N}}$ are both Cauchy sequences. Since $(X,s)$ is a complete strong partial metric space then by \Cref{D4b3}  and \Cref{L4b5} $\{f^i(x_o)\}_{i\in\mathbb{N}}$ and $\{g^i(y_o)\}_{i\in\mathbb{N}}$ both have the same special limit  $a$ in $X$. Since $f$ and $g$ are both non-expansive then by \Cref{L6b2}
$$p(a,a)=p(a,f(a))\text{ and }p(a,a)=p(a,g(a))$$
and, hence, by (s-lbnd) $f(a)=a=g(a)$. Therefore, by \Cref{D7-2} $a$ is a common fixed point of $f$ and $g.\hspace{6ex}\square$
\begin{thm}{\textbf{(Common fixed point and Weak orbital continuity):}}
\label{T7c4}
\\Let $(X,s)$ be a complete strong partial metric space with $x_o$ and $y_o$ in $X$. Let  $f:X \to X$ and $g: X \to X$ be two functions that form a Cauchy pair over $(x_o,y_o)$. If $f$ and $g$ are weakly orbitally continuous on $x_o$ and $y_o$ respectively then $f$ and $g$ have a common fixed point.
\end{thm}
\textbf{Proof:} Since $f$ and $g$  form a Cauchy pair at $(x_o,y_o)$ then by \Cref{D5-3} $\{f^i(x_o)\}_{i\in\mathbb{N}}$ and $\{g^i(y_o)\}_{i\in\mathbb{N}}$ form a Cauchy pair. By \Cref{L4b5} $\{f^i(x_o)\}_{i\in\mathbb{N}}$ and $\{g^i(y_o)\}_{i\in\mathbb{N}}$ are both Cauchy sequences. Since $(X,s)$ is a complete strong partial metric space then by \Cref{D4b3} and \Cref{L4b5} $\{f^i(x_o)\}_{i\in\mathbb{N}}$ and $\{g^i(y_o)\}_{i\in\mathbb{N}}$ both have the same special limit  $a$ in $X$. Since $f$ and $g$ are weakly orbitally continuous on $x_o$ and $y_o$ respectively then 
$$p(f(a),f(a))=p(a,f(a))\text{ and }p(g(a),g(a))=p(a,g(a))$$
and, hence, by (s-lbnd) $f(a)=a=g(a)$. Therefore, by \Cref{D7-2} $a$ is a common fixed point of $f$ and $g.\hspace{6ex}\square$
\begin{thm}{\textbf{(Common fixed point and Mixed criteria):}}
\label{T7c5}
\\Let $(X,s)$ be a complete strong partial metric space with $x_o$ and $y_o$ in $X$. Let  $f:X \to X$ and $g: X \to X$ be two functions that form a Cauchy pair over $(x_o,y_o)$. If $f$ is non-expansive and $g$ is weakly orbitally continuous on $y_o$  then $f$ and $g$ have a common fixed point.
\end{thm}
\textbf{Proof:} Since $f$ and $g$  form a Cauchy pair at $(x_o,y_o)$ then by \Cref{D5a4} $\{f^i(x_o)\}_{i\in\mathbb{N}}$ and $\{g^i(y_o)\}_{i\in\mathbb{N}}$ form a Cauchy pair. By \Cref{L4b5} $\{f^i(x_o)\}_{i\in\mathbb{N}}$ and $\{g^i(y_o)\}_{i\in\mathbb{N}}$ are both Cauchy sequences. Since $(X,p)$ is a complete partial metric space then by \Cref{D4b3} and \Cref{L4b5} $\{f^i(x_o)\}_{i\in\mathbb{N}}$ and $\{g^i(y_o)\}_{i\in\mathbb{N}}$ both have the same special limit  $a$ in $X$. 
\\Since $f$ is non-expansive then by \Cref{L6b1}
$$p(a,a)=p(a,f(a))$$
and, hence by (s-lbnd)
$$f(a)=a.$$
Since  $g$ is weakly orbitally continuous on  $y_o$ then by \Cref{L6b2}
$$p(g(a),g(a))=p(a,g(a))$$
and, hence, by (s-lbnd) $f(a)=a=g(a)$. Therefore, by \Cref{D7-2} $a$ is a common fixed point of $f$ and $g.\hspace{6ex}\square$
\begin{cor}\label{C7c2}
Let $(X,s)$ be a complete strong partial metric space with $x_o$ and $y_o$ in $X$. Let  $f:X \to X$ be a non-expansive  function or  a weakly orbitally continuous function at $x_o$. Similarly, let  $g: X \to X$ be a non-expansive  function or a weakly orbitally continuous function at  $y_o$. If either one of the below criteria holds true: 
\\a) $f$ and $g$ are $f-$pairwise $c_r-$contractive over $(x_o,y_o)$.
\\b) $f$ and $g$ are $g-$pairwise $c_r-$contractive over $(x_o,y_o)$.
\\then $f$ and $g$ have a common fixed point.
\end{cor}
\textbf{Proof:} From \Cref{T5b1} , if $f$ and $g$ are $f-$pairwise $c_r-$contractive (similarly $g-$pairwise $c_r-$contractive) over $(x_o,y_o)$ then $f$ and $g$ form a Cauchy pair over $(x_o,y_o)$.
\begin{thm}{\textbf{(Coincidence Point Theorem):}}
\label{T7c6}
\\Let $(X,p)$ be a complete  partial metric space and let $(Y,s)$ be a strong partial metric space. Let $f:X \to Y$ and $g: X \to Y$ two sequentially continuous functions on $X$. If $f$ and $g$ are $f-$mutually $c_r$-contractive and $g$ is consistent   then $f$ and $g$ have a coincidence point.
\end{thm}
\textbf{Proof:} Since $f$ and $g$ are $f-$mutually $c_r$-contractive then by \Cref{T5b2} there exists a Cauchy sequence $\{x_i\}_{i \in \mathbb{N}}$ in $X$. Let $r$ be the central distance of $\{x_i\}_{i \in \mathbb{N}}$ then again by \Cref{T5b2} for all natural numbers $i$,
$$r\le p(x_i,x_i).$$ 
Since $(X,p)$ is complete partial metric space then $\{x_i\}_{i \in \mathbb{N}}$ has a special limit $a$ in $X$. From \Cref{D4b2} for all natural numbers $i$,
$$r=p(a,a)\le p(x_i,x_i).$$
$g$ is consistent then
$$s(g(a),g(a))\le s(g(x_i),g(x_i))$$
and, hence,
$$-s(g(x_i),g(x_i))\le -s(g(a),g(a)).$$
For every positive real number $\epsilon$ we know that:\\ since $g$ is sequentially continuous, by \Cref{D6-2} $g(a)$ is limit of  $\{g(x_i)\}_{i \in \mathbb{N}}$. Therefore, there exists a natural number $N_1$ such that for all $i>N_1$,
$$s(g(a),g(x_i))<s(g(a),g(a))+\frac{\epsilon}{3}.$$
Similarly, $f$ is sequentially continuous hence, by \Cref{D6-2} $f(a)$ is limit of $\{f(x_i)\}_{i \in \mathbb{N}}$. Therefore, there exists a natural number $N_2$ such that for all $i>N_2$,
$$s(f(a),f(x_i))<s(f(a),f(a))+\frac{\epsilon}{3}.$$ 
By \Cref{T5b2}, there exists a natural number $N_3$ such that for all $i>N_3$,
$$s(f(x_i),g(x_i))-s(f(x_i),f(x_i))<\frac{\epsilon}{3}.$$
Hence, for every positive real number $\epsilon$ there exists a natural number $N=\max\{N_1,N_2,N_3\}$ such that for all $i>N$, using (s-inq)
twice$$ h(f(a),g(a))\le s(f(a),f(x_i))-s(f(x_i),f(x_i))+s(f(x_i),g(x_i))-s(g(x_i),g(x_i))+s(g(x_i),g(a))$$
by (s-sym)
$$=s(f(a),f(x_i))+s(f(x_i),g(x_i))-s(f(x_i),f(x_i))-s(g(x_i),g(x_i))+s(g(a),g(x_i))$$
$$<h(f(a),f(a))+\frac{\epsilon}{3}+\frac{\epsilon}{3}-h(g(a),g(a))+h(g(a),g(a))+\frac{\epsilon}{3}$$
$$=s(f(a),f(a))+\epsilon.$$
Therefore, $s(f(a),g(a))\le s(f(a),f(a)).$ Hence, by (s-lbnd) $f(a)=g(a)$. Therefore, by \Cref{D7-3} $a$ is a coincidence point of $f$ and $g.\hspace{6ex}\square$
\section{Partial $n-\mathfrak{M}$etric Space}
\label{C7d}
\begin{thm}{\textbf{(Fixed point and Partial $n-\mathfrak{M}$etrics\cite{Ass20152}):}}
\label{T7d1}
\\Let $(X,P)$ be a complete partial $n-\mathfrak{M}$etric space with $x_o$ in $X$. Let $f:X \to X$ be a Cauchy function at $x_o$. If $f$ is non-expansive and  weakly orbitally continuous at $x_o$ then $f$ has a fixed point.
\end{thm}
\begin{cor}\label{C7d1}
Let $(X,P)$ be a complete  partial $n-\mathfrak{M}$etric space with $x_o$ in $X$. Let $f:X \to X$ be a non-expansive weakly orbitally continuous function at $x_o$. If either one of the below criteria holds true: 
\\a) $f$ is orbitally $c_r-$contractive at $x_o$.
\\b) $f$ is orbitally  $\varphi_r-$contractive at $x_o$.
\\then $f$ has a fixed point.
\end{cor}
\begin{thm}{\textbf{(Common fixed point and Partial $n-\mathfrak{M}$etrics):}}
\label{T7d2}
\\Let $(X,P)$ be a complete partial $n-\mathfrak{M}$etric  space with $x_o$ and $y_o$ in $X$. Let  $f:X \to X$ and $g: X \to X$ be two functions that form a Cauchy pair over $(x_o,y_o)$. If $f$ and $g$ are non-expansive and $f$ and $g$ are weakly orbitally continuous on $x_o$ and $y_o$ respectively then $f$ and $g$ have a common fixed point.
\end{thm}
\begin{cor}\label{C7d2}
Let $(X,P)$ be a complete partial $n-\mathfrak{M}$etric  space with $x_o$ and $y_o$ in $X$. Let  $f:X \to X$ and $g: X \to X$ be two  non-expansive  functions with $f$ and $g$  weakly orbitally continuous at $x_o$ and $y_o$ respectively. If either one of the below criteria holds true: 
\\a) $f$ and $g$ are $f-$pairwise $c_r-$contractive over $(x_o,y_o)$.
\\b) $f$ and $g$ are $g-$pairwise $c_r-$contractive over $(x_o,y_o)$.
\\then $f$ and $g$ have a common fixed point.
\end{cor}
\begin{thm}{\textbf{(Coincidence Point Theorem):}}
\label{T7d3}
\\Let $(X,P)$ be a complete partial $n-\mathfrak{M}$etric   space and let $(Y,H)$ be a partial $n-\mathfrak{M}$etric   space. Let $f:X \to Y$ and $g: X \to Y$ be two sequentially continuous and consistent functions on $X$. If $f$ and $g$ are $(f,g)-$mutually $c_r$-contractive then $f$ and $g$ have a coincidence point.
\end{thm}
\section{Strong Partial $n-\mathfrak{M}$etric Space}
\label{C7e}
\indent \indent In this section again, we remind our reader that ($S_n$-lbnd) is a stronger version of ($P_n$-sep). Hence, if $(X,S)$ is a strong partial  $n-\mathfrak{M}$etric space, for any two element $x$ and $z$ in $X$, it is enough to have $S(\langle x\rangle^{n-1},z)\le S(\langle x \rangle^n)$ to deduce that $x=z$. Therefore, we are able to relax the requirements on the functions studied in \Cref{C7d} to assert the existence of the fixed point, common fixed point or coincidence point in question. 
\begin{thm}{\textbf{(Fixed point and Non-expansive \cite{Ass20152}):}}
\label{T7e1}
\\Let $(X,S)$ be a complete strong partial   $n-\mathfrak{M}$etric space with $x_o$ in $X$. Let $f:X \to X$ be a Cauchy function at $x_o$. If $f$ is non-expansive  then $f$ has a fixed point.
\end{thm}
\begin{thm}{\textbf{(Fixed point and Weak orbital continuity \cite{Ass20152}):}}
\label{T7e2}
\\Let $(X,S)$ be a complete strong partial   $n-\mathfrak{M}$etric space with $x_o$ in $X$. Let $f:X \to X$ be a Cauchy function at $x_o$. If $f$ is  weakly orbitally continuous at $x_o$ then $f$ has a fixed point.
\end{thm}
\begin{cor}\label{C7e1}
Let $(X,S)$ be a complete strong partial   $n-\mathfrak{M}$etric space with $x_o$ in $X$. Let $f:X \to X$ be a non-expansive function or a weakly orbitally continuous function at $x_o$. If either one of the below criteria holds true: 
\\a) $f$ is orbitally $c_r-$contractive at $x_o$.
\\b) $f$ is orbitally  $\varphi_r-$contractive at $x_o$.
\\then $f$ has a fixed point.
\end{cor}
\begin{thm}{\textbf{(Common fixed point and Non-expansive):}}
\label{T7e3}
\\Let $(X,S)$ be a complete strong partial   $n-\mathfrak{M}$etric space with $x_o$ and $y_o$ in $X$. Let  $f:X \to X$ and $g: X \to X$ be a Cauchy pair over $(x_o,y_o)$. If $f$ and $g$ are non-expansive  then $f$ and $g$ have a common fixed point.
\end{thm}
\begin{thm}{\textbf{(Common fixed point and Weak orbital continuity):}}
\label{T7e4}
\\Let $(X,S)$ be a complete strong partial   $n-\mathfrak{M}$etric space with $x_o$ and $y_o$ in $X$. Let  $f:X \to X$ and $g: X \to X$ be two functions that form a Cauchy pair over $(x_o,y_o)$. If $f$ and $g$ are weakly orbitally continuous on $x_o$ and $y_o$ respectively then $f$ and $g$ have a common fixed point.
\end{thm}
\begin{thm}{\textbf{(Common fixed point and Mixed criteria):}}
\label{T7e5}
\\Let $(X,S)$ be a complete strong partial   $n-\mathfrak{M}$etric space with $x_o$ and $y_o$ in $X$. Let  $f:X \to X$ and $g: X \to X$ be two functions that form a Cauchy pair over $(x_o,y_o)$. If $f$ is non-expansive and $g$ is weakly orbitally continuous on $y_o$  then $f$ and $g$ have a common fixed point.
\end{thm}
\begin{cor}\label{C7e2}
Let $(X,S)$ be a complete strong partial   $n-\mathfrak{M}$etric space with $x_o$ and $y_o$ in $X$. Let  $f:X \to X$ be a non-expansive  function or  a weakly orbitally continuous function at $x_o$. Similarly, let  $g: X \to X$ be a non-expansive  function or a weakly orbitally continuous function at  $y_o$. If either one of the below criteria holds true: 
\\a) $f$ and $g$ are $f-$pairwise $c_r-$contractive over $(x_o,y_o)$.
\\b) $f$ and $g$ are $g-$pairwise $c_r-$contractive over $(x_o,y_o)$.
\\then $f$ and $g$ have a common fixed point.
\end{cor}
\newpage
\begin{thm}{\textbf{(Coincidence Point Theorem):}}
\label{T7e6}
\\Let $(X,P)$ be a complete  partial   $n-\mathfrak{M}$etric space and let $(Y,S)$ be a strong partial   $n-\mathfrak{M}$etric space. Let $f:X \to Y$ and $g: X \to Y$ be two sequentially continuous functions on $X$. If $f$ and $g$ are $f-$mutually $c_r$-contractive and $g$ is consistent   then $f$ and $g$ have a coincidence point.
\end{thm}

\bibliographystyle{amsplain}
\bibliography{Generalizedmetrics}

\end{document}